\documentclass[11pt]{article}

\usepackage{amsfonts}
\usepackage{graphics}
\usepackage{amsmath}
\usepackage{amscd}
\usepackage{amssymb}
\usepackage{rlepsf}
\usepackage{makeidx}

\usepackage{euscript}

\newtheorem{Theorem}{Theorem}

\newtheorem{Definition}[Theorem]{Definition}

\newtheorem{Problem}{Problem}

\newtheorem{Exercise}{Exercise}

\newtheorem{Fundamental Theorem}{Fundamental Theorem}

\newenvironment{proof}[1][Proof]{\textbf{#1.} }{\ \rule{0.5em}{0.5em}}

\def \A {\mathcal{A}}

\def \d {\partial}
\def \f {\phi}
\def \G {\mathcal{G}}
\def \H {\EuScript{H}}

\def \Hom {\mathrm{Hom}}
\def \k {\kappa}

\def \Q {\mathbb{Q}}
\def \R {\mathbb{R}}
\def \ra {\xrightarrow}

\def \s {\scriptstyle}
\def \ss {\scriptscriptstyle}
\def \S {\Sigma}
\def \t {\triangleright}

\def \V {\EuScript{V}}
\def \Z {\mathbb{Z}}

\def \Hom {\mathrm{Hom}}
\def \Alex {\mathrm{Alex}}
\def \CM {\mathrm{CM}}
\def \k {\kappa}

\begin{document}

\title{Invariants of Welded Virtual Knots Via Crossed Module Invariants of Knotted Surfaces}

\author{ Louis H. Kauffman \\ \footnotesize \it {Department of Mathematics, Statistics, and Computer Science,} \\ \footnotesize \it {
University of Illinois at Chicago,}\\ {\footnotesize 851 South Morgan St., Chicago, IL 60607-7045, USA} \\ \footnotesize \it {kauffman@uic.edu} \\ \\ Jo\~{a}o  Faria Martins\footnote{{Also at Departamento de Matem\'{a}tica, Universidade Lus\'{o}fona de Humanidades e Tecnologia, Av. do Campo Grande, 376, 1749-024, Lisboa, Portugal.}}\\ \footnotesize\it  {Departamento de Matem\'{a}tica},\\ \footnotesize\it {Instituto Superior T\'{e}cnico (Universidade T\'{e}cnica de Lisboa)}\\ {\footnotesize\it Av. Rovisco Pais, 1049-001 Lisboa, Portugal}\\  {\footnotesize\it jmartins@math.ist.utl.pt}}

\maketitle

\begin{abstract}
We define an invariant of welded virtual knots  from each finite crossed
module by considering crossed module invariants of ribbon knotted surfaces
which are naturally associated with them. We elucidate that the invariants
obtained are {non-trivial} by calculating explicit  examples. {We define welded virtual graphs and consider invariants of them defined in a similar way}.

\end{abstract}

\noindent {\bf  2000 Mathematics Subject Classification:} {\it  {57M25 
 (primary), 57Q45 
(secondary).}}

\section{Introduction}
 Welded virtual knots were defined in \cite{K1},  by allowing one extra move
 in addition {to} the moves appearing in the definition of a virtual knot. This
 extra move preserves the (combinatorial) fundamental group of {the} complement,
 which is therefore an invariant of welded virtual knots (the knot group).  Given a finite
 group $G$, one can therefore define a welded virtual  knot invariant $\H_G$,
 by considering the number of morphisms from the fundamental group of the 
 complement into $G$. The Wirtinger presentation of knot groups enables a
 quandle type calculation of this ``Counting Invariant'' $\H_G$.

Not a lot of welded virtual knot invariants are known. The aim of this
 article is to introduce a new one, the ``Crossed Module Invariant'' $\H_\G$, which depends on a finite automorphic crossed module
 $\G=(E,G,\t)${, in other words} on a pair of groups $E$ and $G$, with $E$
 abelian, and a  left action of $G$ on $E$ by automorphisms.

The Crossed Module Invariant $\H_\G$ reduces to the Counting Invariant $\H_G$ when $E=0$. However, the Crossed Module Invariant distinguishes, in some
 cases, between welded virtual links  with the same knot group, and therefore it is  strictly stronger than the Counting Invariant.  {We will assert this fact by calculating explicit examples.}

 {Let $\G=(E,G,\t)$ be an automorphic crossed module.} Note that the
 Counting Invariant ${\H_G}$ is trivial whenever $G$ is abelian. However, taking
 $G$ to be abelian and 
 $E$ to be {non-trivial,} 
 yields a {non-trivial} invariant $\H_\G$, which is, as a rule, much easier to calculate than
 the Counting Invariant $\H_G$ where  $G$ is  generic group, and it is strong enough to tell apart some pairs of links with the same knot group.
{Suppose that the welded virtual link $K$ has $n$-components. Let $\k_n=\Z[X_1,X_1^{-1},\ldots , X_n,X_n^{-1}]$. We will define a $k_n$-module $\CM(K)$, depending only on $K$, up  to isomorphism and permutations of the variables $X_1,\ldots, X_n$. If  $G$ is abelian, then $\H_\G$ simply counts the number of crossed module morphisms $\CM(K) \to \G$. We prove in this article that if $K$ is classical then $\CM(K)$ coincides with the Alexander module $\Alex(K)$ {of $K$}. However, this is not the case if $K$ is not classical. We will give  examples of   pairs of welded virtual links $(K,K')$ with the same knot group (thus the same Alexander module) but with {$\CM(K) \ncong \CM(K')$}. This will happen when $K$ and $K'$ have the same knot group, but are distinguished by their crossed module invariants for $G$ abelian.}

Let us explain the construction of the Crossed Module Invariant $\H_\G$. Extending a previous construction due to T. Yagima,  Shin Satoh defined in \cite{S} a map which associates an oriented
knotted torus $T(K)$, the ``tube of {$K$}'', to each oriented welded virtual knot
$K$. The map $K \mapsto
T(K)$ preserves knot groups. 
In the case when $K$
is a classical knot, then $T(K)$ coincides with the torus spun of $K${, obtained by spinning $K$ $4$-dimensionally, in order to obtain an embedding of the {torus $S^1 \times S^1$ into} $S^4$.}

 The existence of the {tube} map $K\mapsto T(K)$  makes it natural to define
 invariants of welded virtual  knots by considering invariants of knotted surfaces.
 We will consider this construction for the case of the crossed module
 invariants $I_\G(\S)$ of knotted surfaces $\S$, defined in \cite{FM1,FM2}.  Here $\G=\left (E \ra{\d} G,\t\right) $ is a finite
 crossed module. {Note that the invariant $I_\G$  on a knotted surface coincides with Yetter's Invariant (see \cite{Y3,P1,FMP}) of the complement of it.}   We can  thus define a welded
 virtual knot invariant by considering $\H_\G(K)\doteq I_\G(T(K))$, where $K$ is a {welded virtual} knot. 

A straightforward analysis of the crossed module invariant of the {tube  $T(K)$ of  the welded virtual knot  $K$} permits
the evaluation of $ \H_\G(K)$ in a quandle type way, albeit the biquandle we define
is sensitive to maximal and minimal points, so it should probably be called a
``Morse biquandle''. 

  A proof of the existence of the invariant $\H_\G$, where $\G$ is {a {finite crossed module,} can be done directly,} from the Morse biquandle obtained. In fact all the results of this article are {fully} independent of the 4-dimensional picture, {and can be given a direct proof. Moreover, they confirm} the results obtained previously for the crossed module invariants $I_\G$  of knotted surfaces in $S^4$.

As we have referred to {above}, the {tube} map $K \mapsto T(K)$ preserves the fundamental group of the
complements. We prove that $\H_\G$ is powerful enough to distinguish between
distinct welded virtual links with the same knot group.  {For example, we will construct an infinite  set of pairs $(P_i,{c_1(P_i')})$, where $i$ is an odd integer,  of welded virtual links with the following properties:}
\begin{enumerate}
\item {$P_i$ and ${c_1(P_i')}$ each have two components for all $i$.}
\item {$P_i$ and ${c_1(P_i')}$ have isomorphic knot groups for each $i$.}
\item {$P_i$ and ${c_1(P_i')}$ can be  distinguished by their crossed module invariant for each $i$.}
\end{enumerate}
{In fact $P_i$ and ${c_1(P_i')}$ will be distinguished by their crossed module invariant $\H_\G$ with $\G=(E,G,\t)$ being an automorphic crossed module with $G$ abelian.}
 This in particular proves that the Crossed Module Invariant of knotted surfaces $I_\G$ defined in \cite{FM1,FM2,FM3}  sees beyond  the fundamental group of their complement, {in an infinite number of cases.}

{We will also give examples of pairs of 1-component welded virtual knots with the same knot group, but separated by their crossed module invariants. However,  we will need to  make use of computer based calculations in this case.}

In this article we will  propose a definition of Welded Virtual Graphs. {The Crossed Module Invariant of welded virtual links extends naturally to them.}

\section{An Invariant of Welded Virtual Knots}

\subsection{Welded virtual knots}\label{r1}
\begin{figure}
\centerline{\relabelbox 
\epsfysize 1.5cm
\epsfbox {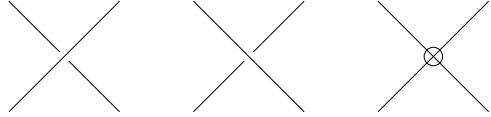}
\endrelabelbox }
\caption{\label{Cross} Classical and virtual crossings.  }
\end{figure}

\begin{figure}
\centerline{\relabelbox 
\epsfysize 2.5cm
\epsfbox {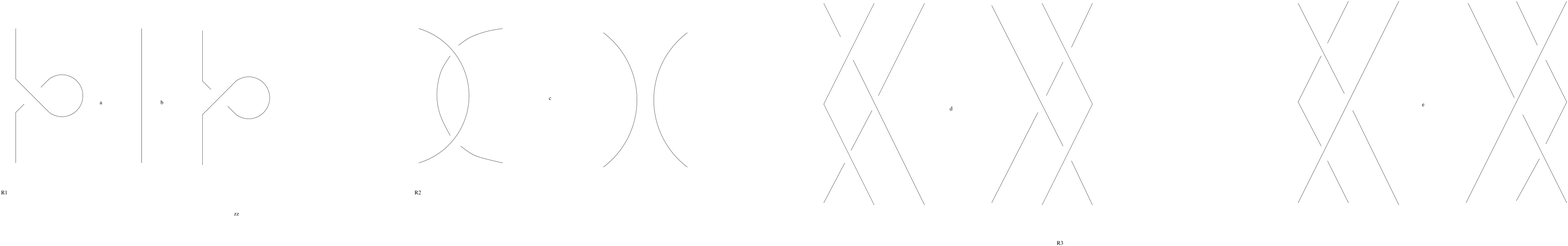}
\relabel {a}{$\s{\leftrightarrow}$}
\relabel {b}{$\s{\leftrightarrow}$}
\relabel {c}{$\s{\leftrightarrow}$}
\relabel {d}{$\s{\leftrightarrow}$}
\relabel {e}{$\s{\leftrightarrow}$}
\relabel {R1}{$\s{\textrm{Reidemeister-I Move}}$}
\relabel {R2}{$\s{\textrm{Reidemeister-II Move}}$}
\relabel {R3}{$\s{\textrm{Reidemeister-III Move}}$}
\endrelabelbox }
\caption{\label{Reidemeister} Reidemeister Moves I, II and III.}
\end{figure}

\begin{figure}
\centerline{\relabelbox 
\epsfysize 1.4cm
\epsfbox {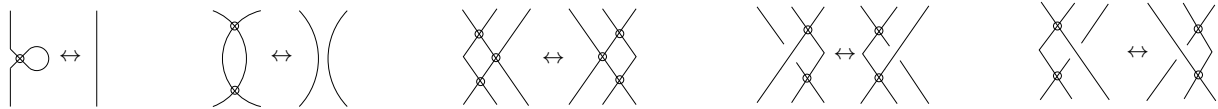}
\relabel {a}{$\s{\leftrightarrow}$}
\relabel {b}{$\s{\leftrightarrow}$}
\relabel {c}{$\s{\leftrightarrow}$}
\relabel {d}{$\s{\leftrightarrow}$}
\relabel {e}{$\s{\leftrightarrow}$}
\endrelabelbox }
\caption{\label{VReidemeister} Virtual Reidemeister Moves.}
\end{figure}

Recall that a virtual knot diagram is, by definition, an immersion of a disjoint
union of circles into the plane $\R^2$, where the 4-valent vertices of the immersion can represent either classical or virtual crossing{; see} figure \ref{Cross}. The definition of an oriented virtual knot diagram is the obvious one.
We say that two virtual knot diagrams are equivalent if they can be related by
the moves of figures \ref{Reidemeister} and \ref{VReidemeister}, as well as
planar isotopy.  It is
important to note that in the oriented case we will need to consider all the
possible orientations of the strands. A virtual knot is an equivalence class
of virtual knot diagrams under the equivalence relation just described{; see} \cite{K1}.

\begin{figure}
\centerline{\relabelbox 
\epsfysize 1.5cm
\epsfbox {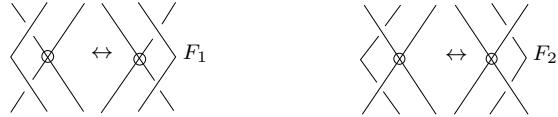}
\relabel {a}{$\s{\leftrightarrow}$}
\relabel {b}{$\s{\leftrightarrow}$}
\relabel {B}{$\s{F_1}$}
\relabel {A}{$\s{F_2}$}
\endrelabelbox }
\caption{\label{Forbiden} The forbidden moves $F_1$ and $F_2$.}
\end{figure}

Observe that, as far as virtual knots are concerned,  we do not allow the moves
shown in figure \ref{Forbiden}, called respectively the forbidden moves $F_1$
and $F_2$. Considering the first forbidden move $F_1$ in
addition  to the ones appearing in the definition of a virtual knot, one
obtains the notion of a ``welded virtual knot'', due {to} the {first} author{; see} {\cite{K1}. }

\subsubsection{The fundamental group of the complement}
\begin{figure}
\centerline{\relabelbox
\epsfysize 1.8cm
\epsfbox{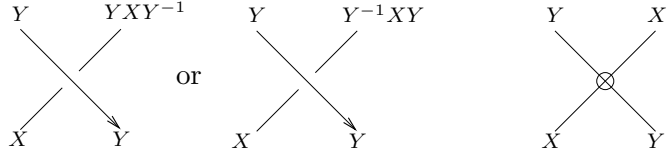}
\relabel {A}{$\s{Y}$}
\relabel {B}{$\s{Y^{-1}XY}$}
\relabel {C}{$\s{X}$}
\relabel {D}{$\s{Y}$}
\relabel {A1}{$\s{Y}$}
\relabel {B1}{$\s{YXY^{-1}}$}
\relabel {C1}{$\s{X}$}
\relabel {D1}{$\s{Y}$}
\relabel {E}{$\s{Y}$}
\relabel {F}{$\s{X}$}
\relabel {G}{$\s{X}$}
\relabel {H}{$\s{Y}$}
\relabel {or}{$\textrm{or}$}
\endrelabelbox }
\caption{\label{Wirtinger} Wirtinger Relations. {The first two are called ``Left Handed'' and ``Right Handed'' Wirtinger  Relations, respectively.} }
\end{figure}

The (combinatorial) fundamental group of the complement of a
virtual knot diagram (the knot group) is, by definition, generated by all the arcs of a diagram of it, considering the relations  {(called Wirtinger Relations)} of figure \ref{Wirtinger} at each
crossing.  {It is understood that in each calculation of a knot group from a virtual knot diagram  we will use either the ``Left Handed'' or the ``Right Handed'' Wirtinger Relation.} {The final result will not depend on this choice.}

In the case of classical {knots or links,} this does coincide with the
fundamental group of the complement, so we can drop the prefix ``combinatorial''. This combinatorial fundamental group is in fact an invariant of welded virtual knots. {This can be proved easily.}

\subsection{Virtual knot presentations of knotted surfaces}\label{presentation}

By definition, a  torus link\footnote{Not to be confused with the 3-dimensional notion of a torus link.} in $S^4$  is an embedding of a disjoint union of
tori $S^1 \times S^1$ into $S^4$, considered up to ambient isotopy. A knotted torus is an embedding of a torus $S^1 \times S^1$ into $S^4$, considered up to ambient isotopy. The definition of an oriented knotted torus or torus link is the obvious one.

As proved in \cite{S,Ya,CKS}, it is possible to  associate  an oriented
torus link  $T(K)\subset S^4$, the ``tube of $K$'',  to each oriented welded virtual {link} $K$. This correspondence was defined first in \cite{Ya}, for the
case of classical knots. The extension to welded virtual knots was completed
in  \cite{S}.

\begin{figure}
\centerline{\relabelbox 
\epsfysize 3cm
\epsfbox {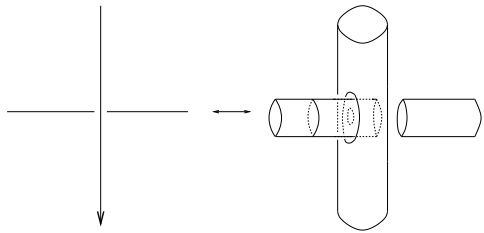}
\endrelabelbox }
\caption{\label{ShinTube} {The tube of a virtual knot at the vicinity of a classical crossing.}}
\end{figure}

\begin{figure}
\centerline{\relabelbox 
\epsfysize 3cm
\epsfbox {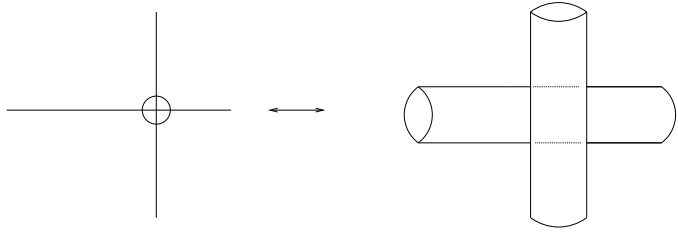}
\endrelabelbox }
\caption{\label{ShinMap2} {The tube of a virtual knot at the vicinity of a  virtual  crossing.}}
\end{figure}

 {The tube map is very easy to define. {Given a virtual {link} diagram, we define the tube  of it by  considering the broken surface diagram obtained by doing the transition of figures \ref{ShinTube} and \ref{ShinMap2}.  For the representation of knotted surfaces in $S^4$ in the form of broken surface diagrams, we refer the reader to \cite{CKS}.} {The tube  of a} virtual knot diagram  has a natural orientation determined by the orientation of a ball in $S^3$. It is proved in \cite{S} that if $K$ and $L$ are diagrams of the same welded virtual knot {then} it follows that $T(K)$ and $T(L)$ are isotopic knotted surfaces in $S^4$. This defines the tube of a welded virtual knot.}

\begin{figure}
\centerline{\relabelbox 
\epsfysize 7cm
\epsfbox {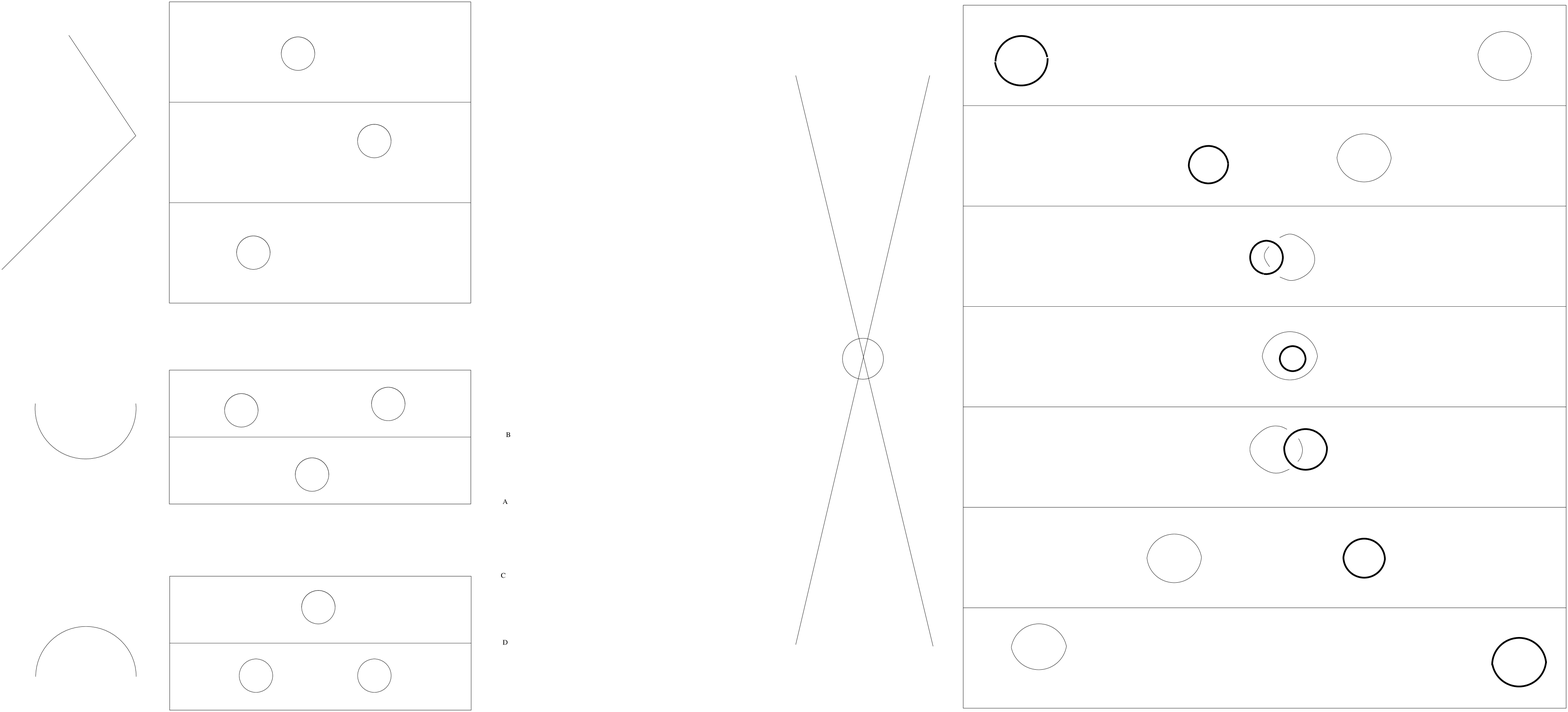}
\relabel {C}{{$\s{\textrm{{death of a circle}}}$}}
\relabel {D}{{$\s{\textrm{{saddle point}}}$}}
\relabel {B}{{$\s{\textrm{{saddle point}}}$}}
\relabel {A}{{$\s{\textrm{{birth of a circle}}}$}}
\endrelabelbox }
\caption{ Associating a knotted torus to a virtual knot: edges, minimal and maximal points and virtual crossings. All circles are oriented counterclockwise. {Note that the movies should be read from bottom to top.}}
\label{Tor1}
\end{figure}

\begin{figure}
\centerline{\relabelbox 
\epsfysize 7cm
\epsfbox {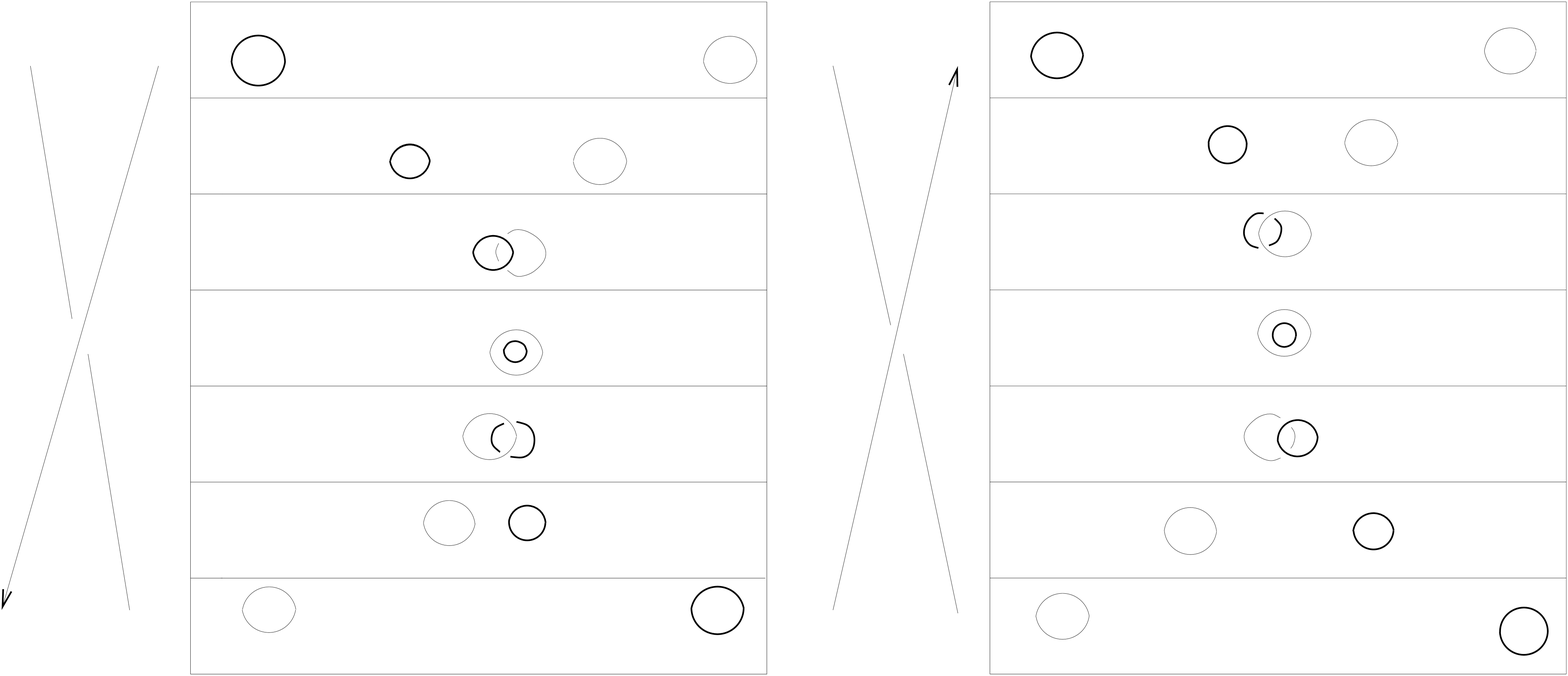}
\endrelabelbox }
\caption{\label{Tor2} Associating a knotted torus to a virtual knot: {classical crossing points}, first case. {All circles are oriented counterclockwise.} }
\end{figure}

\begin{figure}
\centerline{\relabelbox 
\epsfysize 7cm
\epsfbox {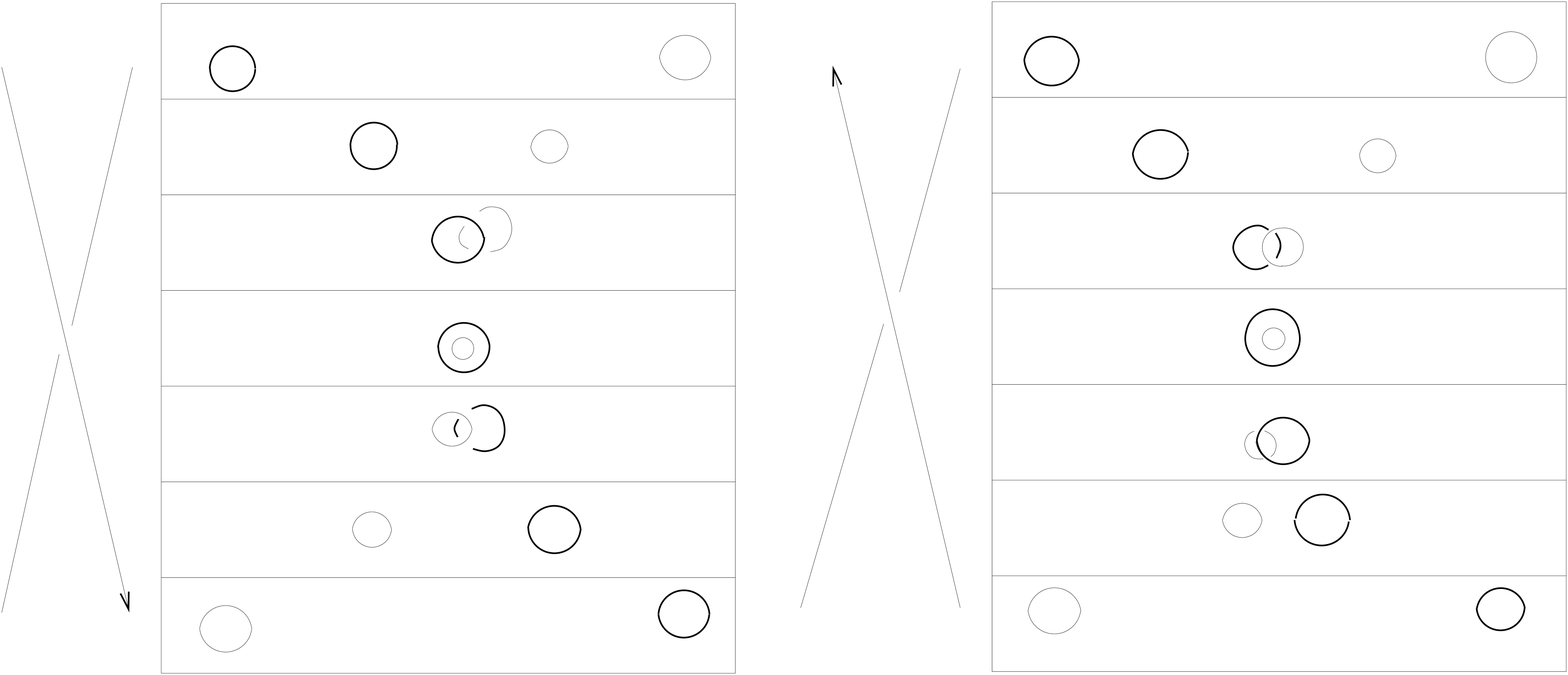}
\endrelabelbox }
\caption{\label{Tor3} Associating a knotted torus to a virtual knot: 
{classical crossing  points,} second case. {All circles are oriented counterclockwise.}}
\end{figure}

{For calculation purposes, however, it is important {to have a definition of the  ``Tube Map'' in terms of movies}.} 
Let $D\subset \R^2$ be an oriented virtual knot diagram. We can suppose, apart from planar
isotopy, that the projection on the second variable is a Morse function on
$D$. Define a movie of a knotted surface by using the correspondence of figures
{\ref{Tor1},  \ref{Tor2} and \ref{Tor3}.}
 Note our convention of reading movies of knotted surfaces from the  bottom to
 the top.  {This  yields an alternative way for describing the   tube $T(K)$ of the  virtual knot $K$, if we are provided a diagram of it.} 

It was proved in \cite{S,Ya} that the  correspondence $K \mapsto T(K)$, where $K$ is a  welded virtual knot,  preserves the  fundamental groups of the complement {(the knot {groups})}.

Given a (classical) link $K$ with $n$ components sitting in the interior of the semiplane
$\{(x,y,z) \in \R^3\colon z\ge 0\}$, we define the torus spun of $K$ by rotating $K$ 
4-dimensionally around the {plane $\{z=0\}$.} Therefore, we obtain an embedding
of the disjoint union of $n$ tori $S^1 \times S^1$ into $S^4$. It was shown
in \cite{S} that the torus spun of $K$ is in fact isotopic {to} the tube
$T(K)$ of $K$.

The correspondence $K \mapsto T(K)$ actually sends welded virtual links to
ribbon  torus links. In fact, any ribbon  torus link is of the form $T(K)$ for
some welded virtual knot $K$. However, it is an open problem  whether the map
$K \mapsto T(K)$ is faithful{; see} {\cite[problems (1) and (2) of 2.2.2]{CKS}.}

\subsubsection{Welded virtual arcs}\label{refer}
\begin{figure}
\centerline{\relabelbox 
\epsfysize 1.5cm
\epsfbox {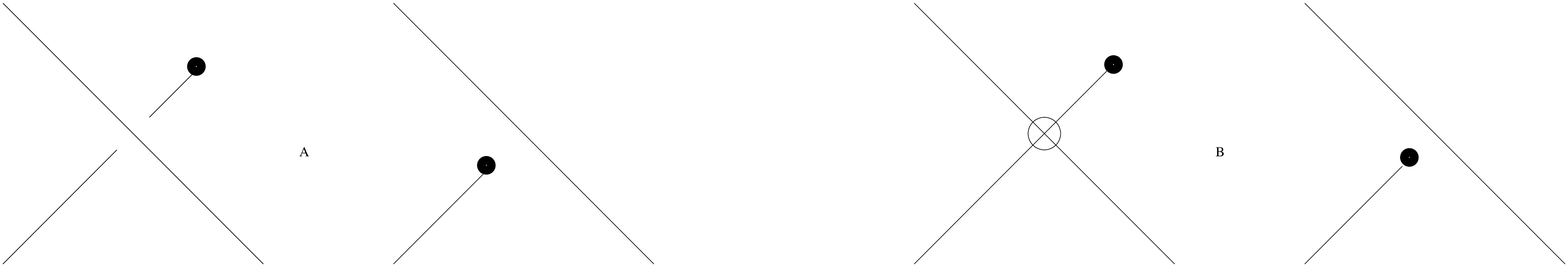}
\relabel{A}{$\leftrightarrow$}
\relabel{B}{$\leftrightarrow$}
\endrelabelbox }
\caption{\label{arc} Moves on {welded} virtual arc diagrams.}
\end{figure}

A virtual arc diagram is, by definition, an immersion of a disjoint union
of intervals $[0,1]$ into the plane $\R^2$, {where the 4-valent vertices of the immersion can represent either classical of virtual crossings.} The definition of a  welded virtual arc
is similar to the definition of a  welded virtual knot, but considering in
addition the moves of figure \ref{arc}{; see} \cite{S}.

\begin{figure}
\centerline{\relabelbox 
\epsfysize 2cm
\epsfbox {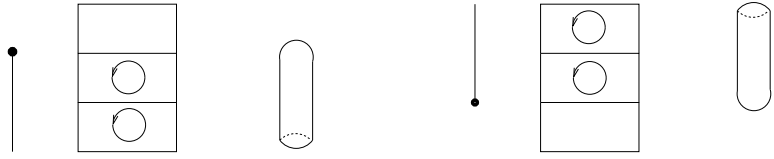}
\endrelabelbox }
\caption{ The tube of a welded {virtual} arc close to {the} endpoints. }
\label{ends}
\end{figure}

A  sphere link  is, by definition, an embedding of a disjoint union of spheres
$S^2$ into $S^4$, considered up to ambient isotopy. 
Similarly to ribbon torus links in $S^4$, any ribbon  sphere link admits a presentation as
the tube  $T(A)$, where $A$ is a welded virtual arc. {Here $T(A)$ is defined in
the same way as the tube of a welded virtual knot,}  considering
additionally  the movies of figure
\ref{ends} at the {end-points} of the arcs of $A$ . Therefore $T(A)$ is an embedding of a disjoint union of spheres $S^2$ into $S^4$.

Suppose that the arc $A$ is classical, and that it  sits inside the semiplane
$\{z\ge 0\}$ of $\R^3$, intersecting the plane $\{z=0\}$ at the end-points of
$A$, transversally. Then in fact $T(A)$ is the spun knot of $A${; see} \cite{R,S,CKS}.

We can define the knot group of a welded virtual  arc exactly in the same way  as we {defined} the combinatorial fundamental group of the complement of a welded virtual knot. {As in the case of welded virtual knots,  the map $A \mapsto T(A)$ preserves knot groups; see \cite{S}.} 

Suppose that $A$ is a classical arc {(with one component)} sitting in the semiplane
$\{z\ge 0\}$ of $\R^3$, intersecting the plane $\{z=0\}$ at the end-points of
$A$. Let $K$ be the obvious closure of $A$. Then it is easy to see that $A$ and $K$ have the same knot groups. Note that the fact that $A$ is classical is essential for this to hold. 
This is also true if $A$ may have some {$S^1$} components, even though it is strictly necessary that $A$ have only one  {component homeomorphic to $[0,1]$}.

\subsection{Crossed module invariants of knotted surfaces}\label{cmiks}

A crossed module (see \cite{B1}) $\G=\left (E \ra{\d} G,\t\right )$ is given by a group
morphism $\d\colon E \to G$ together with a left action $\t$ of $G$ on $E$ by automorphisms. The conditions on $\d$ and $\t$ are:
\begin{enumerate}
\item $\d(X \t e)=X \d(e) X^{-1}, \forall X \in G, \forall e \in E$,
\item $\d(e)\t f =e f e^{-1}, \forall e,f \in E$.
\end{enumerate}
Note that the second condition implies that the subgroup $\ker{\d}$ of $E$ is
central in $E$, whereas the {first} implies that $\ker \d$ is $G$-invariant.

\begin{figure}
\centerline{\relabelbox 
\epsfysize 2cm
\epsfbox{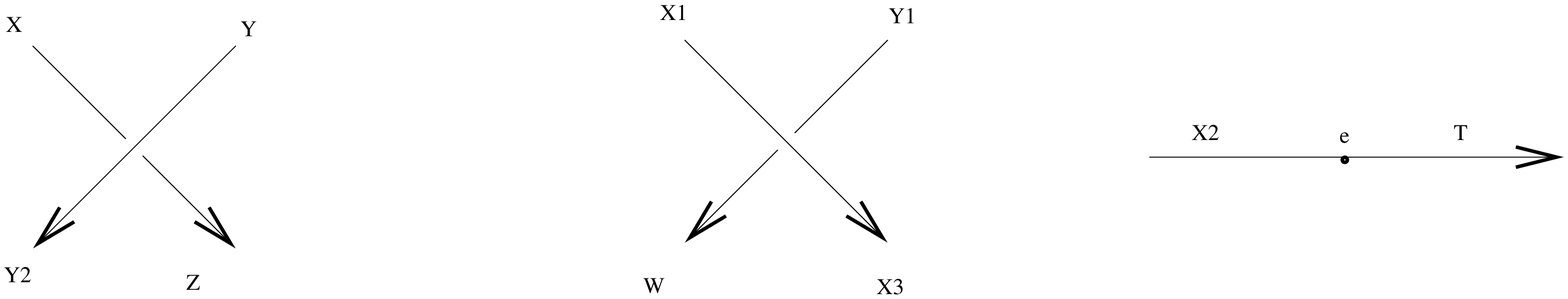} 
\relabel {X1}{$\s{X}$}
\relabel {X2}{$\s{X}$}
\relabel {X}{$\s{X}$}
\relabel {X3}{$\s{X}$}
\relabel {Y1}{$\s{Y}$}
\relabel {Y}{$\s{Y}$}
\relabel {Y2}{$\s{Y}$}
\relabel {Z}{$\s{Y^{-1}XY}$}
\relabel {W}{$\s{XYX^{-1}}$}
\relabel {e}{$\s{e}$}
\relabel {T}{$\s{\d(e)X}$}
\endrelabelbox }
\caption{Definition of a colouring of a dotted knot diagram.}
\label{Colour}
 \end{figure}

 A  dotted knot diagram is,  by definition, a regular projection of a bivalent
 graph, in other words of a link, possibly with some extra bivalent vertices inserted. 
{Let $D$ be a dotted knot diagram,} which we suppose to be oriented. 
 Let also  $\G=\left ( E \ra{\d} G,\t\right)$ be a  finite crossed module.
\begin{Definition}
A colouring of $D$ is an assignment of an element of $G$ to each arc of $D$
and of an element  of $E$ to each bivalent vertex of $D$ {satisfying} the conditions of figure \ref{Colour}.\end{Definition}

\begin{Definition}
Let $D$ be a knot diagram {(without vertices)}. A dotting of $D$ is an insertion of bivalent
vertices {in $D$,} considered up to a planar isotopy sending $D$ to $D$,
setwise. If $D$ is an oriented knot diagram, let $V(D)$ be the free $\Q$-vector space on the set of all colourings of all dottings of $D$.  \end{Definition}

\begin{figure}
\centerline{\relabelbox 
\epsfysize 2cm
\epsfbox{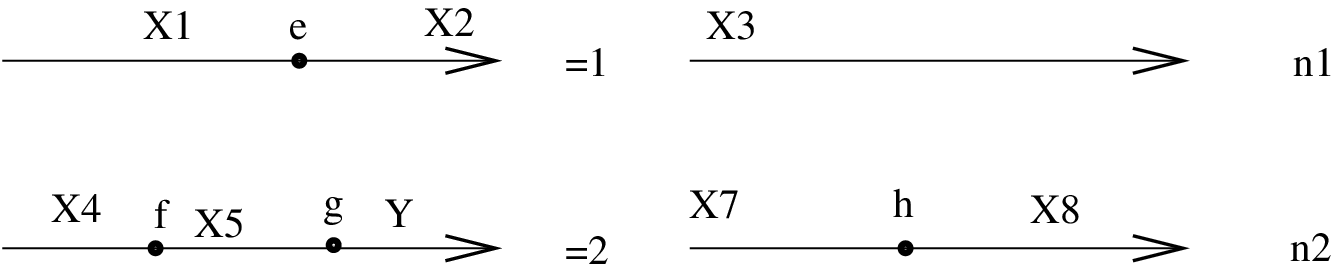} 
\relabel {X1}{$\s{X}$}
\relabel {X2}{$\s{X}$}
\relabel {X3}{$\s{X}$}
\relabel {X4}{$\s{X}$}
\relabel {X5}{$\s{\d(e)X}$}
\relabel {Y}{$\s{\d(fe)X}$}
\relabel {X7}{$\s{X}$}
\relabel {X8}{$\s{\d(fe)X}$}
\relabel {e}{$\s{1_E}$}
\relabel {f}{$\s{e}$}
\relabel {g}{$\s{f}$}
\relabel {h}{$\s{fe}$}
\relabel{=2}{$\s{=}$}
\relabel{=1}{$\s{=}$}
\relabel{n1}{$\s{R1}$}
\relabel{n2}{$\s{R2}$}
\endrelabelbox }
\caption{Relations on colourings.}
\label{Relations1}
\end{figure}
\begin{figure}
\centerline{\relabelbox 
\epsfysize 4cm
\epsfbox{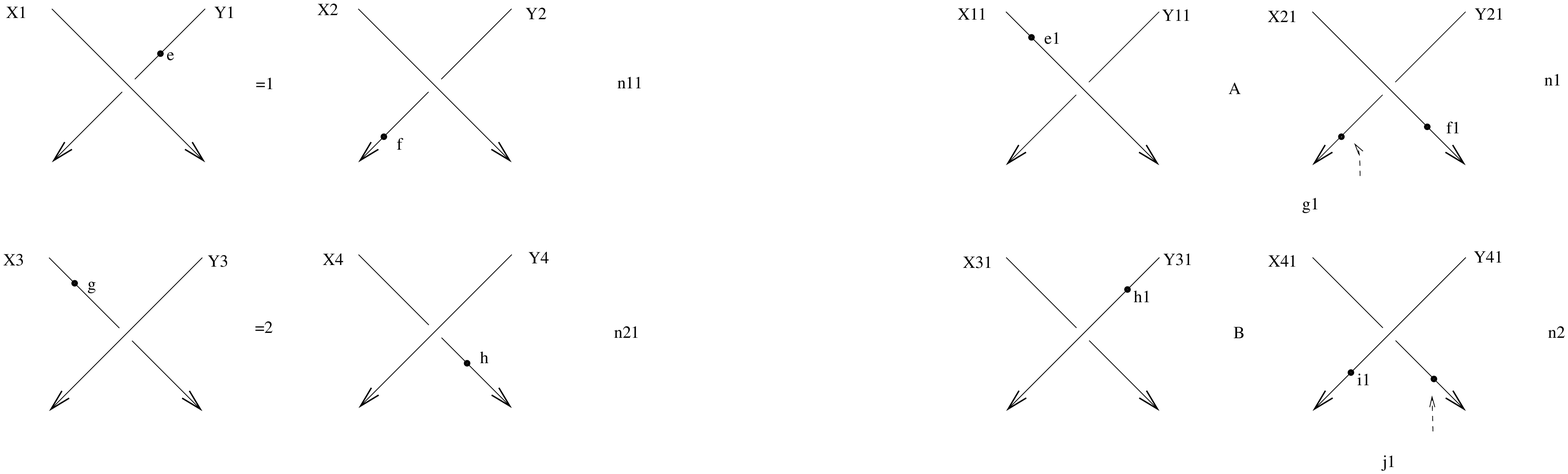} 
\relabel {X1}{$\s{X}$}
\relabel {X2}{$\s{X}$}
\relabel {X3}{$\s{X}$}
\relabel {X4}{$\s{X}$}
\relabel {Y1}{$\s{Y}$}
\relabel {Y2}{$\s{Y}$}
\relabel {Y3}{$\s{Y}$}
\relabel {Y4}{$\s{Y}$}
\relabel {e}{$\s{e}$}
\relabel {f}{$\s{X\t e}$}
\relabel {g}{$\s{e}$}
\relabel {h}{$\s{Y^{-1} \t e}$}
\relabel{=2}{$\s{=}$}
\relabel{=1}{$\s{=}$}
\relabel{n1}{$\s{R5}$}
\relabel{n2}{$\s{R6}$}
\relabel {X11}{$\s{X}$}
\relabel {X21}{$\s{X}$}
\relabel {X31}{$\s{X}$}
\relabel {X41}{$\s{X}$}
\relabel {Y11}{$\s{Y}$}
\relabel {Y21}{$\s{Y}$}
\relabel {Y31}{$\s{Y}$}
\relabel {Y41}{$\s{Y}$}
\relabel {e1}{$\s{e}$}
\relabel {f1}{$\s{e}$}
\relabel {g1}{$\s{eXYX^{-1}\t e^{-1}}$}
\relabel {h1}{$\s{e}$}
\relabel {i1}{$\s{e}$}
\relabel {j1}{$\s{Y^{-1}\t e^{-1}Y^{-1}X \t e}$}
\relabel {A}{$\s{=}$}
\relabel {B}{$\s{=}$}
\relabel {n11}{$\s{R3}$}
\relabel {n21}{$\s{R4}$}
\endrelabelbox }
\caption{Relations on colourings.}
\label{Relations3}
\end{figure}

Consider now the relations of figures \ref{Relations1} and \ref{Relations3}. It is straightforward to see that they are local on the knot
diagrams and that they transform colourings into colourings.

\begin{Definition}
Let $D$ be an oriented knot diagram (without vertices). The vector space $\V(D)$ is defined as the vector space obtained from $V(D)$ by {modding out by} the relations $R1$ to $R6$.
\end{Definition}

\begin{figure}
\centerline{\relabelbox 
\epsfysize 1.5cm
\epsfbox{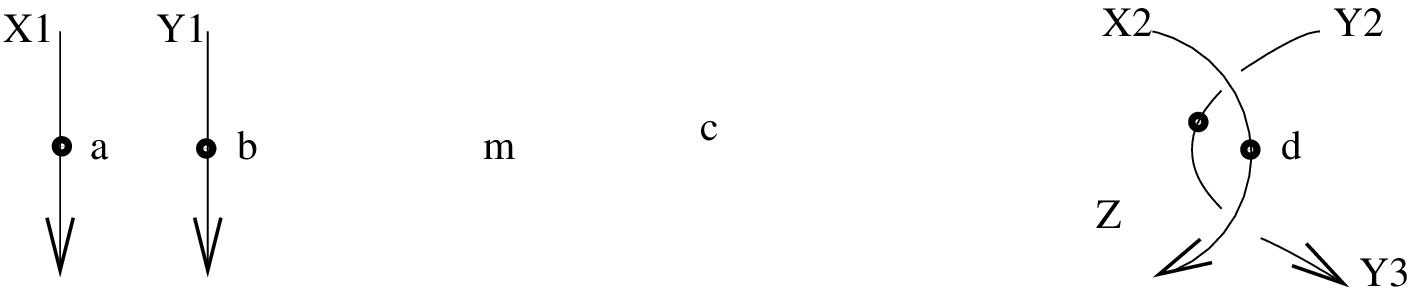} 
\relabel{X1}{$\s{X}$}
\relabel{X2}{$\s{X}$}
\relabel{Y1}{$\s{Y}$}
\relabel{Y2}{$\s{Y}$}
\relabel{m}{$\longmapsto$}
\relabel{a}{$\s{e}$}
\relabel{b}{$\s{f}$}
\relabel{c}{$\s{e X \t f XYX^{-1} \t e^{-1}}$}
\relabel{d}{$\s{e}$}
\endrelabelbox }
\caption{Map assigned to positive Reidemeister-II move.}
\label{MapR2}
\end{figure}

\begin{figure}
\centerline{\relabelbox 
\epsfysize 1.5cm
\epsfbox{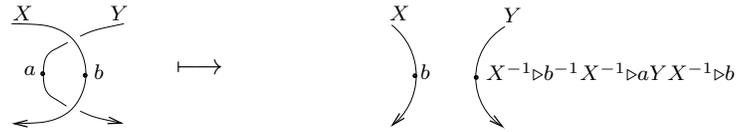} 
\relabel{A}{$\s{X}$}
\relabel{B}{$\s{Y}$}
\relabel{e}{$\s{a}$}
\relabel{f}{$\s{b}$}
\relabel{E}{$\s{X}$}
\relabel{F}{$\s{Y}$}
\relabel{c}{$\s{b}$}
\relabel{d}{$\s{X^{-1} \t b^{-1} X^{-1} \t a YX^{-1} \t b  }$}
\relabel{t}{$\longmapsto$}
\endrelabelbox }
\caption{Map assigned to negative Reidemeister-II move.}
\label{MapR2neg}
\end{figure}

\begin{figure}
\centerline{\relabelbox 
\epsfysize 1.5cm
\epsfbox{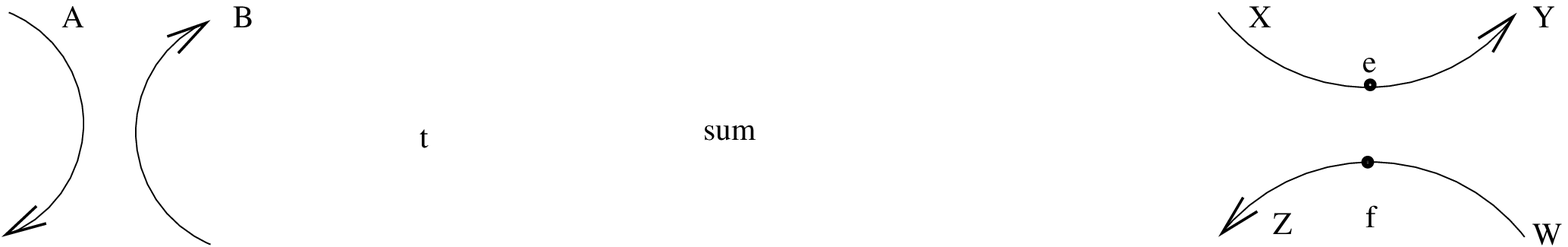}
\relabel{t}{$\longmapsto$}
\relabel{A}{$\s{X}$}
\relabel{B}{$\s{Y}$}
\relabel{X}{$\s{X}$}
\relabel{Y}{$\s{\d(e)X}$}
\relabel{Z}{$\s{X}$}
\relabel{W}{$\s{\d(e)X}$}
\relabel{e}{$\s{e}$}
\relabel{f}{$\s{e^{-1}}$}
\relabel{sum}{$\displaystyle{\frac{1}{\# E}{\sum_{e \in E} \delta(Y,\d(e)X)}}$}
\endrelabelbox}
\caption{Map associated to saddle point moves.}
\label{Mapsaddle}
\end{figure}

\begin{figure}
\centerline{\relabelbox 
\epsfysize 1.5cm
\epsfbox{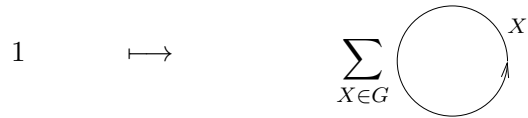}
\relabel{1}{$1$}
\relabel{t}{$\longmapsto$}
\relabel{s}{$\displaystyle{\sum_{X \in G}}$}
\relabel{X}{$\s{X}$}
\endrelabelbox}
\caption{Map associated with births of a circle.}
\label{Mapcap}
\end{figure}

\begin{figure}
\centerline{\relabelbox 
\epsfysize 1.5cm
\epsfbox{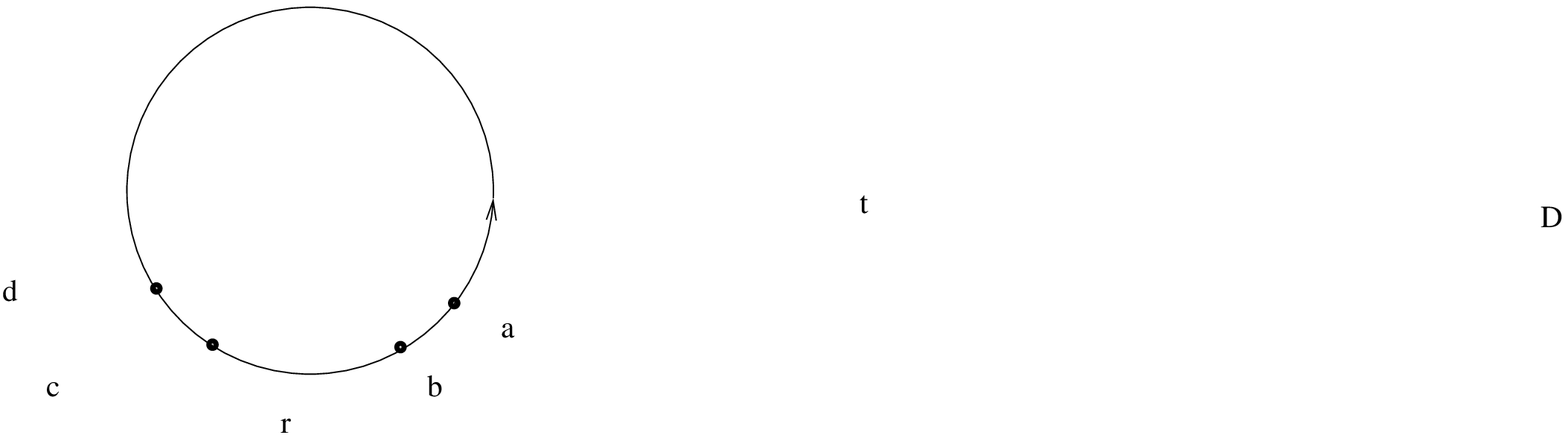}
\relabel{t}{$\longmapsto$}
\relabel{a}{$\s{x_1}$}
\relabel{b}{$\s{x_2}$}
\relabel{c}{$\s{x_{n-1}}$}
\relabel{d}{$\s{x_n}$}
\relabel{r}{$\dots$}
\relabel{D}{$\# E \delta(x_1x_2...x_{n-1}x_n,1_E)$}
\endrelabelbox}
\caption{Map associated with deaths of a circle.}
\label{Mapcup}
\end{figure}

\begin{figure}
\centerline{\relabelbox 
\epsfxsize 12cm
\epsfbox{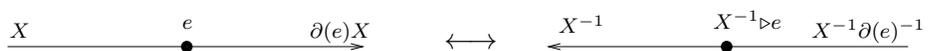}
\relabel{l}{$\longleftrightarrow$}
\relabel{X}{$\s{X}$}
\relabel{Y}{$\s{\d(e)X}$}
\relabel{Z}{$\s{X^{-1}}$}
\relabel{W}{$\s{X^{-1} \d (e)^{-1}}$}
\relabel{e}{$\s{e}$}
\relabel{f}{$\s{X^{-1} \t e}$}
\endrelabelbox}
\caption{Inversion of strands.}
\label{inversion}
\end{figure}

Let $D$ and $D'$ be oriented knot diagrams. If $D$ and $D'$ differ by planar
isotopy, then there exists an obvious map $\V(D) \to \V(D')$. In fact, if $D$
and $D'$ differ by a Reidemeister move or a Morse move (in other words a
birth/death of a circle or a saddle point), then there also exists a well
defined map $\V(D) \to \V(D')$. All this is explained in \cite{FM1}. In
figures \ref{MapR2}, \ref{MapR2neg}, \ref{Mapsaddle}, \ref{Mapcap} and \ref{Mapcup} we display
the definition of these maps for the case of the Reidemeister-II move and the
Morse moves, which we are going to need in this article.  The remaining cases of these moves 
can be dealt with by doing the transition shown in figure \ref{inversion}, and using the relations  $R1$ to $R6$. {In figure {\ref{Mapsaddle}}, $\delta$ is a Kronecker delta.}

Therefore, any movie of an oriented knotted surface $\S$ can be evaluated to give an element $I_\G(\S) \in \Q$. 
\begin{Theorem}
The evaluation $I_\G$ of a movie of an oriented knotted surface defines  an
isotopy invariant of oriented knotted surfaces.
\end{Theorem}
This is shown in \cite{FM1}. The homotopy theoretical interpretation of {the isotopy invariant $I_\G$} is
discussed in \cite{FM2,FM3,FMP}. {The construction of {the invariant $I_\G$} was initially inspired by Yetter's Invariant of manifolds; see \cite{Y3,P1,P2}.} 

Actually $I_\G$ defines an embedded TQFT, in other words, an invariant of link cobordisms considered up to ambient isotopy fixing both ends. 

\subsubsection{The case of ribbon knotted torus}

As we have seen, if $\S$ is a ribbon knotted surface, which topologically is
the disjoint union of tori $S^1 \times S^1$ or spheres $S^2$, then we can
represent it as the tube $T(K)$ of  welded virtual knot $K$, in the first case, or the tube $T(A)$ of a welded virtual
arc $A$, in the second case.

\begin{figure}
\centerline{\relabelbox 
\epsfysize 1.5cm
\epsfbox{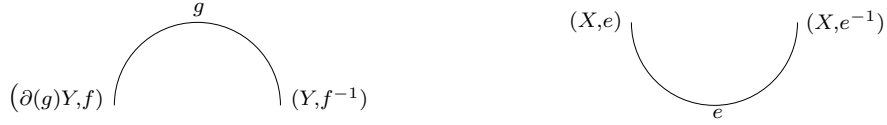}
\relabel{A}{$(\s{\d(g)Y,f)}$}
\relabel{B}{$\s{(Y,f^{-1})}$}
\relabel{D}{$\s{(X,e)}$}
\relabel{E}{$\s{(X,e^{-1})}$}
\relabel{e}{$\s{e}$}
\relabel{g}{$\s{g}$}
\endrelabelbox}
\caption{Relations at maximal and minimal points.}
\label{minmax}
\end{figure}

\begin{figure}
\centerline{\relabelbox 
\epsfysize 10cm
\epsfbox{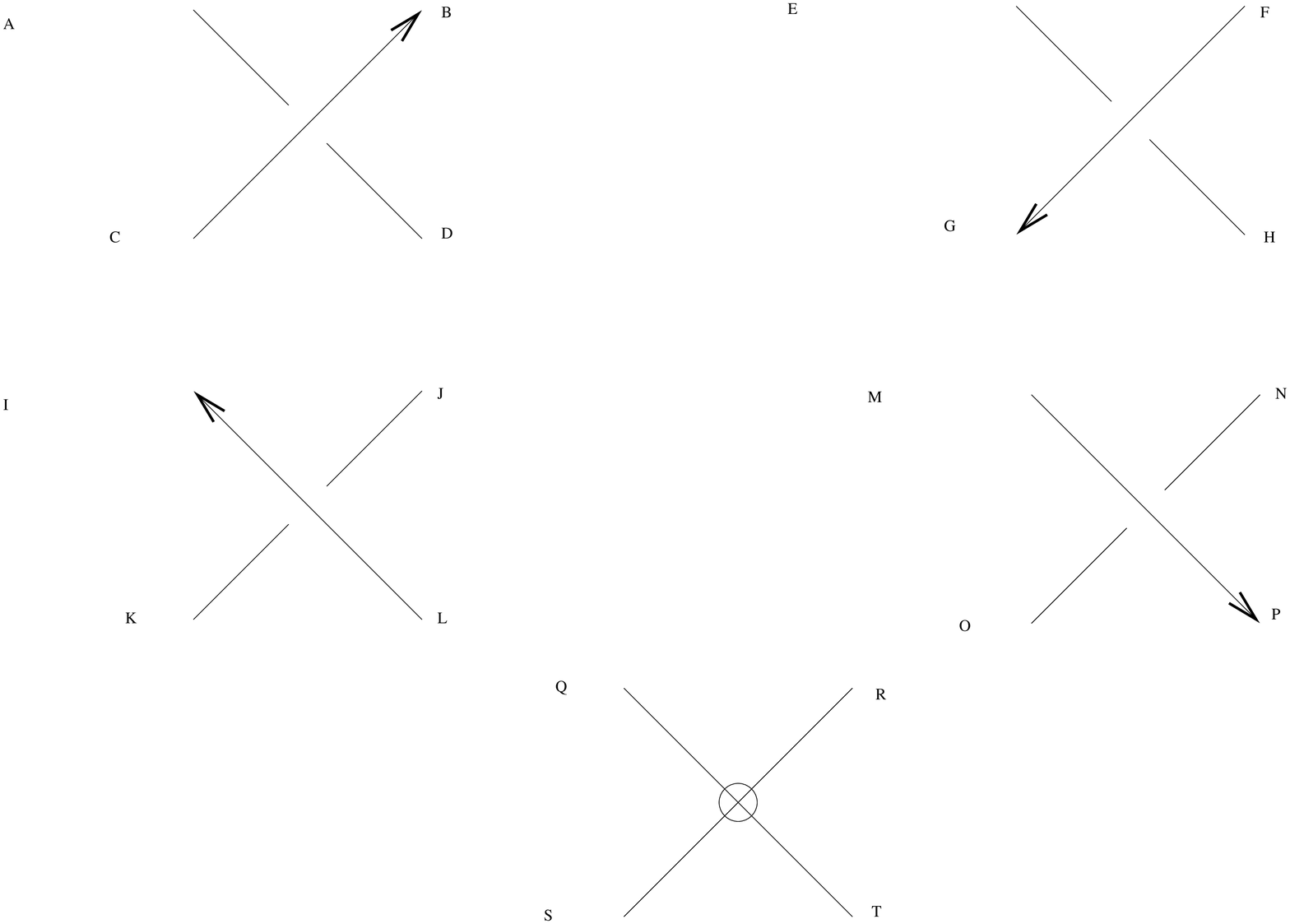}
\relabel{A}{$\s{(XYX^{-1},X \t f)}$}
\relabel{B}{$\s{(X,efX \t f^{-1})}$}
\relabel{C}{$\s{(X,e)}$}
\relabel{D}{$\s{(Y,f)}$}
\relabel{E}{$\s{(X^{-1}YX, X^{-1} \t f) }$}
\relabel{F}{$\s{(X,X^{-1} \t f^{-1} e f )}$}
\relabel{G}{$\s{(X,e)}$}
\relabel{H}{$\s{(Y,f)}$}
\relabel{I}{$\s{(Y, Y^{-1} \t e^{-1} e f)}$}
\relabel{J}{$\s{(Y^{-1}XY,Y^{-1} \t e) }$}
\relabel{K}{$\s{(X,e)}$}
\relabel{L}{$\s{(Y,f)}$}
\relabel{M}{$\s{(Y,feY \t e^{-1})}$}
\relabel{N}{$\s{(YXY^{-1}, Y \t e)}$}
\relabel{O}{$\s{(X,e)}$}
\relabel{P}{$\s{(Y,f)}$}
\relabel{Q}{$\s{(Y,f)}$}
\relabel{R}{$\s{(X,e)}$}
\relabel{S}{$\s{(X,e)}$}
\relabel{T}{$\s{(Y,f)}$}
\endrelabelbox}
\caption{Relations at crossings.  }
\label{Crossings}
\end{figure}

 We want to find an algorithm {for calculating} $I_\G(T(K))$, where $K$ is a {welded virtual knot}, directly from {a diagram of $K$ itself,} and analogously for a welded virtual arc $A$. 
A careful look at the definition {of the invariant $I_\G$} together with {the definition of the tube map in \ref{presentation}} leads
to the following definition:

\begin{Definition}
Let $\G=\left(E \ra{\d} G,\t\right )$ be a crossed module. Let also $D$ be a welded
virtual knot  diagram. Suppose that the projection on the second
variable defines  a Morse function on $D$.  A $\G$-colouring\footnote{{This should not be confused with the notion of a colouring which was considered in the definition of the invariant $I_\G$, above.}} of $D$ is an assignment of a
pair $(X,f)$, where $X \in G$ and $f \in \ker \d$, to each connected component
of $D$ minus its set of crossings and extreme points; of an element $e \in \ker \d$ to each
minimal point; and an element $g \in E$ to each {maximal point,} {satisfying} the
conditions shown in figures \ref{minmax} and \ref{Crossings}.
\end{Definition}

\begin{figure}
\centerline{\relabelbox 
\epsfysize 5cm
\epsfbox{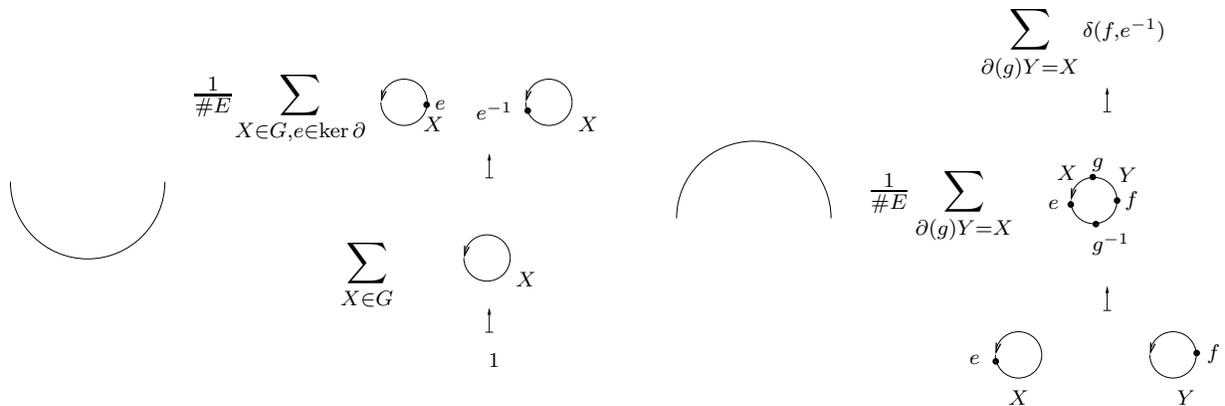}
\relabel {m}{$\s{e}$}
\relabel {n}{$\s{f}$}
\relabel {1}{$\s{1}$}
\relabel {X}{$\s{X}$}
\relabel {sum2}{$\s{\displaystyle{{\sum_{X \in G}}}}$}
\relabel {sum1}{$\frac{1}{\# E}{\displaystyle{{\sum_{X \in G, e \in \ker \d}}}}$}
\relabel {Y}{$\s{X}$}
\relabel {Z}{$\s{X}$}
\relabel {e}{$\s{e}$}
\relabel {f}{$\s{e^{-1}}$}
\relabel {A}{$\s{X}$}
\relabel {B}{$\s{Y}$}
\relabel {p}{$\s{e}$}
\relabel {q}{$\s{f}$}
\relabel {s}{$\s{g^{-1}}$}
\relabel {r}{$\s{g}$}
\relabel {C}{$\s{X}$}
\relabel {D}{$\s{Y}$}
\relabel {rel}{${\frac{1}{\# E}\displaystyle\sum_{\d(g)Y=X}}$}
\relabel {11}{$\displaystyle{\sum_{\d(g)Y=X} \s{\delta(f,e^{-1})}}$}
\endrelabelbox}
\caption{ Calculation of $I_\G$ of the tube of a welded virtual knot: minimal and maximal points. }
\label{cupcapcalc}
\end{figure}

\begin{figure}
\centerline{\relabelbox 
\epsfysize 10cm
\epsfbox{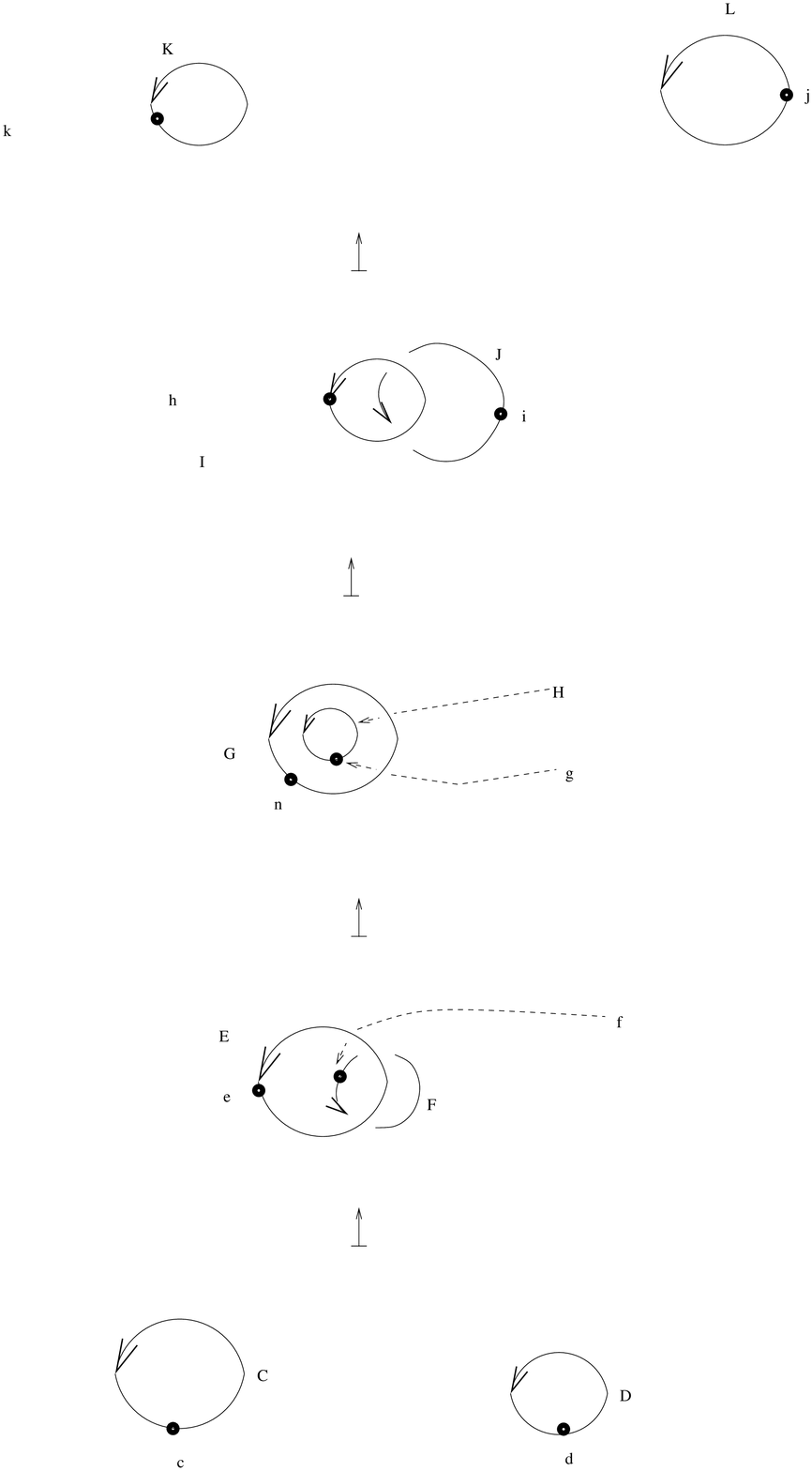}
\relabel {C}{$\s{X}$}
\relabel {c}{$\s{e}$}
\relabel {D}{$\s{Y}$}
\relabel {d}{$\s{f}$}
\relabel {E}{$\s{X}$}
\relabel {F}{$\s{Y}$}
\relabel {e}{$\s{e}$}
\relabel {f}{$\s{X^{-1} \t f }$}
\relabel {G}{$\s{X}$}
\relabel {H}{$\s{X^{-1} Y X}$}
\relabel {g}{$\s{X^{-1} \t f }$}
\relabel {h}{$\s{X^{-1} \t f}$}
\relabel {I}{$\s{X^{-1}YX}$}
\relabel {J}{$\s{X}$}
\relabel {i}{$ \s{eX^{-1} \t f^{-1} f}  $}
\relabel {j}{$\s{eX^{-1} \t f^{-1} f}$}
\relabel {K}{$\s{X^{-1}YX}$}
\relabel {L}{$\s{X}$}
\relabel {k}{$\s{X^{-1} \t f}$}
\relabel {n}{$\s{e}$}
\endrelabelbox}
\caption{ Calculation of $I_\G$ of the tube {$T(D)$} of a welded virtual knot {$D$: the type of crossings} relative to figure \ref{Tor2}.}
\label{CrossCalc}
\end{figure}

\begin{figure}
\centerline{\relabelbox 
\epsfysize 1.5cm
\epsfbox{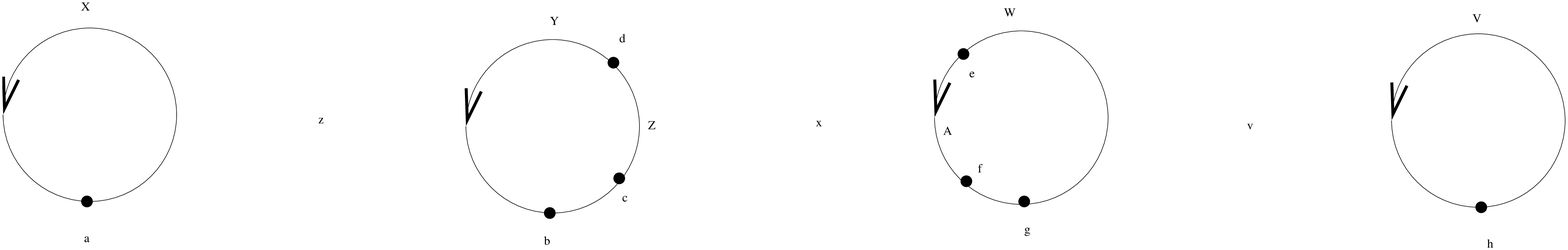}
\relabel {X}{$\s{X}$}
\relabel {a}{$\s{e}$}
\relabel {Y}{$\s{X}$}
\relabel {Z}{$\s{\d(g)X}$}
\relabel {d}{$\s{g^{-1}}$}
\relabel {c}{$\s{g}$}
\relabel {b}{$\s{e}$}
\relabel {W}{$\s{\d(g)X}$}
\relabel {e}{$\s{g^{-1}}$}
\relabel {f}{$\s{g}$}
\relabel {g}{$\s{e}$}
\relabel {A}{$\s{X}$}
\relabel {V}{$\s{\d(g)X}$}
\relabel {h}{$\s{e}$}
\relabel {z}{$\s{=}$}
\relabel {x}{$\s{=}$}
\relabel {v}{$\s{=}$}
\endrelabelbox}
\caption{An identity. Here $e \in \ker \d$.}
\label{change}
\end{figure}

The reason {for  considering} these relation is obvious from figure \ref{cupcapcalc},
and figure \ref{CrossCalc}, and its counterparts for different types of
crossings. {Note that $\ker \d \subset E$ is central in $E$.} However, for this calculus to approximate the definition of  $I_\G(T(D))$, {for $D$ a virtual knot diagram,} the
relation of figure \ref{change} still needs to be incorporated  {into} the
calculations.  To avoid needing to involve this relation, we consider the
following restriction on the crossed modules with which  we work.

\begin{Definition} [Automorphic Crossed Module]
A crossed module $\G=\left ( E \ra{\d} G,\t \right )$ is called automorphic if
  $\d(e)=1, \forall e \in E$. Therefore, an automorphic crossed module is
  given simply by two groups $G$ and $E$, with $E$ abelian, and a left  action
  $\t$ of
  $G$ on $E$ by automorphisms.
\end{Definition}

\begin{figure}
\centerline{\relabelbox 
\epsfysize 1.5cm
\epsfbox{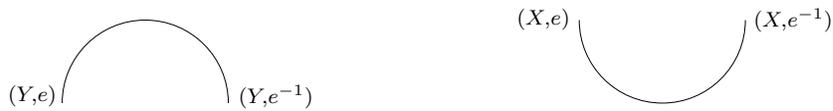}
\relabel{A}{$\s{(Y,e)}$}
\relabel{B}{$\s{(Y,e^{-1})}$}
\relabel{D}{$\s{(X,e)}$}
\relabel{E}{$\s{(X,e^{-1})}$}
\endrelabelbox}
\caption{Reduced $\G$-colouring at extreme points.}
\label{extred}
\end{figure}

\begin{Definition}[Reduced $\G$-Colourings]
Let $\G=(E,G,\t)$ be an automorphic crossed module. Let also $D\subset \R^2$ be a virtual
knot diagram, such that the projection on the second variable is a Morse
function on $D$. A reduced  $\G$-colouring of $D$ is given by an assignment
of a pair $(X,e)\in G \times E$ to each connected component of $D$ minus its set of crossings and extreme points, {satisfying} the relations
of figures \ref{Crossings} and \ref{extred}.
\end{Definition}

{The following result is easy to prove  by using all the information we provided,} and the fact that{, for any knot diagram,} the number of minimal points of it equals {the} number of maximal points. 

\begin{Theorem}
Let $D$ be a  virtual knot diagram, such that the projection on the second
variable is a Morse function on $D$. Let also $\G=\left (E, G,\t \right)$  be
a finite automorphic crossed module. Consider the quantity:
$$\H_\G(D)={\# \{\textrm{reduced }\G\textrm{-colourings of } D\}}.$$
Then {$\H_\G(D)$} is an invariant of welded  {virtual} knots. In fact:
$${\H_\G(D)=I_\G(T(D)).}$$
Here $I_\G$ is the Crossed Module Invariant of oriented knotted surfaces defined in \cite{FM1}.

\end{Theorem}

\begin{figure}
\centerline{\relabelbox 
\epsfysize 1.7cm
\epsfbox{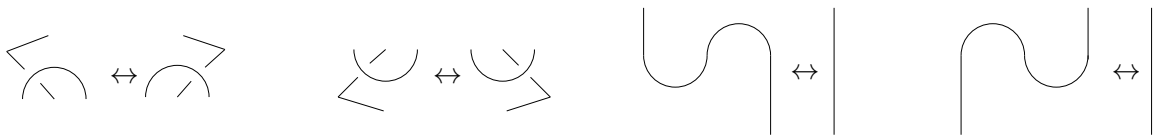}
\relabel{a}{$\leftrightarrow$}
\relabel{b}{$\leftrightarrow$}
\relabel{c}{$\leftrightarrow$}
\relabel{d}{$\leftrightarrow$}
\endrelabelbox}
\caption{Sample of Yetter's moves capturing planar isotopy.}
\label{Yetter}
\end{figure}

\begin{Exercise}\label{Ex1}
Check directly that $\H_\G$ {(where $\G$ is an automorphic {finite} crossed module)} is an invariant of welded virtual knots. Note that
together with the moves defining welded virtual knots, we still need to check
invariance under planar isotopy, thus enforcing us to check invariance under
the moves  of the type depicted in figure \ref{Yetter}, usually called
Yetter's Moves{;  see \cite{Y2,FY}.} It
is important to note that we need to consider all the possible different
crossing informations, and, since we are working in the oriented case, all
the possible orientations of the strands.     
\end{Exercise}

Let $\G=\left (E \ra{\d} G,\t \right)$ be a crossed module. Define $\pi_1(\G)=\mathrm{coker}(\d)$ and
$\pi_2(\G)=\ker \d$, which is an abelian group. Then $\pi_1(\G)$ has a natural
left action $\t'$ on
$\pi_2(\G)$ by automorphisms.  In particular
$\Pi(\G)=(\pi_2(\G),\pi_1(\G),\t')$ is an automorphic crossed module.
In fact $\G$ also determines a cohomology class $k^3 \in
{H^3(\pi_1(\G), \pi_2(\G))}$, called the $k$-invariant of $\G$.

It is not difficult  to extend the invariant $\H_\G(D)$, where $D$ is a welded
virtual knot, to handle non-automorphic crossed modules $\G$, so that
$\H_\G(D)=I_\G(T(D))$. We do this  by incorporating the
relation in figure \ref{change} into the notion of a $\G$-colouring of a
virtual knot diagram. However, it
is possible to prove that for any {welded virtual} knot $D$ and any finite crossed module $\G$ we have that 
$I_\G(T(D))$ equals $I_{\Pi(\G)}(T(D))$, apart from {normalisation} factors. {This can be proved by using the graphical framework presented in this article.}
 Hence, we do not lose generality if we restrict our attention only to  automorphic crossed modules.

\begin{Problem} Let $\G=(E,G,\t)$ be an automorphic crossed module. Find a
  ribbon Hopf algebra $\A_\G$ acting on the  vector space freely generated by $G
  \times E$ such that $\H_\G$ is the Reshetikhin-Turaev invariant of knots
  associated to it {(see \cite{RT}),} and so that the case of welded virtual knots also follow from this
  Hopf algebra framework in a natural way. Note that in the case {when} $E=0$, we can take
  $\A_\G$ to be the quantum double of the function algebra on $G$. The solution to this problem would be somehow the quantum double of a finite categorical group, and therefore would be of considerable importance.
\end{Problem}

\subsubsection{The case of welded virtual arcs}\label{weldedarcs}

\begin{figure}
\centerline{\relabelbox 
\epsfysize 1.5cm
\epsfbox{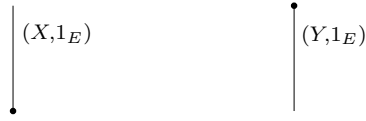}
\relabel{X}{$\s(X,1_E)$}
\relabel{Y}{$\s(Y,1_E)$}
\endrelabelbox}
\caption{Reduced {$\G$-colourings} of welded virtual arcs at end-points. Here $X,Y \in
  G$.} 
\label{endscoloured}
\end{figure}

Let $\G=(E,G,\t)$ be {a} finite automorphic crossed module. Let also $A$ be a virtual
arc diagram.  The notion of a reduced $\G$-colouring of $A$ is totally analogous to
the concept of a reduced $\G$-colouring of a virtual knot diagram, considering that if
an {arc} of $A$ has a free end then it must be coloured by $(X,1_E)$, where $X
\in G${; see} figure \ref{endscoloured}. One can see this from figure
\ref{ends}. We have:

\begin{Theorem}\label{arccalc}
Let $A$ be a virtual arc diagram. The quantity:
$$\H_\G(A)=\frac{\#\{\textrm{reduced } \G\textrm{-colourings of } A\}}{\# E^{\#\{\textrm{cups}\}  -  \#\{\textrm{caps}\} -  \#\{\textrm{pointing upwards ends of }A\}  }}$$
is an invariant of $A$ as a {welded virtual arc.} In fact
$$\H_\G(A)=I_\G(T(A)).$$

\end{Theorem}
Therefore, the graphical  framework presented in this article is also a
 calculational device {for calculating} {the
crossed module invariant} of spun knots, accordingly to \ref{refer}.

The invariant $\H_\G$  of  Theorem \ref{arccalc}   actually is an invariant of virtual arcs of which some components may be circles. In fact, it also naturally  extends to an invariant of welded virtual graphs, to be defined in \ref{FMK}.

\section{Examples}

\subsection{Virtual and Classical Hopf Link}\label{vchl}

\begin{figure}
\centerline{\relabelbox 
\epsfysize 1.5cm
\epsfbox{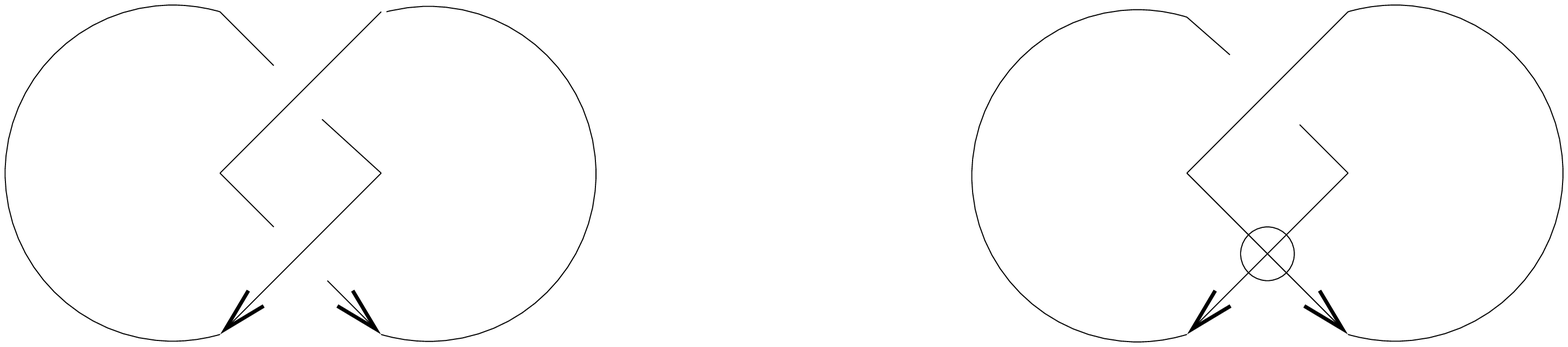}
\endrelabelbox}
\caption{Classical and Virtual Hopf links.}
\label{vhl}
\end{figure}
\subsubsection{Virtual Hopf Link}
The simplest non-trivial  welded virtual link  is the Virtual Hopf Link $L$,
depicted in  figure \ref{vhl}. Note that $L$ is linked since its {knot group is $\{X,Y\colon XY=YX\}\cong \Z^2$.}

\begin{figure}
\centerline{\relabelbox 
\epsfysize 3.2cm
\epsfbox{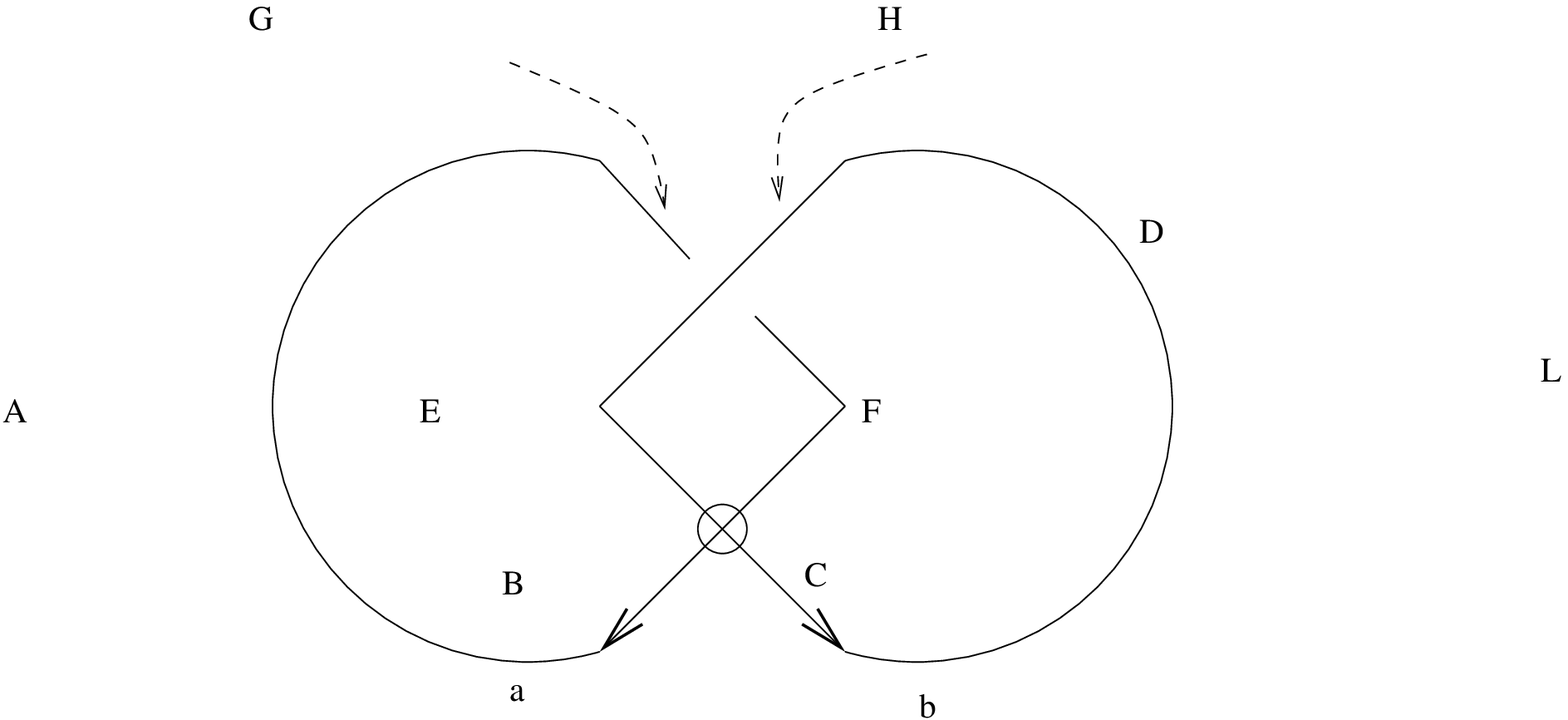}
\relabel {A}{$\s{(X,e^{-1})}$}
\relabel {B}{$\s{(X,e)}$}
\relabel {C}{$\s{(Y,f)}$}
\relabel {D}{$\s{(Y,f^{-1})}$}
\relabel {E}{$\s{(Y,f)}$}
\relabel {F}{$\s{(X,e)}$}
\relabel {G}{$\s{(Y^{-1}XY,Y^{-1} \t e)}$}
\relabel {H}{$\s{(Y,Y^{-1}\t e^{-1} ef)}$}
\relabel {L}{$\left\{ \begin{CD}& \s{Y^{-1}\t e=e} \\ &\s{Y^{-1}XY=X} \\ &\s{Y^{-1} \t
      e^{-1}ef=f }\end{CD}\right. $}
 \endrelabelbox}
\caption{Calculation  of the crossed module invariant of the Virtual Hopf Link $L$.} 
\label{vlhcalc}
\end{figure}

Let $\G=(E,G,\t)$ be a finite automorphic crossed module. Let us calculate the
crossed module invariant $\H_\G$ of the Virtual Hopf Link $L$. This calculation appears in figure
\ref{vlhcalc}.  From this we can conclude that:
\begin{align}
\H_\G(L)&={\#\{X,Y \in G; e,f \in E| XY=YX,  Y^{-1} \t e =e \}}\\ \label{vhlv}
&={\# E}\# \{X,Y \in G; e \in E| XY=YX,  Y^{-1} \t e =e \}.
\end{align}

Note that the previous equation  simplifies to $${\H_\G(L)=\#E\#G{\#\{Y \in G; e \in E|  Y^{-1} \t e =e \}}},$$
when the group $G$ is abelian.
On the other hand it is easy to see that if $O^2$ is a pair of unlinked unknots
then we have:
\begin{equation}\label{tlv}
\H_\G(O^2)=\#G^2 \#E^2.
\end{equation}
 From
 equations (\ref{vhlv}) and (\ref{tlv}), it thus follows that any finite  automorphic crossed module
 $(E,G,\t)$ with $G$ abelian sees the knotting of the Virtual Hopf Link if
 there exists $Y \in G$ and $e \in E$ such that $Y^{-1} \t e \neq e$. This is
 verified in  any automorphic crossed module $(E,G,\t)$ with $\t$ being a
 non-trivial action of $G$ on $E$.

Consider  the automorphic crossed module $\A=(E=\Z_3,G=\Z_2,\t)$ such that $1 \t
a=a$ and $-1 \t a =-a$, where $a \in \Z_3$ and $\Z_2=(\{1,-1\},\times)${; see}
\cite{BM}. Then this crossed module detects the knottedness of the Virtual
Hopf Link $L$. If fact  $\H_\A(L)=6\#\{Y \in \Z_2; e \in \Z_3|  Y^{-1} \t e =e \}  =24,$
whereas $\H_\A(O^2)=36$.

\subsubsection{The Hopf Link}

The  Hopf Link $H$ is depicted in figure \ref{vhl}. Note that the fundamental
group of the complement of it is, similarly with the Virtual Hopf Link $L$,
isomorphic with $\Z^2$.

\begin{figure}
\centerline{\relabelbox 
\epsfysize 5cm
\epsfbox{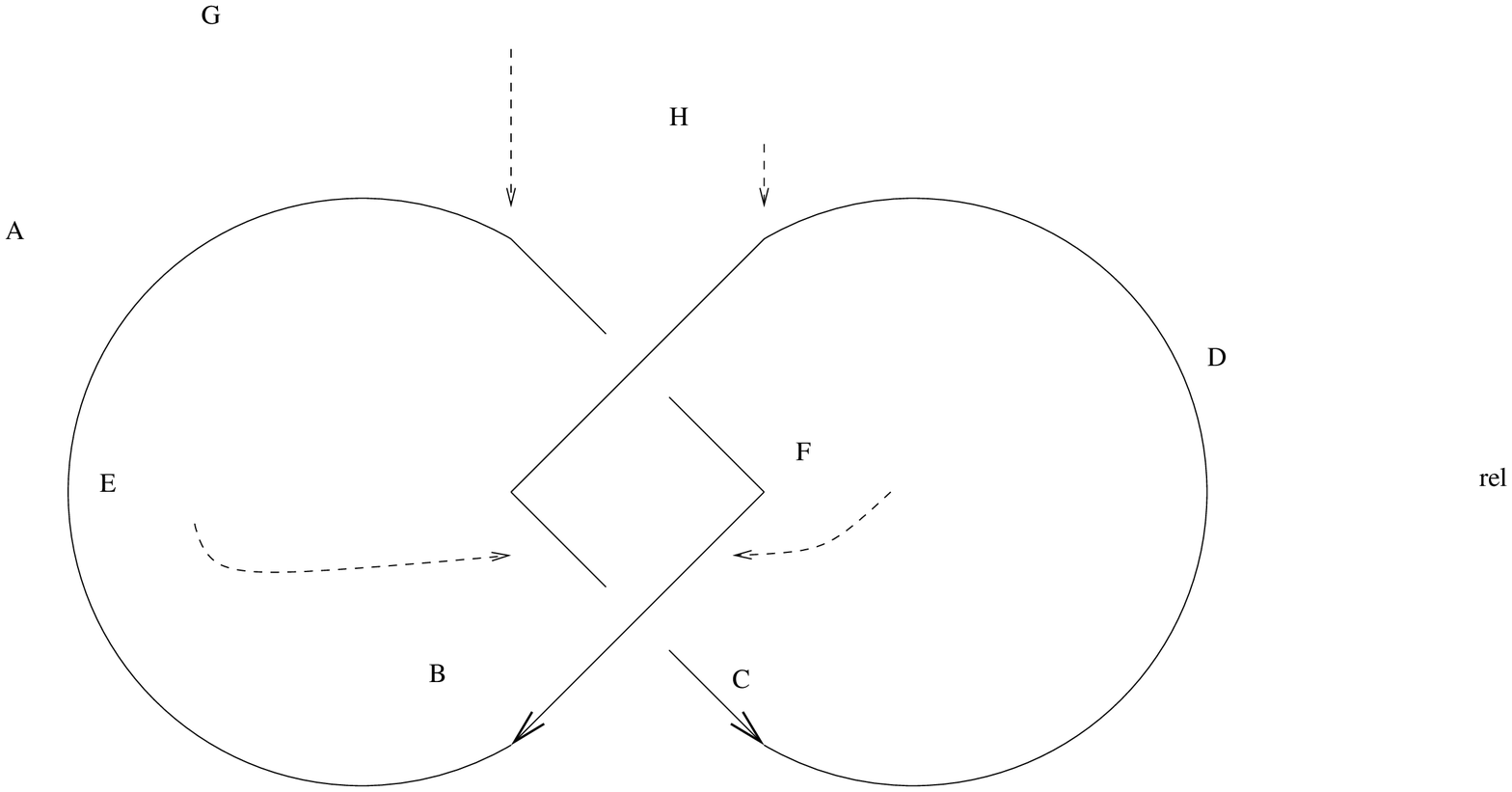}
\relabel {A}{$\s{(X,e^{-1})}$}
\relabel {B}{$\s{(X,e)}$}
\relabel {C}{$\s{(Y,f)}$}
\relabel {D}{$\s{(Y,f^{-1})}$}
\relabel {F}{$\s{(X,X^{-1}\t f^{-1} e f)}$}
\relabel {E}{$\s{(X^{-1}YX,X^{-1} \t f  )}$}
\relabel {G}{$\s{\left (X^{-1}Y^{-1}XYX,X^{-1}Y^{-1} \t f^{-1} X^{-1}Y^{-1}X \t (ef)\right)}$}
\relabel {H}{$\s{\left (X^{-1}YX, X^{-1} Y^{-1} \t f X^{-1}Y^{-1}X \t (e^{-1} f^{-1})
      ef\right )}$}
\relabel {rel}{$\left\{ \begin{CD}&  \s{X^{-1}Y^{-1}XYX=X} \\ &\s{X^{-1} YX=Y} 
\\ & \s{X^{-1}Y^{-1} \t f^{-1} X^{-1}Y^{-1}X \t (ef)=e }\\&  \s{X^{-1} Y^{-1} \t f X^{-1}Y^{-1}X \t (e^{-1} f^{-1})
      ef =f} \end{CD}\right. $}
 \endrelabelbox}
\caption{Calculation of the crossed module invariant of the Hopf Link.} 
\label{hlcalc}
\end{figure}

Let us calculate the 
crossed module invariant of the Hopf Link $H$. To this end, let $\G=(E,G,\t)$ be a finite
automorphic crossed module.  We display the calculation of $\H_\G(H)$ in
figure \ref{hlcalc}.
{This permits us to conclude that:}

$$\H_\G(H)=\# \left \{X,Y \in G; e,f \in E\left |\begin{CD} &XY=YX \\  & X^{-1}Y^{-1} \t f^{-1} X^{-1}Y^{-1}X \t(e f) =e\end{CD} \right. \right \} ,$$
which particularises to 
\begin{equation}
\H_\G(H)= \# \left \{X,Y \in G; e,f \in E:  X^{-1}Y^{-1} \t f^{-1} Y^{-1} \t(e f) =e \right \},
\end{equation}
in the case when $G$ is abelian. This is in agreement with the calculation in \cite{FM2}.

Let us see that the Hopf Link $H$ is not equivalent to the Virtual Hopf Link
$L$ as {a welded virtual link.} Consider the automorphic crossed module
$\A=(E=\Z^3,G=\Z^2, \t)$ defined above. We have (note that we switched to additive
notation, more adapted to this example):
\begin{align*}
\H_\A(H)&= \# \left \{X,Y \in \Z_2; e,f \in \Z_3: 
     -XY \t f+ Y \t (e +  f) =e \right \}\\
        &= \# \left \{X,Y \in \Z_2; e,f \in \Z_3: 
     -XY \t f+  Y \t f =e-Y \t e \right \}.
\end{align*}
In the case $Y=1$, we are led  to the equation $-X \t f+f=0$, which has
$4\times 3 $ solutions in $\Z_2 \times \Z_3 \times \Z_3$. In the case $Y=-1$,
we get the equation $e=2^{-1}( X \t f- f)$, which has $3\times 2$ solutions in
$\Z_2 \times \Z_3 \times \Z_3$.  Therefore, we obtain $\H_\A(H)=18$.

Therefore, we have proved that the Virtual Hopf Link is not equivalent to the
Hopf Link as {a welded virtual link,} and also {that the} Hopf Link is knotted, by using the crossed module invariant. 

As we have referred to before, the {knot groups of the} Hopf Link and the  Virtual Hopf Link are both isomorphic with
$\Z^2$. Therefore,  we {have} proved that the crossed module
invariant $\H_\G$ sees beyond the fundamental group of the complement of a
welded virtual knot.

Since the correspondence $K \mapsto T(K)$, where $K$ is a welded virtual link,
preserves the fundamental groups of the complement we have also proved:
\begin{Theorem}
The Crossed Module Invariant $I_\G$ of knotted surfaces defined in
\cite{FM1,FM2} is powerful enough to distinguish  {between} knotted surfaces  $\S,\S' \subset S^4$, with $\S$ diffeomorphic with $\S'$, whose complements have isomorphic fundamental groups, at least in a particular case.   
\end{Theorem}
Therefore, one of the main open problems about {the Crossed Module Invariant $I_\G$} of
knotted surfaces that prevails is whether the invariant $I_\G$ can
distinguish between knotted surfaces whose complements have isomorphic
fundamental groups and second {homotopy} groups, seen  as   $\pi_1$-modules,
but have  distinct Postnikov invariants $k^3 \in H^3(\pi_1,\pi_2)$. This problem was referred to  in
\cite{FM2}.  Examples of pairs of  knotted surfaces like this do exist; {see  \cite{PS}.}

\begin{Exercise}\label{E1}
 Consider the Hopf Arc $HA$ depicted in figure \ref{Trefoil}. Prove that $\H_\G(HA)=\H_\G(L)$, where $L$ is the Virtual Hopf Link.  Here $\G=(E,G,\t)$ is any finite automorphic crossed module. In fact, cf. \ref{ADD}, $T(L)$ is obtained from $T(HA)$ by adding a trivial 1-handle, which explains this identity. We  will  go back to this later in   \ref{FMK}. 
\end{Exercise}

\subsection{Trefoil Knot and Trefoil Arc}\label{TKTA}

\begin{figure}
\centerline{\relabelbox 
\epsfysize 1.5cm
\epsfbox{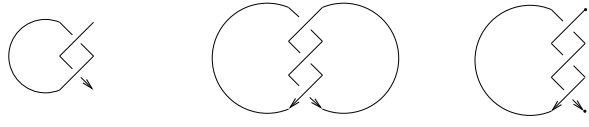}
\endrelabelbox}
\caption{The Hopf Arc $HA$, the Trefoil Knot ${{3_1}}$ and the  Trefoil Arc {${{3_1}'}$}.}
\label{Trefoil}
\end{figure}

The Trefoil Knot ${3_1}$ and the Trefoil Arc ${{3_1}'}$ are depicted in figure \ref{Trefoil}.

\begin{figure}
\centerline{\relabelbox 
\epsfysize 6.7cm
\epsfbox{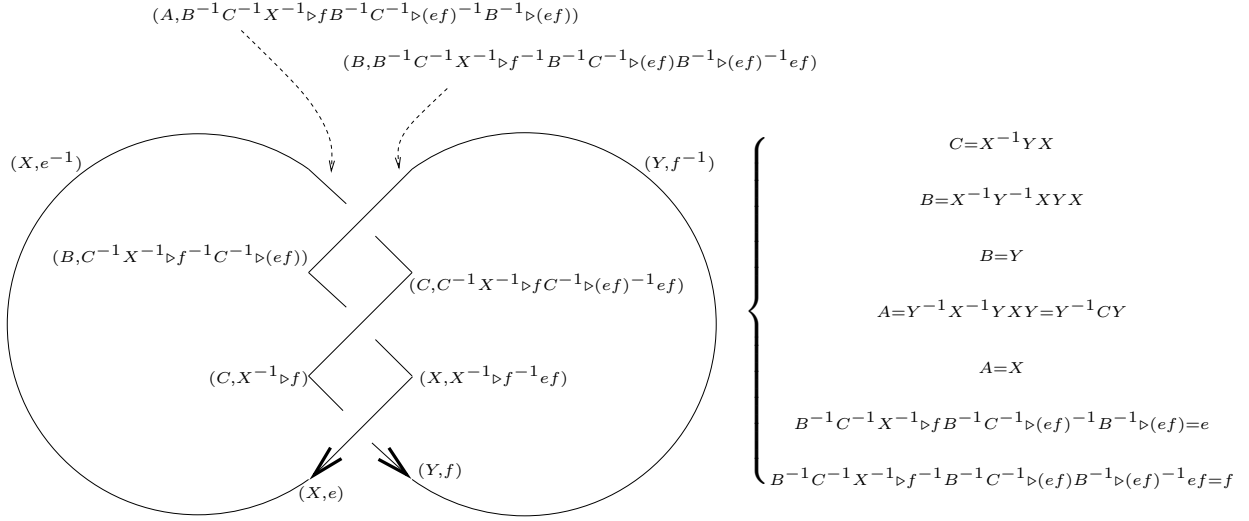}
\relabel {A}{$\ss{(X,e^{-1})}$}
\relabel {B}{$\ss{(X,e)}$}
\relabel {C}{$\ss{(Y,f)}$}
\relabel {D}{$\ss{(Y,f^{-1})}$}
\relabel {E}{$\ss{(X,X^{-1} \t f^{-1} e f )}$}
\relabel {F}{$\ss{(C, X^{-1}\t f)}$}
\relabel {G}{$\ss{(C,C^{-1} X^{-1} \t f C^{-1}{\t}(ef)^{-1} ef)}$}
\relabel {H}{$\ss{(B,C^{-1} X^{-1} \t f^{-1} C^{-1} \t (ef))}$}
\relabel {I}{$\ss{(A,B^{-1} C^{-1} X^{-1} \t f B^{-1}C^{-1} \t (ef)^{-1} B^{-1} \t (ef) )}$}
\relabel {J}{$\ss{(B,B^{-1} C^{-1} X^{-1} \t f^{-1} B^{-1}C^{-1} \t (ef) B^{-1} \t (ef)^{-1}ef )  }$}
\relabel {rel}{$\left\{ \begin{CD}& \ss{C=X^{-1}YX }\\ & \ss{B=X^{-1}Y^{-1}XYX }\\& \ss{B=Y} \\ &\ss{A=Y^{-1}X^{-1}YXY=Y^{-1}CY}\\&\ss{ A=X}\\ &\ss{B^{-1} C^{-1} X^{-1} \t f B^{-1}C^{-1} \t (ef)^{-1} B^{-1} \t (ef)  =e}\\ & \ss{B^{-1} C^{-1} X^{-1} \t f^{-1} B^{-1}C^{-1} \t (ef) B^{-1} \t (ef)^{-1}ef  =f  }\end{CD}\right. $}
 \endrelabelbox}
\caption{ \label{TrefoilCalc} {Calculation of the crossed module invariant of the Trefoil Knot ${3_1}$.}} 
\end{figure}

{Let us calculate the crossed module invariant of the Trefoil Knot ${3_1}$. Let $\G=(E,G,\t)$ be a finite automorphic crossed module. 
The calculation of $\H_\G({3_1})$ appears in figure \ref{TrefoilCalc}. }
This permits us to conclude that:
\begin{align}\label{exp1}
\H_\G({3_1})&={\# \left \{X,Y \in G; e,f \in E \left | \begin{CD} & \s{X^{-1} Y^{-1}X^{-1}=Y^{-1}X^{-1}Y^{-1}} \\&\s{Y^{-1}X^{-1}Y^{-1} \t fY^{-1}X^{-1}Y^{-1}X \t (ef)^{-1}Y^{-1} \t (ef)=e} \end{CD} \right. \right \}}\\   
&{=\# \left \{X,Y \in G; e,f \in E \left | \begin{CD} & \s{X^{-1} Y^{-1}X^{-1}=Y^{-1}X^{-1}Y^{-1}} \\&\s{Y^{-1}X^{-1}Y^{-1} \t fX^{-1}Y^{-1} \t (ef)^{-1}Y^{-1} \t (ef)=e} \end{CD} \right. \right \}.}\label{exp}
\end{align}
{This simplifies to:}
\begin{equation}\label{refer1}
{\H_\G({3_1})=\#\{X \in G; e,f \in E | X^{-3}\t f X^{-2} \t (ef)^{-1} X^{-1}
\t (ef)=e\}, }
\end{equation}
 when $\G=(E,G,\t)$ is an automorphic crossed module with $G$ abelian{; see \ref{Relation}.}

Note that the crossed module invariant of the {Trefoil Arc ${{3_1}'}$} can also be obtained from this calculation, by making $f=1_E$, and inserting the {necessary} {normalisation} factors{; see \ref{weldedarcs}.} This  yields:
\begin{equation}\label{exp3}
\H_\G({{{3_1}'}})=\#E\# \left \{X,Y \in G; e \in E \left | \begin{CD} & \s{{X^{-1} Y^{-1}X^{-1}=Y^{-1}X^{-1}Y^{-1}}} \\&\s{{X^{-1}Y^{-1} \t e^{-1}Y^{-1} \t e=e}} \end{CD} \right. \right \},
\end{equation}
{which simplifies to:}
\begin{equation}\label{refer2}
\H_\G({{3_1}'})=\#E \#\{X \in G; e \in E\colon  X^{-2} \t e^{-1} X^{-1} \t e=e\}, 
\end{equation}
whenever $G$ is abelian. This is coherent with the calculation in \cite{FM1,FM2}.

Observe that from equations  (\ref{refer1})  and (\ref{refer2}) it follows
that (we switch to additive notation):
\begin{align*}
{\H_\G({3_1})}&=\#\{X \in G; e,f \in E\colon  X^{-3}\t f -X^{-2} \t (e +f) + X^{-1} \t (e+f)=e\}\\
      &=\#\{X \in G; e,f \in E\colon X^{-2} \t \left (X^{-1} \t f-e \right)-X^{-1}  \t \left (X^{-1} \t f-e \right) +  \left (X^{-1} \t f-e \right) =0 \}\\
     &=\H_\G({{3_1}'}).
\end{align*}
Thus:
\begin{equation}\label{equal}
\H_\G({3_1})=\H_\G({{3_1}'}),
\end{equation}
 whenever $\G=(E,G,\t)$ is an automorphic crossed
module with $G$ abelian. {An analogous identity holds for any classical 1-component knot, {see \ref{Relation}.}}

{We will consider the crossed module invariants of the Trefoil Knot  and the Trefoil Arc for the case when $\G=(E,G,\t)$ is an automorphic crossed module with $G$ being a non-abelian group in \ref{nonab}. In this case the previous identity does not hold.}

Let us see that $\H_\G$ detects the knottedness of the Trefoil Knot ${3_1}$. The
crossed module $\A=(\Z_3,\Z_2,\t)$ defined previously  detects it. In fact it
is easy to see that $\H_\A({3_1})=12$. On the other  hand, if $O$ is the unknot,
we have that $\H_\G(O)=\#E \#G$,  for any automorphic crossed module
$\G=(E,G,\t)$.  Thus ${3_1}$ is knotted. Analogously we can prove that the {Trefoil Arc}
${{3_1}'}$ is knotted.

\begin{figure}
\centerline{\relabelbox 
\epsfysize 1.5cm
\epsfbox{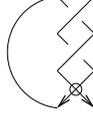}
\endrelabelbox}
\caption{A non trivial {welded virtual} arc whose closure is trivial.}
\label{VA}
\end{figure}

\begin{Exercise}\label{E2}
Consider the virtual arc $A$ of figure \ref{VA}. Prove that if $\G=(E,G,\t)$ is an automorphic finite crossed module {with $G$ abelian} then:
 $$\H_\G(A)=\#E\#\{X \in G;e \in E{|}X^{-2} \t e^{-1} X^{-1} \t e e^{-1}=1\}. $$
Thus the crossed module $\A=(\Z_3,\Z_2,\t)$ defined previously detects that it is knotted. However, it is easy to show that the closure of $A$ is {the} trivial welded virtual  knot, a fact confirmed by the crossed module invariant.  
\end{Exercise}

\subsection{{Universal module constructions}} \label{universal}

Let $G$ be an abelian group. Suppose that $\G=(E,G,\t)$ is an automorphic crossed module, where $E$ is an abelian group.
 Consider   a welded virtual link $K$.  Suppose that $K$ has $n$-components $S^1$, where $n$ is a positive integer. 
Let $\k_n=\Z[X_1,X_1^{-1},\ldots,X_n,X_n^{-1}]$ be the ring of Laurent polynomials on  the formal variables $X_1, \ldots X_n$. We can assign to $K$ a $\k_n$-module, so that $\H_\G(K)$ will satisfy:
$$\H_\G(K)=\# \Hom(\CM(K),\G),$$
where $\Hom(\CM(K),\G)$ denotes the set of all crossed module morphism $\CM(K) \to \G$. 

\subsubsection{The definition of the module $\CM(K)$}

\begin{figure}
\centerline{\relabelbox 
\epsfysize 7cm
\epsfbox{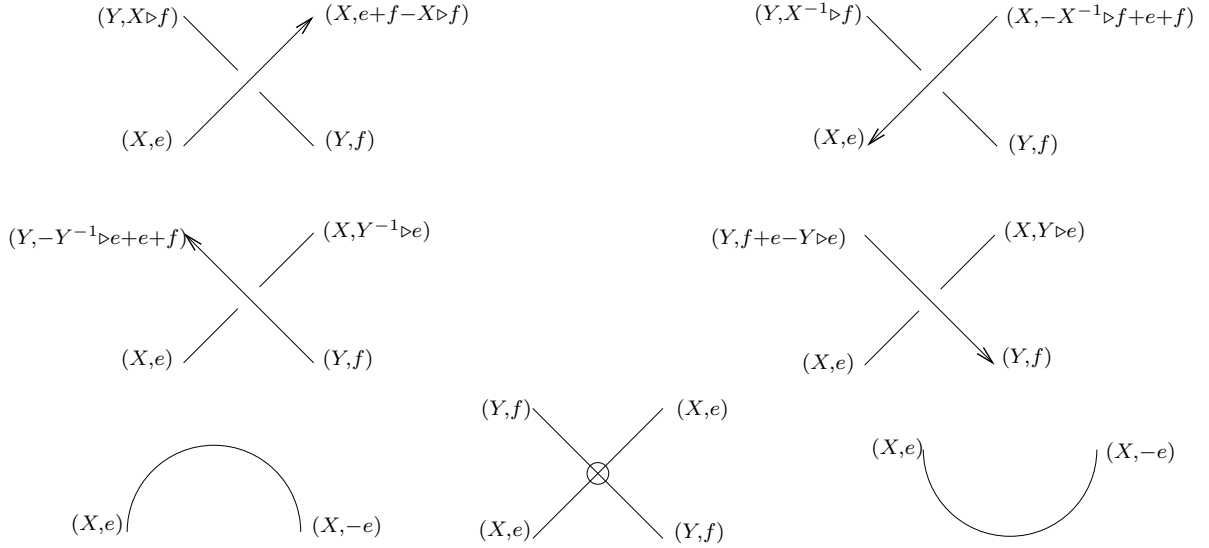}
\relabel{A}{$\s{(Y,X \t f)}$}
\relabel{B}{$\s{(X,e+f-X \t f)}$}
\relabel{C}{$\s{(X,e)}$}
\relabel{D}{$\s{(Y,f)}$}
\relabel{E}{$\s{(Y, X^{-1} \t f) }$}
\relabel{F}{$\s{(X,-X^{-1} \t f+ e +f )}$}
\relabel{G}{$\s{(X,e)}$}
\relabel{H}{$\s{(Y,f)}$}
\relabel{I}{$\s{(Y, -Y^{-1} \t e+ e+ f)}$}
\relabel{J}{$\s{(X,Y^{-1} \t e) }$}
\relabel{K}{$\s{(X,e)}$}
\relabel{L}{$\s{(Y,f)}$}
\relabel{M}{$\s{(Y,f+e-Y \t e)}$}
\relabel{N}{$\s{(X, Y \t e)}$}
\relabel{O}{$\s{(X,e)}$}
\relabel{P}{$\s{(Y,f)}$}
\relabel{Q}{$\s{(Y,f)}$}
\relabel{R}{$\s{(X,e)}$}
\relabel{S}{$\s{(X,e)}$}
\relabel{T}{$\s{(Y,f)}$}
\relabel{X}{$\s{(X,e)}$}
\relabel{Y}{$\s{(X,-e)}$}
\relabel{Z}{$\s{(X,e)}$}
\relabel{W}{$\s{(X,-e)}$}
\endrelabelbox}
\caption{Defining relations for the module $\CM(K)$.}
\label{CMrel}
\end{figure}

\begin{Definition}
Let $K$ be a welded virtual link diagram. Suppose that $K$ is an immersion of a disjoint union of $n$ circles $S^1$ into the plane, each of which is assigned a variable $X_i$, where $i\in \{1,\ldots,n\}$; in other words, suppose that we have a total order on the set of all $S^1$-components of $K$. The module $\CM(K)$ is defined as the $\k_n$-module generated by all the connected components of $K$ minus the set of crossings of $K$  and extreme points of $K$, modding out by the relations of figure \ref{CMrel}. It is understood that any connected component is assigned a pair $(X,e)$, where $e \in \CM(K)$ is the module element that  the connected components defines, whereas $X \in  \{X_1,\ldots,X_n\}$ is the labelling of the $S^1$-component of $K$ in which the connected component is included. 
\end{Definition}

By using the same technique as in Exercise \ref{Ex1} we can prove:
\begin{Theorem}
Let $K$ be a  welded virtual link diagram with $n$ $S^1$-components. The isomorphism class of the $\k_n$-module $\CM(K)$ depends only on the welded virtual {link} determined by $K$, up to  reordering of the $S^1$-components  of $K$. In addition, if $\G=(E,G,\t)$ is an automorphic {finite} crossed module with $G$ abelian we have:
$$\H_\G(K)=\# \Hom(\CM(K),\G).$$
\end{Theorem}

\subsubsection{Relation with the Alexander Module}\label{Relation}

\begin{figure}
\centerline{\relabelbox 
\epsfysize 2.5cm
\epsfbox{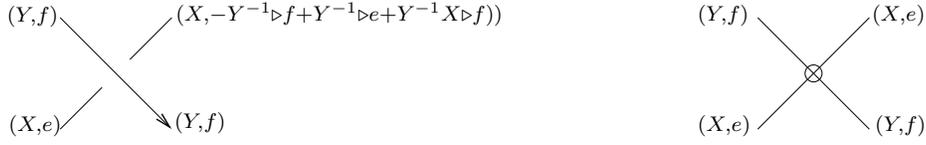}
\relabel{A}{$\s{(X,e)}$}
\relabel{B}{$\s{(Y,f)}$}
\relabel{D}{$\s{(Y,f)}$}
\relabel{C}{$\s{(X,-Y^{-1} \t f+Y^{-1} \t e + Y^{-1}X\t f))}$}
\relabel{X}{$\s{(Y,f)}$}
\relabel{Y}{$\s{(X,e)}$}
\relabel{Z}{$\s{(X,e)}$}
\relabel{W}{$\s{(Y,f)}$}
\endrelabelbox}
\caption{Relations at crossings for the Alexander Module $\Alex(K)$.}
\label{Alex}
\end{figure}

Let $K$ be a welded virtual link diagram with $n$ $S^1$-components, each labelled with an $X_i \in  \{X_1,...,X_n\}$. We can define the Alexander module $\Alex(K)$
 of $K$, defined as the module over $\k_n$ with a generator {for}  each  connected component of $K$ minus its set of crossings, modulo the relations of figure \ref{Alex}, obtained from the right handed Wirtinger relations of figure \ref{Wirtinger} by applying Fox derivatives; see \cite[Chapter 9]{BZ}, \cite[Chapter XI]{K2} or \cite{F}.  Therefore, if  $K$ is a classical 1-component knot, then $\Alex(K)\cong \Z[X,X^{-1}]/\left <\Delta(K)=0\right >\oplus \Z[X,X^{-1}]$, where $\Delta(K)$ denotes the Alexander polynomial of $K$; see for example \cite[9 C]{BZ}.

Let $K$ be a welded virtual link diagram. The Alexander module  $\Alex(K)$ depends only on the knot group of the  welded virtual link  defined by $K$, up to isomorphism and reordering of the $S^1$-components of $K$. 

\begin{figure}
\centerline{\relabelbox 
\epsfysize 2.5cm
\epsfbox{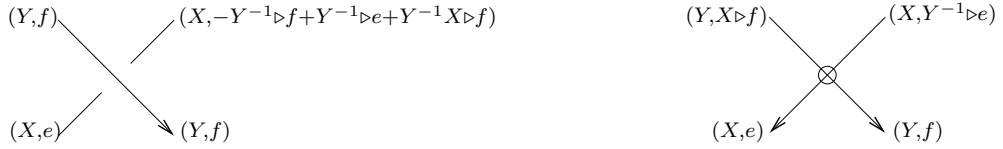}
\relabel{A}{$\s{(X,e)}$}
\relabel{B}{$\s{(Y,f)}$}
\relabel{D}{$\s{(Y,f)}$}
\relabel{C}{$\s{(X,-Y^{-1} \t f+Y^{-1} \t e + Y^{-1}X\t f)}$}
\relabel{X}{$\s{(Y,{X \t f})}$}
\relabel{Y}{$\s{(X,e)}$}
\relabel{Z}{$\s{(X,{Y^{-1} \t e)}}$}
\relabel{W}{$\s{(Y,f)}$}
\endrelabelbox}
\caption{Relations at crossings for the module $\Alex'(K)$.}
\label{Alexvar}
\end{figure}

The module $\Alex(K)$ admits a variant $\Alex'(K)$ whose defining relations appear in figure \ref{Alexvar}. Note that the $\k_n$-module $\Alex(K)$ is isomorphic to $\Alex'(K)$ whenever $K$ is a classical link diagram.

The  module $\Alex'(K)$ is invariant under virtual and classical Reidemeister moves. However, $\Alex'(K)$ is not invariant under the first forbidden move $F_1$; rather it is invariant under the second forbidden move $F_2$; see subsection 
\ref{r1}.

 Given a virtual link diagram $K$, we can define the mirror image $K^*$ of it by switching positive to negative crossings, and vice-versa, and leaving virtual crossings unchanged.  Therefore, the module $\Alex'(K^*)$ depends only on the welded virtual knot defined by $K$, up to isomorphism and reordering of the components of $K$.  

\begin{Theorem}\label{relation}
Let $K$ be a welded virtual link diagram.  There exists an isomorphism $$\f\colon \CM(K) \to \Alex'(K^*).$$
\end{Theorem}

\begin{proof}
We can suppose that $K$ is the closure of a virtual braid $B$; see \cite{KL,Ka}. This avoids needing to deal with the {defining} relations of $\CM(K)$ at maximal and minimal points. Let $b$ be a connected component of the braid $B$ minus its set of crossings,  defining therefore an element $b \in \CM(K)$. The isomorphism $\f\colon \CM(K) \to \Alex'(K^*)$ sends $b$ to $Z^{-1} \t b$, where $Z$ is the product of all the elements $X_i$ assigned to the strands of $B$ on the left of $b$ (each belonging to a certain $S^1$-component of $K$). The remaining details are left to the reader.
\end{proof}

The Alexander module of the Trefoil Knot {$3_1$} is the module over $Z[X,X^{-1}]$ with generators $e$ and $f$ and the relation $X^2\t (e+f)-X \t ( e+f) +(e+f)=0$, thus  we have $\Alex({3_1})=\Z[X,X^{-1}]/\left <X^2-X+1=0\right > \oplus \Z[X,X^{-1}]$. In particular, it follows equation (\ref{refer1}). 

Let $K$ be a classical {1-component knot.} By using Theorem \ref{relation}, we can prove  that for any automorphic crossed module $\G=\left ( E, G, \t\right )$, with $G$ abelian, the invariant $\H_\G(K)$ is determined {by} the Alexander module $\Alex(K)$ of $K$, and thus from the  Alexander polynomial $\Delta(K)$ of $K$.  This is not the case for non classical links, since the crossed module invariants of {the}  virtual and classical Hopf links $L$ and $H$; see {subsection \ref{vchl}} are different,
 {even though} they have isomorphic Alexander modules. In fact we have:
 $$\Alex(H),\CM(H), \Alex(L)=\frac{\Z[X,X^{-1},Y,Y^{-1}] \t e  \oplus \Z[X,X^{-1},Y,Y^{-1}]  \t f} {\left < (X-1) \t f=(Y -1) \t e\right >},$$
the module over the ring $Z[X,X^{-1},Y,Y^{-1}]$ with two generators $e$ and $f$, and the relation $(X-1) \t f=(Y -1) \t e$, whereas
$$\CM(L)=\frac{\Z[X,X^{-1},Y,Y^{-1}] \t e  \oplus \Z[X,X^{-1},Y,Y^{-1}]  \t f} {\left < Y \t f=f\right > }.$$ 
These last two modules are not isomorphic, as the calculations in {subsection \ref{vchl}} certify.

\subsubsection{Welded virtual arcs}

Let $A$ be a welded virtual arc with a single component. The {$\Z[X,X^{-1}]$-modules}  $\Alex(A),\Alex'(A)$ and $\CM(A)$ defined above can still be assigned to $A$, considering the analogue of the relations in figure  \ref{endscoloured} at the end-points of $A$, so that the elements of $\Alex(A),\Alex'(A)$ and $\CM(A)$ assigned to the edges of $A$ incident to {its end-points} are zero.

Any welded virtual  arc $A$ can be obtained as the (incomplete) closure of some braid. Therefore the proof of Theorem \ref{relation} gives an isomorphism $\f\colon \CM(A) \to \Alex'(A)$. 

Suppose that $A$ is a classical arc  sitting in the semiplane
$\{z\ge 0\}$ of $\R^3$, intersecting the plane $\{z=0\}$ at the end-points of
$A$, only. Since $A$ is classical we have $\Alex(A)=\Alex'(A)$. Let $K$ be the obvious closure of $A$.  Then $\Alex(K)=\Z[X,X^{-1}]/\left <\Delta(K)=0\right >\oplus \Z[X,X^{-1}]$, where $\Delta(K)$ is the Alexander polynomial of $K$. Choosing a  connected component of $K$ minus its set of crossings, and sending the generator of $\Alex(K)$ it defines to zero yields a presentation of $\Z[X,X^{-1}]/\left <\Delta(K)=0\right >$; see \cite[Theorem 9.10]{BZ}. Comparing with the definition {of} $\Alex(A)$, proves that $\Alex(A)=\Z[X,X^{-1}]/\left <\Delta(K)=0\right >$.

Therefore it follows that $\CM(A)\cong Z[X,X^{-1}]/\left <\Delta(K)=0\right >$ if $A$ is a classical arc and $K$ is the closure of $A$. The discussion above also implies that  if $\G=(E,G,\t)$ is an automorphic crossed module with $G$ abelian then  {$\H_\G(K)=\H_\G(A)$} whenever $A$ is a classical 1-component arc and $K$ is the 1-component knot obtained by closing $A$. This is not the case if $G$ is not abelian.

\begin{figure}
\centerline{\relabelbox 
\epsfysize 4cm
\epsfbox{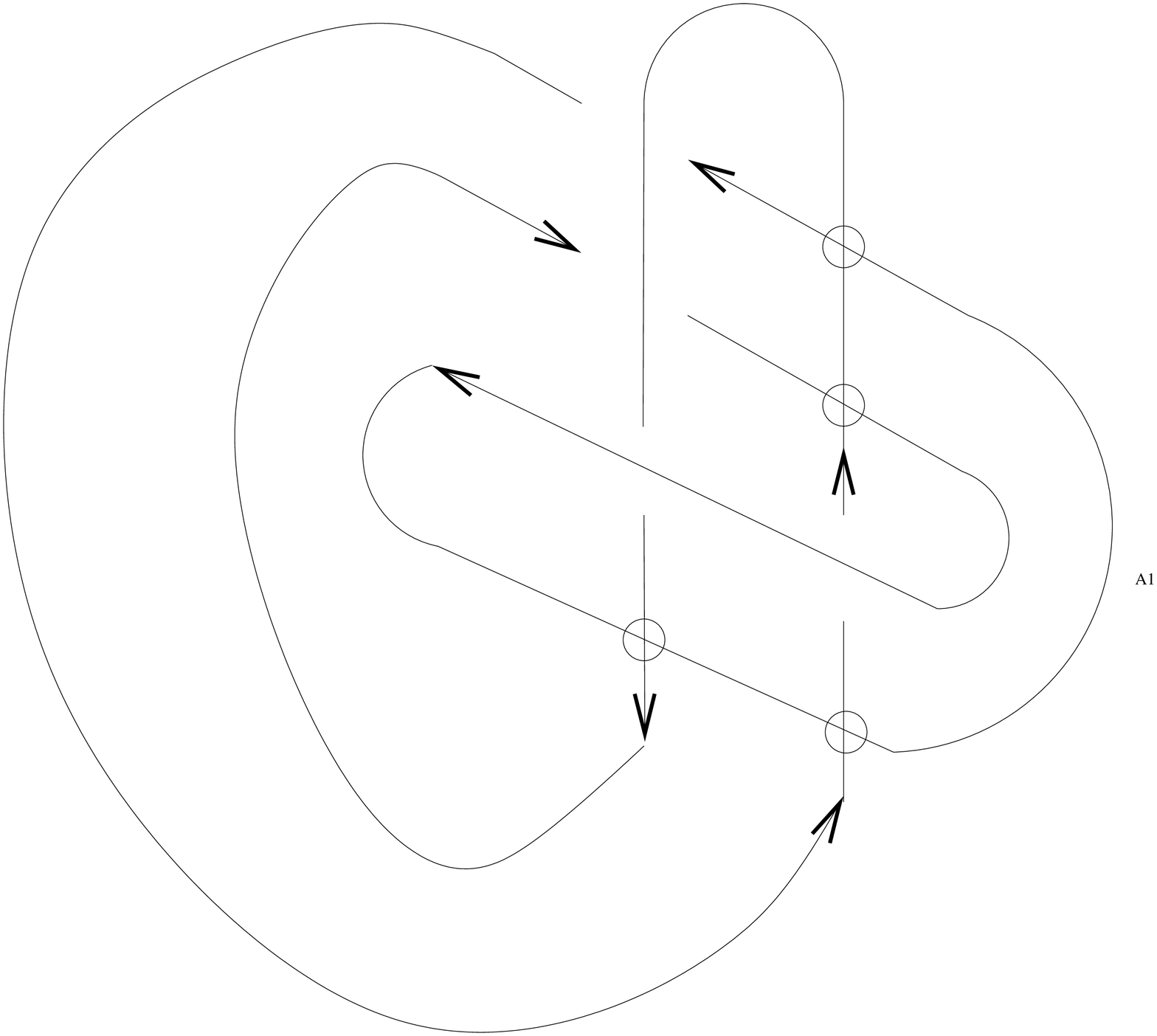}
 \endrelabelbox}
\caption{ \label{satohexample} Shin Satoh's Knot {$S$}.} 
\end{figure}
\begin{Problem}
{Let $K$ be a welded virtual link. What is the algebraic topology interpretation of the module $\CM(K)$ in terms of the tube $T(K)\subset S^4$ of $K$.}
\end{Problem}

\subsection{Shin Satoh's Knot}\label{ss}
{In \cite{S}, Shin Satoh considered the welded virtual link $S$ displayed in figure \ref{satohexample}. It is a welded virtual knot whose knot group is isomorphic with the knot group of the
Trefoil Knot $3_1$.}  It is possible to prove that $S$ is not equivalent to any
classical knot as a welded virtual knot, see \cite{S}, thus the Shin Satoh's Knot $S$ is not equivalent to the Trefoil.  {See also \ref{nonab}.}

Let us calculate the crossed module invariant of the Shin Satoh's Knot $S$. Let $\G=(E,G,\t)$ be a finite automorphic crossed module. {We consider in this case that $G$ is an abelian group, which makes the calculations much easier, since we simply need to calculate the $\Z[X,X^{-1}]$-module $\CM(S)$. {The case when $G$ is non-abelian is considered  in \ref{nonab}.}  }
Figure \ref{SatohCalc} permits us to conclude that:

\begin{figure}
\centerline{\relabelbox 
\epsfysize 9cm
\epsfbox{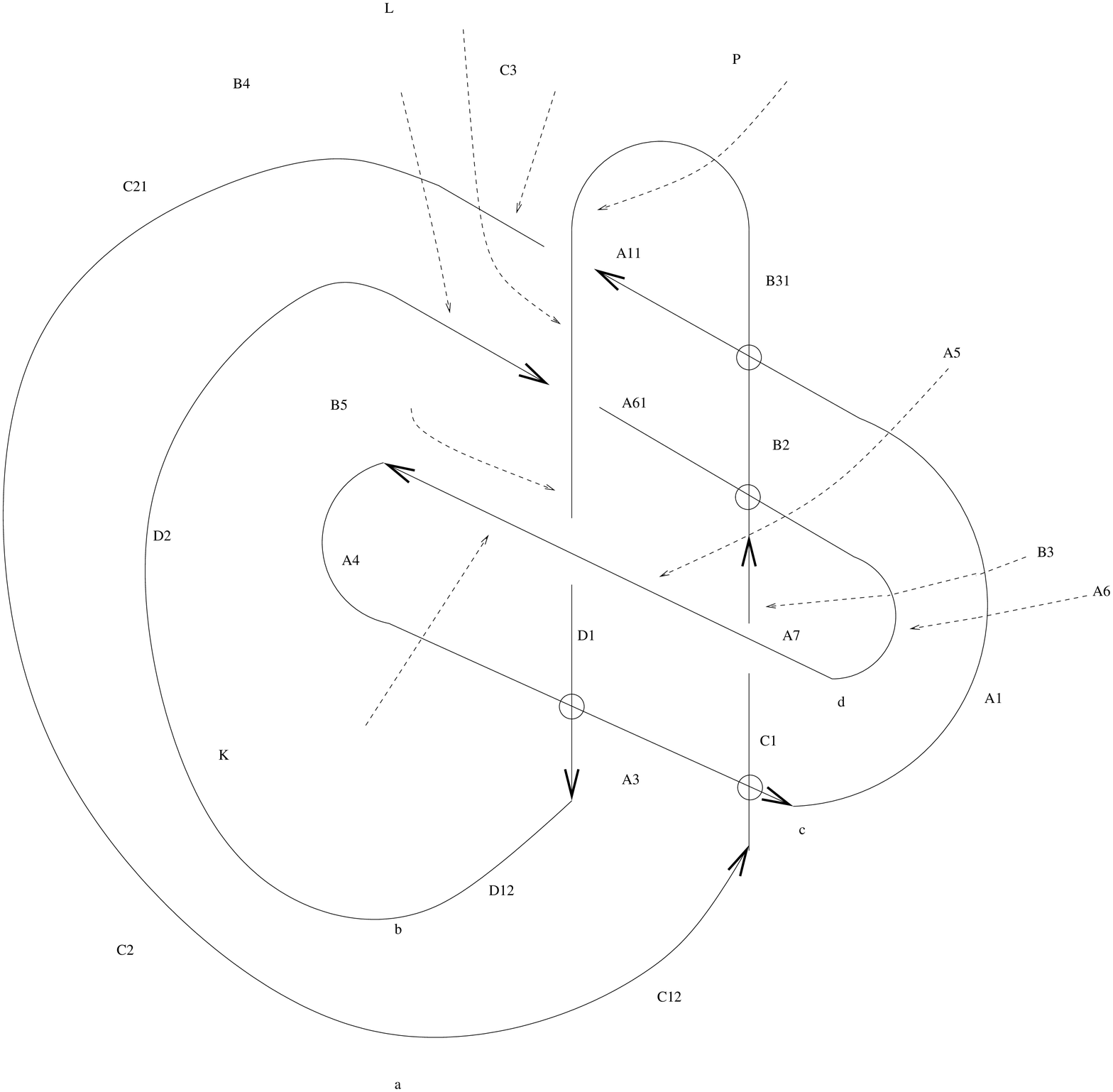}
\relabel {A1}{$\ss{\left ( X,c^{-1}     \right) }$}
\relabel {A11}{$\ss{\left ( X,c^{-1}     \right) }$}
\relabel {C2}{$\ss{\left ( X,a     \right) }$}
\relabel {C21}{$\ss{\left ( X,a     \right) }$}
\relabel {C1}{$\ss{\left ( X,a^{-1}      \right) }$}
\relabel {C12}{$\ss{\left ( X,a^{-1}      \right) }$}
\relabel {D2}{$\ss{\left ( X,b                \right) }$}
\relabel {A3}{$\ss{\left ( X,c                \right) }$}
\relabel {D1}{$\ss{\left ( X,b^{-1}           \right) }$}
\relabel {D12}{$\ss{\left ( X,b^{-1}           \right) }$}
\relabel {A4}{$\ss{\left ( X,c                \right) }$}
\relabel {A7}{$\ss{\left ( X,d                \right) }$}
\relabel {A6}{$\ss{\left ( X,d^{-1}           \right) }$}
\relabel {A61}{$\ss{\left ( X,d^{-1}           \right) }$}
\relabel {B3}{$\ss{\left ( X,X^{-1} \t a^{-1}         \right) }$}
\relabel {B31}{$\ss{\left ( X,X^{-1} \t a^{-1}         \right) }$}
\relabel {A5}{$\ss{\left ( X,X^{-1} \t a a^{-1} d    \right) }$}
\relabel {B5}{$\ss{\left ( X,X^{-1} \t b^{-1}   \right) }$}
\relabel {K}{$\ss{\left ( X,X^{-1} \t (ab) a^{-1} d b^{-1}  \right)}$}
\relabel {B4}{$\ss{\left ( X,X^{-1} \t d^{-1}    \right) }$}
\relabel {L}{$\ss{\left ( X,X^{-1} \t d X^{-1} \t b^{-1} d^{-1}    \right) }$}
\relabel {C3}{$\ss{\left ( X,X ^{-1} \t c^{-1}  \right) }$}
\relabel {P}{$\ss{\left ( X,X ^{-1} \t c  c^{-1}X^{-1} \t d X^{-1} \t b^{-1} d^{-1} \right) }$}
\endrelabelbox}
\caption{ \label{SatohCalc} Calculation of the crossed module invariant  of
  the Shin Satoh's Knot $S$ for $G$ abelian.} 
\end{figure}
\begin{align*}
\H_\G(S)&=\# \left \{X \in G; a,b,c,d \in E\left | \begin{CD} & X^{-1} \t (ab) a^{-1} d b^{-1}=c^{-1}  \\ &X^{-1} \t d^ {-1}=b^{-1} \\ & X^{-1} \t c^{-1} =a^{-1}\\ & X ^{-1} \t c  c^{-1}X^{-1} \t d X^{-1} \t b^{-1} d^{-1} =X^{-1} \t a \end{CD} \right .
\right\} \\ & \\
&=\# \left \{X \in G; a,d \in E \left | \begin{CD} & X^{-1} \t a X^{-2} \t d  a^{-1} d X^{-1} \t d^{-1}=X \t a^{-1}  \\ & a X \t a^{-1} X^{-1} \t d X^{-2} \t d^{-1} d^{-1} =X^{-1} \t a \end{CD}\right .
\right\}.
\end{align*}
The two equations in the final expression are equivalent. We obtain, switching to additive notation:
 \begin{equation}
\H_\G(S)=\#\left \{X \in G; a,d \in E| X^{-1} \t a -a+X \t a=X^{-1} \t d -d-X^{-2}\t d\right \}.
\end{equation}
This should be compared with the crossed module invariant of the Trefoil Knot ${3_1}$,
for $G$ abelian:
\begin{align*}
\H_\G({3_1})&=\#\{X \in G; e,f \in E| X^{-3}\t f - X^{-2} \t (e+f)+ X^{-1}
\t (e+f)=e\}\\
&=\#\{X \in G; e,f \in E| X^{-3}\t f -X^{-2} \t f + X^{-1} \t f=e-X^{-1}
\t e +X^{-2}\t e\}\\
&=\#\{X \in G; e,f \in E |  X^{-2}\t f -X^{-1} \t f +   f=X \t e- e +X^{-1}\t e\}.
\end{align*}
Therefore it follows that if $\G=({E,G},\t)$ is an automorphic crossed module
with $G$ abelian then:
\begin{equation}\label{KS}
\H_\G({3_1})=\H_\G(S).
\end{equation}
{We present in the following subsection {(see \ref{ADD})} an alternative  proof of this fact, which should reassure the
reader that the calculations in this article are correct, despite this being somehow a negative example.
{We will also see below (see \ref{nonab}) that if we take $G$ to be non-abelian, then we can prove that the Trefoil Knot is not equivalent to the Shin Satoh's Knot,  by using the crossed module invariant.} }

\subsection{Welded Virtual Graphs}
\subsubsection{Crossed module invariants of knotted surfaces obtained by
  adding trivial 1-handles}\label{ADD}

\begin{figure}\label{main}
\centerline{\relabelbox 
\epsfysize 4cm
\epsfbox{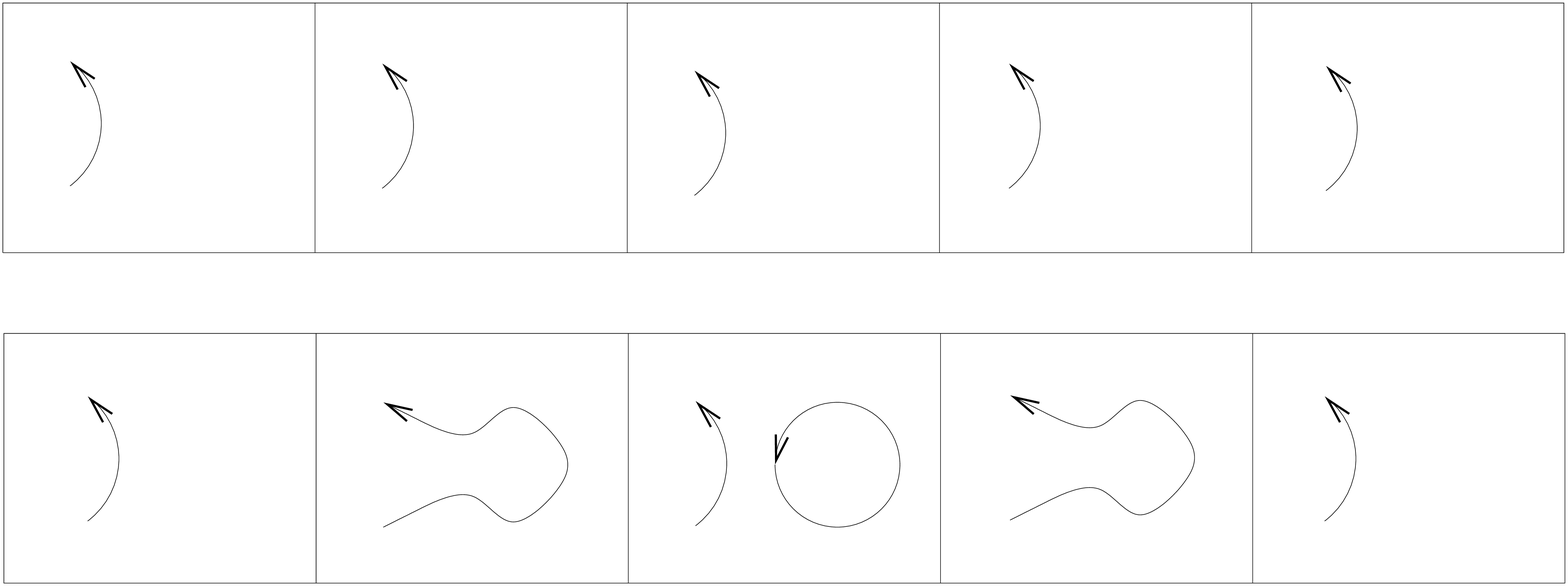}
\endrelabelbox}
\caption{Adding a trivial 1-handle to a knotted surface. On the top  we display the original
  movie, and on the bottom the new movie, both read from left to right. A concise description of this modification is fission saddle, fusion saddle.}
\label{add}
\end{figure}

Let $\S\subset S^4$ be a knotted surface which we suppose to be connected. The knotted surface $\S'$ obtained from $\S$ by
adding a trivial 1-handle is defined simply as the connected sum $\S'=\S \# T^2$, where $T^2$
is a torus $S^1 \times S^1$, trivially embedded in $S^4$. The {non-connected} case is totally analogous, but a connected component of $\S$ must be chosen.
 A movie of $\S'$ is
obtained from a movie of $\S$ by choosing a strand of the  movie of $\S$ belonging to the chosen component of $\S$,  and making the
modification shown in figure \ref{add}. The straightforward proof of the following theorem is
left to the reader. 
\begin{Theorem}\label{addtheo}
Let ${\G}=\left( E \ra{\d} G,\t \right)$ be a finite  crossed module. If the {oriented}
knotted surface $\S'$ is
obtained from the oriented knotted surface $\S$ by adding a trivial 1-handle then:
$$I_\G(\S')=\frac{(\# \ker \d)^2}{(\# E)^2}I_\G(\S),$$
thus in particular $I_\G(\S)=I_\G(\S')$ whenever $\G$ is automorphic. 
\end{Theorem}

The tube $T(S)$ of the Shin Satoh's Knot $S$ is obtained from the Spun Trefoil (the tube $T({3_1'})$ of the Trefoil Arc ${3_1'}$) by adding a trivial 1-handle; see \cite{S} or \ref{FMK}. This fact
together with equation (\ref{equal}) proves that {$\H_\G({3_1})=\H_\G(S)$,} whenever 
$\G=(E,G,\t)$ is a finite automorphic crossed module with $G$ abelian, as  already proved by
other means{; see subsection \ref{ss}. Here {$3_1$} is the Trefoil Knot.}

\subsubsection{Definition of welded virtual graphs}\label{FMK}
\begin{figure}
\centerline{\relabelbox 
\epsfysize 2cm
\epsfbox{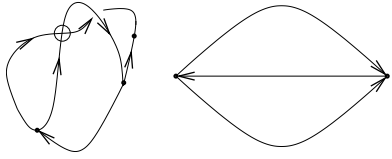}
\endrelabelbox}
\caption{A welded virtual graph.}
\label{wvgraph}
\end{figure}

Let $K$ be an oriented virtual graph diagram.  Note that $K$ may have some bivalent vertices where the orientation of an edge of $K$ may change; however, there cannot be a change of orientation {of a strand} at a crossing{; see} figure \ref{wvgraph}. 

\begin{figure}
\centerline{\relabelbox 
\epsfysize 2cm
\epsfbox{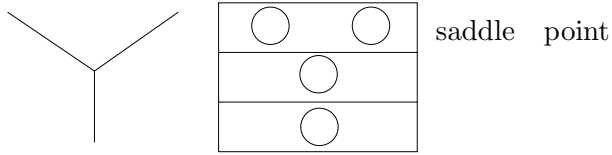}
\relabel{a}{$\textrm{saddle}\quad\textrm{point}$}
\endrelabelbox}
\caption{The tube of a  virtual graph at a 3-valent vertex {(movie version). As usual, all circles are oriented counterclockwise.}}
\label{vertex}
\end{figure}

\begin{figure}
\centerline{\relabelbox 
\epsfysize 2cm
\epsfbox{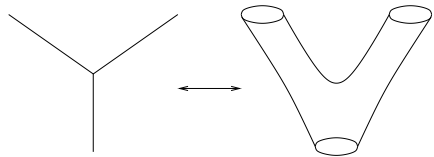}
\endrelabelbox}
\caption{The tube of a  virtual graph at a 3-valent vertex; broken surface {diagram} version of the movie of figure \ref{vertex}.}
\label{threetube}
\end{figure}

Given a virtual graph diagram $K$, we can define the tube $T(K)$ of it exactly in the same way as the tube of a virtual link or arc is defined. {We consider the type of  movie of figure \ref{vertex} at the 3-valent vertices. For  the   broken surface  diagram version of this see figure  \ref{threetube}. We  proceed analogously for $n$-valent vertices if $n>3$.}  The 2-valent vertices do not affect the calculation of $T(K)$. On the other hand $1$-valent vertices were already considered in the case of virtual arcs. 

\begin{figure}
\centerline{\relabelbox 
\epsfysize 9cm
\epsfbox{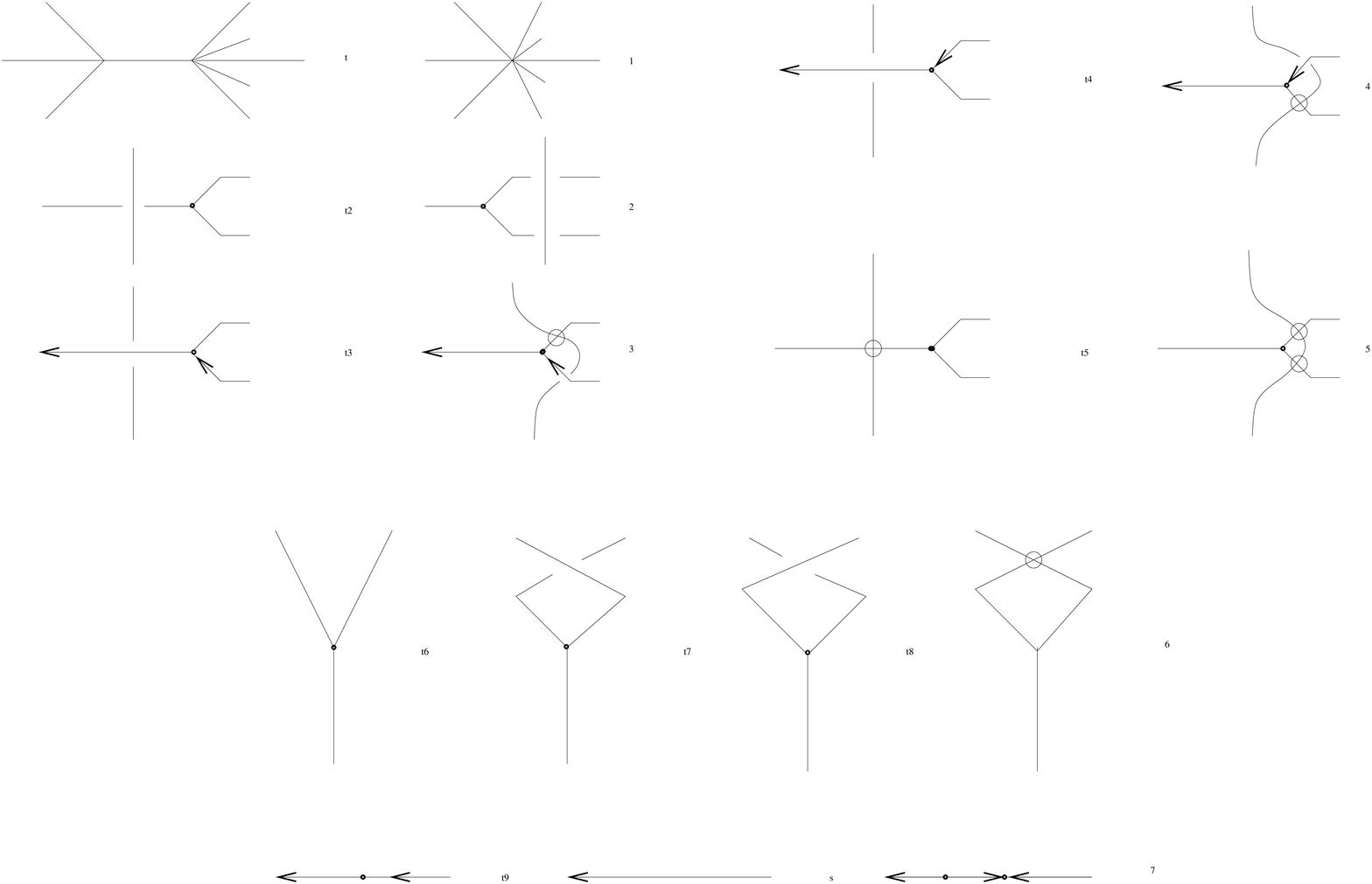}
\relabel{t}{$\leftrightarrow$}
\relabel{t2}{$\leftrightarrow$}
\relabel{t3}{$\leftrightarrow$}
\relabel{t4}{$\leftrightarrow$}
\relabel{t5}{$\leftrightarrow$}
\relabel{t6}{$\leftrightarrow$}
\relabel{t7}{$\leftrightarrow$}
\relabel{t8}{$\leftrightarrow$}
\relabel{t9}{$\leftrightarrow$}
\relabel{1}{$1$}
\relabel{2}{$2$}
\relabel{3}{$3$}
\relabel{4}{$4$}
\relabel{5}{$5$}
\relabel{6}{$6$}
\relabel{7}{$7$}
\relabel{s}{$\leftrightarrow$}
\endrelabelbox}
\caption{Moves defining Welded Virtual Graphs. Notice that the third and forth {moves} have a variant for  which the direction of each strand is reversed. However, these moves are a consequence of the remaining.}
\label{WVG}
\end{figure}

{
It is easy to see that the tube $T(K)$  of a virtual graph is invariant under the moves defining welded virtual knots and arcs{; see subsection} \ref{r1} and \ref{refer}. In addition, $T(K)$ is invariant under the moves shown in figure \ref{WVG}. } {Note that if a strand in figure \ref{WVG} is drawn without orientation, then this means that the corresponding identity is valid for any choice of orientation.}

\begin{figure}
\centerline{\relabelbox 
\epsfysize 4cm
\epsfbox{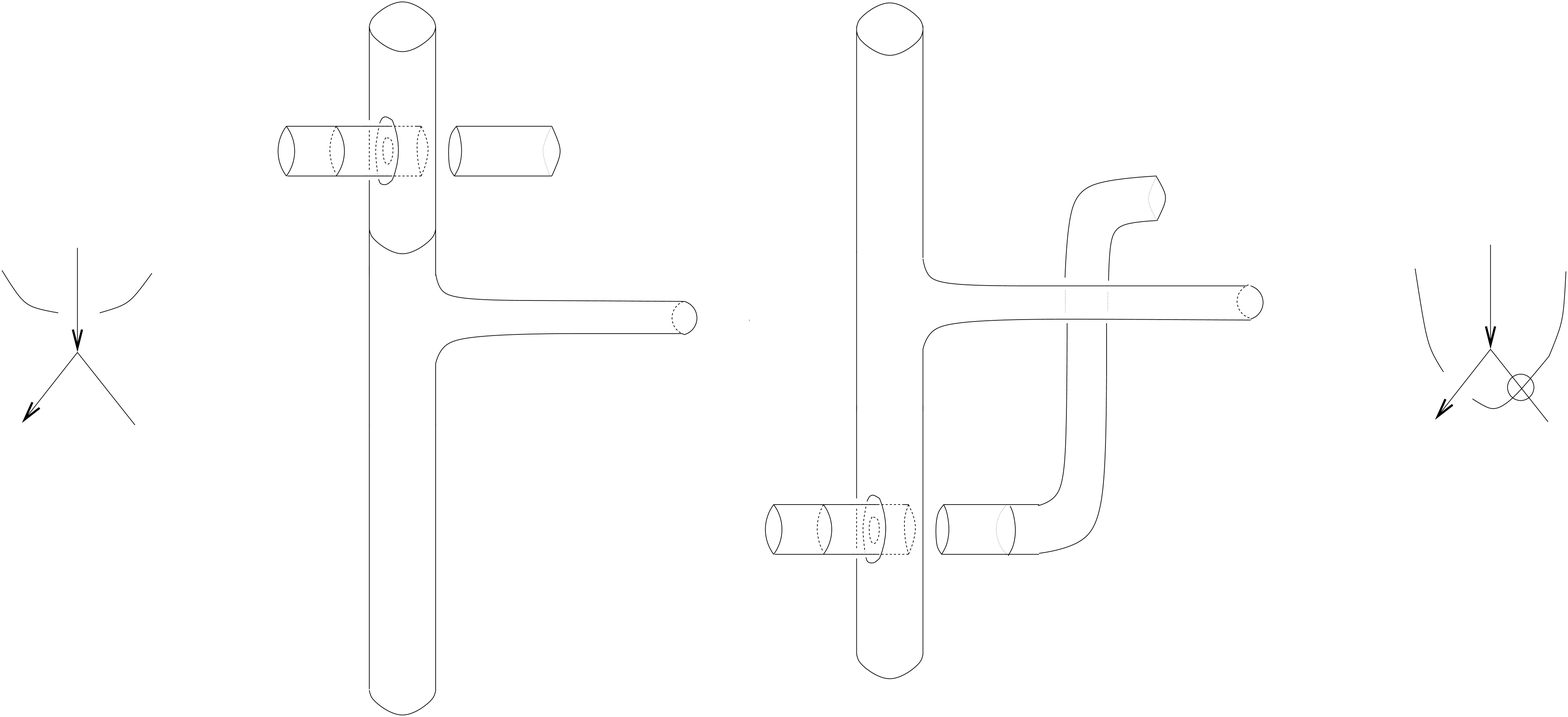}
\relabel{t}{$=$}
\endrelabelbox}
\caption{{An identity between {broken surface diagrams of knotted surfaces} (reverse orientation of the fourth move of figure \ref{WVG}.)}}
\label{slide1}
\end{figure}

{The invariance under the first, second and fifth moves is immediate.}
{ The invariance under the third and  forth  moves follows from figures \ref{ShinTube} and \ref{threetube},  by sliding the cylinder that goes inside the other cylinder towards the  end strand, in the obvious way, as shown in figure \ref{slide1}. It is strictly necessary that the edges incident  to the vertex in cause have compatible orientations in the sense shown in figure \ref{WVG}. Note that otherwise the crossing informations in the corresponding {initial} and {final} broken surface diagrams in figure \ref{slide1} would not be {compatible}.}

{
 The invariance of $T(K)$  under the  penultimate  moves  of figure \ref{WVG} follows from the same argument that proves invariance under the classical and virtual Reidemeister-I {moves.} }

\begin{Definition}[Welded Virtual Graph]
The  moves on oriented {virtual} graph diagrams  of  figure \ref{WVG}, together with the ones defining {welded} virtual knots and {welded virtual} arcs define what we called  a ``welded virtual graph''. 
\end{Definition}

\begin{figure}
\centerline{\relabelbox 
\epsfysize 1.5cm
\epsfbox{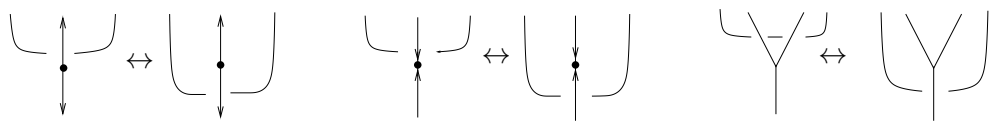}
\relabel{t1}{$\leftrightarrow$}
\relabel{t2}{$\leftrightarrow$}
\relabel{t3}{$\leftrightarrow$}
\endrelabelbox}
\caption{{Fordidden Moves.}}
\label{MoreForbidden}
\end{figure}

{Note that the moves of figure \ref{MoreForbidden} are not allowed.}

\begin{figure}
\centerline{\relabelbox 
\epsfysize 1.5cm
\epsfbox{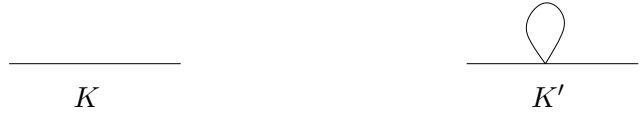}
\relabel{K}{$K$}
\relabel{Kp}{$K'$}
\endrelabelbox}
\caption{Adding a trivial 1-handle to a welded virtual graph. On the left {we display} the original graph.}
\label{addfirst}
\end{figure}

If $K$ is a {welded virtual} graph, then a  {welded virtual} graph $K'$ for which {the tube $T(K')$ of $K'$} is obtained from $T(K)$ by adding a trivial 1-handle is obtained from $K$ by choosing a string of $K$ (in the correct component) and doing the transition shown in figure \ref{addfirst} {(adding a trivial 1-handle to a welded virtual graph)}.

\begin{figure}
\centerline{\relabelbox 
\epsfysize 2.5cm
\epsfbox{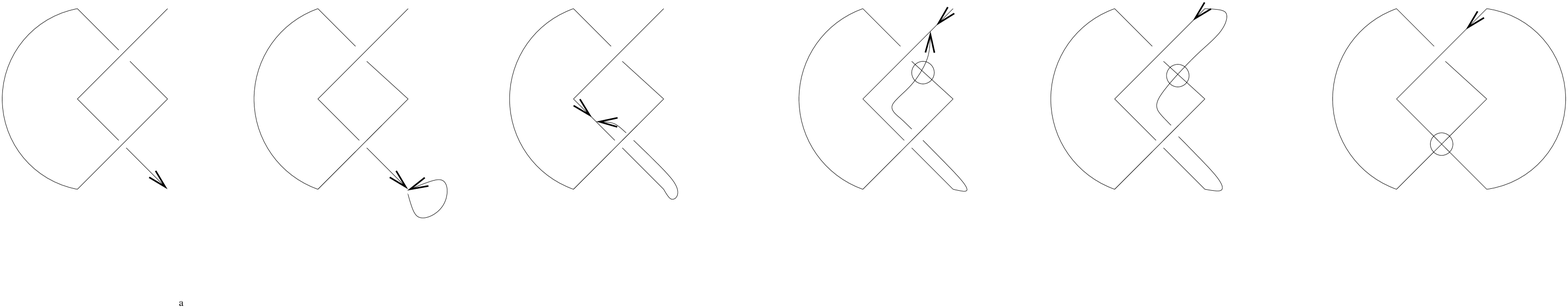}
\endrelabelbox}
\caption{{Adding a trivial 1-handle to the Hopf Arc yields the Virtual Hopf Link.}}
\label{Add}
\end{figure}

For example, consider the Hopf Arc $HA$ defined in Exercise \ref{E1}. Then adding a trivial 1-handle to the {unclosed} component of it yields the Virtual Hopf Link $L${; see} figure \ref{Add}. Note the usage of the moves of figure \ref{WVG}. 

{Let $G_1$ be {a} {welded} virtual graph such that, topologically, $G_1$ is the union of circles $S^1$ and intervals $I=[0,1]$. Suppose that $G_1'$ is obtained from $G_1$ by adding a trivial 1-handle to an $I$-component {of it}. Then we can always use the moves of figure   \ref{WVG} to find a graph $G_2$, equivalent to $G_1'$ as a welded virtual graph, but so that, topologically, $G_2$ is the union of circles $S^1$ and intervals $I$. This was exemplified above for the case of the  Hopf Arc $HA${, and should} be compared with the method indicated in \cite[page 541]{S}.}

It is a good exercise to verify that adding a trivial 1-handle to the Trefoil Arc yields the Shin Satoh's Knot.

\subsubsection{{The} fundamental group of the complement}\label{fgc}
\begin{figure}
\centerline{\relabelbox 
\epsfysize 2.5cm
\epsfbox{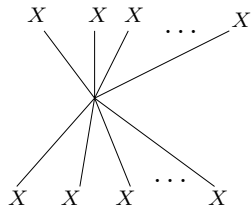}
\relabel{A}{$\s{X}$}
\relabel{B}{$\s{X}$}
\relabel{C}{$\s{X}$}
\relabel{D}{$\s{X}$}
\relabel{E}{$\s{X}$}
\relabel{F}{$\s{X}$}
\relabel{G}{$\s{X}$}
\relabel{H}{$\s{X}$}
\relabel{t}{$\ldots$} 
\relabel{t2}{$\ldots$} 
\endrelabelbox}
\caption{The relations satisfied  {by the} knot group of a welded virtual graph at a vertex.}
\label{ggroup}
\end{figure}

The (combinatorial) fundamental group of a welded virtual graph complement (the knot group) is defined in the same way as the knot group of a virtual knot or arc. {However, we consider the relations} of figure \ref{ggroup} at the vertices of a graph {(the  edges incident to a vertex may carry any orientation).} Note that this is in sharp contrast with the classical fundamental group of graph complements.  In fact, we can easily find examples of welded virtual graphs for which the classical and virtual knot {groups} are different. The $\theta$-graph which appears in figure \ref{wvgraph} is such an example. 

It is not difficult to see that the knot group is an invariant of welded virtual graphs. Moreover the tube map $K \mapsto T(K)$ preserves knot groups.

{Suppose that the graph $K'$ is obtained from $K$ by adding  a trivial 1-handle. We can see that the knot groups of $K$ and $K'$ are isomorphic, thus also that the fundamental groups of the complements of the  tubes $T(K)$ and $T(K')$ in $S^4$ are isomorphic. {This can  easily be proved directly.}}

Given an arc $A$ embedded in the upper semiplane $\{ z \ge 0\}$ of $\R^3$, intersecting the plane $\{z=0\}$ at the end points of $A$, only, there {exist} {two knotted tori} naturally associated to {$A$}. The first one is obtained from {the tube $T(A)$ of $A$} by adding a trivial 1-handle, and a virtual knot $c_1(A)$ representing it can be easily determined from $A$ using the method indicated in \cite{S} and \ref{FMK}. {In the second one,} one simply closes $A$ in the obvious way, obtaining $c_2(A)$, before taking the tube of it. If $A$ is a classical arc, with only one component, then we have that the fundamental groups of the complements of the  knotted surfaces $T(A)$, $T(c_1(A))$ and $T(c_2(A))$ are all isomorphic. {This also happens if we allow $A$ to have more than one component, as long as all {the} other components {are} diffeomorphic to $S^1$.} {However, it is necessary that $A$ be classical.}

{The pairs of {welded virtual knots $\left(c_1(A),c_2(A)\right)$, {one for each classical $1$-component arc $A$,}} provide a family of welded virtual knots with the same knot group. For example if $3_1'$ is the Trefoil Arc, then $c_2(3_1')$ is the {Trefoil Knot $3_1$,} whereas $c_1(3_1')$ is the  Shin Satoh Knot. These two can be proven to be non-equivalent by using the crossed module invariant; see \ref{nonab}. See  also subsections \ref{Fig8} to \ref{Open} for other analogous examples.}

\begin{Problem}
Under which circunstancies are the welded  virtual knots $c_1(A)$ and $c_2(A)$ equivalent? What to say about their tubes in $S^4$.
\end{Problem}

\subsubsection{Crossed module invariants of welded virtual graphs}\label{CA}
{Let $\G=(E,G,\t)$ be an automorphic finite crossed module.}
The invariant $\H_\G$ of welded virtual knots{, or arcs,} extends in a natural way to an invariant of welded virtual graphs {$K$}, by considering: 
$${\H_\G(K)\doteq I_\G(T(K)),}$$ 
{where $I_\G$ is the 4-dimensional invariant defined in subsection \ref{cmiks}.}
As before, $\H_\G(K)$ can be calculated directly from a diagram of $K$.

\begin{figure}
\centerline{\relabelbox 
\epsfysize 3cm
\epsfbox{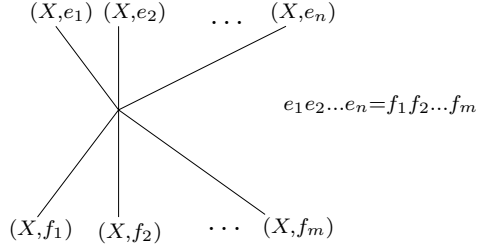}
\relabel{A}{$\s{(X,e_1)}$}
\relabel{B}{$\s{(X,e_2)}$}
\relabel{D}{$\s{(X,e_n)}$}
\relabel{E}{$\s{(X,f_1)}$}
\relabel{F}{$\s{(X,f_2)}$}
\relabel{H}{$\s{(X,f_m)}$}
\relabel{t}{$\ldots$}
\relabel{t2}{$\ldots$}
\relabel{rel}{$\s{e_1e_2\ldots e_n=f_1f_2\ldots 	f_m}$}
\endrelabelbox}
\caption{Reduced {$\G$-colourings} of a welded {virtual} graph diagram {at a vertex}. }
\label{colourG}
\end{figure}

\begin{Definition}
Let $\G=(E,G,\t)$ be a finite automorphic crossed module. 
Let $K$ be an oriented welded virtual graph diagram chosen so that the projection on the second variable is a Morse function in $K$. A reduced $\G$-colouring of  $K$ is given by an assignment of a pair $(X,e)\in G \times E$ to each arc of $G$ minus its set of critical points, crossings and vertices, {satisfying} the conditions already shown for virtual {knot} and arc  diagrams, and the relation {displayed} in figure \ref{colourG}.
\end{Definition}

\begin{figure}			
\centerline{\relabelbox 
\epsfysize 1.5cm
\epsfbox{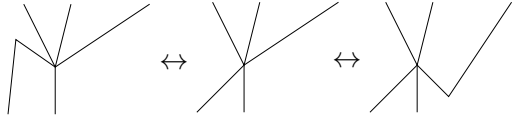}
\relabel{t}{$\leftrightarrow$}
\relabel{t2}{$\leftrightarrow$}
\endrelabelbox}
\caption{One type of Yetter's moves capturing planar isotopy of graph diagrams.}
\label{YetterG}
\end{figure}

We have:

\begin{Theorem}
Let $\G=(E,G,\t)$ be a finite automorphic crossed module. Let also $K$ be  {an oriented welded virtual graph diagram.} The quantity:
\begin{multline}\H_\G(K)=\#\{\textrm{reduced }\G\textrm{  colourings of }K\}\#E^{\#\{\textrm{caps}\}} \#E^{-\#\{\textrm{cups}\}}\\ {\#E^{\#\{\textrm{pointing upward 1-valent vertices of } K\}}}\\
  \prod_{\substack{n\textrm{-valent vertices } v \textrm{ of } K \\{n\ge 2}} } \#E^{1-\#\{\textrm{edges  of } K \textrm{ incident to  } v \textrm{ from above} \}}
\end{multline}
 coincides with $I_\G(T(K))$, and therefore defines an invariant of welded virtual graphs.
\end{Theorem}
\begin{Exercise}
Check directly that $\H_\G$ is a topological invariant {of welded virtual graphs}. Together with the moves of figure \ref{WVG}, as well as the moves defining welded virtual knots and arcs, one still {needs} to check invariance under planar isotopy. Planar isotopies of graph diagrams are captured by Yetter's moves shown in figure \ref{Yetter}, as well as figure \ref{YetterG}; see \cite{Y2} and \cite{FY}.
\end{Exercise}
\begin{Exercise}\label{AAA}
 Check {directly that} $\H_\G$ is invariant under {addition} of trivial $1$-handles, as shown in figure \ref{addfirst}; cf. Theorem \ref{main}.
\end{Exercise}

\subsubsection{The Trefoil Knot is not equivalent to the {Shin-Satoh's Knot}} \label{nonab}

 {We now use the extension of the Crossed Module Invariant to welded virtual graphs to prove} that the Shin Satoh's Knot $S$ is not equivalent to the Trefoil Knot $3_1$. Let $3_1'$ be the Trefoil Arc. Recall that $S$ is obtained from $3_1'$ by adding a trivial 1-handle, in other words $S=c_1(3_1')$; see \ref{fgc}. Therefore, whenever $\G=(E,G,\t)$ is a finite {automorphic} crossed module we have: 
$$\H_\G(3_1')=\H_\G(S).$$  
In particular, from equation {(\ref{exp3}):}
\begin{equation}
\H_\G(S)={\H_\G(3_1')}=\#E\# \left \{\s{X,Y \in G; e \in E} \left | \substack{  XYX=YXY \\ \\-XY \t e+Y \t e=e}  \right. \right \},
\end{equation}
{note that we switched to additive notation.}
Also, from equation (\ref{exp}):
\begin{equation}
\H_\G(3_1)=\# \left \{\s{X,Y \in G; e,f \in E} \left | \substack{ X YX=YXY \\ \\YXY \t f  -XY \t (e+f) +Y \t (e+f)=e}  \right. \right \}.
\end{equation}
{A natural  example of a finite {automorphic} crossed module {$\G=(E,G,\t)$} with $G$ non abelian is constructed by taking $G={\rm GL}_n(\Z_p)$ and $E=(\Z_p)^n$. Here ${\rm GL}_n(\Z_p)$ denotes the group of $n\times n$ matrices in $\Z_p$ with invertible determinant, where $p$ is a {positive integer}. The action of ${\rm GL}_n(\Z_p)$ in $ (\Z_p)^n $ is taken to be the obvious  one. Denote these crossed modules by  $\G{(n,p)}$.}

Computations with {\it Mathematica} prove that {$\H_{\G{(n,p)}}(3_1) \neq \H_{\G{(n,p)}}(S)$} for example for {$p=3,4,5,7$} and $n=2${; see the following table.} 
This proves that the crossed module invariant distinguishes the Trefoil Knot from the Shin Satoh's Knot{, even though they have the same knot group.}

\begin{table}[ht]
\caption{}
\centering
	\begin{tabular}{c c c c c c}
 \hline \hline
knot & $\H_{\G_{2,2}}$ & $\H_{\G_{2,3}}$   & $\H_{\G_{2,4}}$   & $\H_{\G_{2,5}}$  & $\H_{\G_{2,7}}$\\
\hline
$3_1$ & 96 & 4320  & 24576 &132000 &  2272032 \\
$S$ & 96 &  4752 & 27648 &168000 & 2765952  \\
 \end{tabular}
\end{table}

\subsection{Figure of Eight Knot}\label{Fig8}

\begin{figure}
\centerline{\relabelbox 
\epsfysize 8cm
\epsfbox{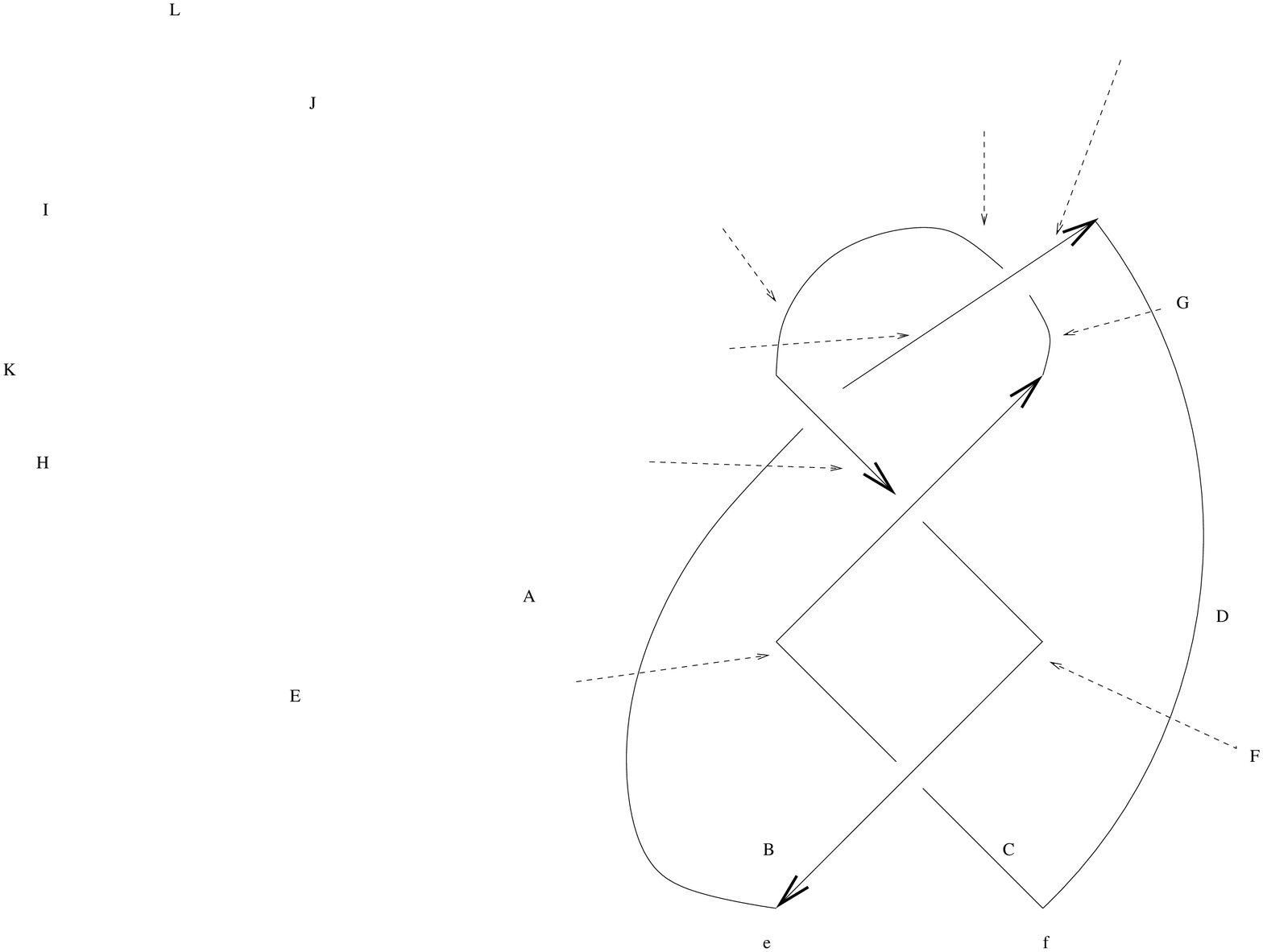}
\relabel{A}{$\ss{\left (X,e^{-1}\right)}$}
\relabel{B}{$\ss{\left (X,e\right )}$}
\relabel{C}{$\ss{\left (Y,f\right )}$}
\relabel{D}{$\ss{\left (Y,f^{-1}\right )}$}
\relabel{E}{$\ss{\left (X^{-1}YX,X^{-1} \t f\right) }$}
\relabel{F}{$\ss{\left (X,X^{-1} \t f^{-1} e f \right )}$}
\relabel{H}{$\ss{\left	 (X^{-1}YXY^{-1}X,X^{-1}Y\t f^{-1}  X^{-1}YX \t (ef)  \right )}$}
\relabel{G}{$\ss{\left (X^{-1}YX,X^{-1}Y \t f X^{-1}YX \t (ef)^{-1}ef \right )}$}
\relabel{K}{$\ss{\left (X^{-1}YXY^{-1}XYX^{-1}Y^{-1}X,X^{-1}YXY^{-1}X \t e^{-1}   \right )}$}
\relabel{I}{$\ss{\left (X^{-1}YXY^{-1}X,X^{-1}Y\t f^{-1}  X^{-1}YX \t (ef)   e^{-1}  X^{-1}YXY^{-1}X \t e 
\right )}$}
\relabel{L}{$\ss{\left (X^{-1}YXY^{-1}XYX^{-1}Y^{-1}X, X^{-1}YXY^{-1}X \t e^{-1}   X^{-1}Y \t f X^{-1}YX \t (ef )^{-1}ef YX^{-1}Y \t f^{-1} Y X^{-1}YX \t (ef)Y\t(ef)^{-1}  \right )}$} 
\relabel{J}{$\ss{\left (YX^{-1}YXY^{-1},YX^{-1}Y \t f Y X^{-1}YX \t (ef)^{-1}Y\t(ef) \right )}$}
\endrelabelbox}
\caption{Calculation of the crossed module invariant of the Figure of Eight Knot $4_1$.}
\label{f8}
\end{figure}

Let $\G=(E,G,\t)$ be a finite automorphic crossed module. Let us calculate the crossed module invariant $\H_\G(4_1)$ of the Figure of Eight Knot $4_1$.  This calculation appears in Figure \ref{f8}.  
This permits us to conclude that{, if $\G=({E,G},\t)$ is an automorphic crossed module,} then:
\begin{equation}\label{previous}
\H_\G(4_1)=\#\left\{\s{X,Y \in G; e,f \in E }
\left | 
\substack {
\s{X^{-1}YXY^{-1}XYX^{-1}Y^{-1}X=Y }\\\\
\s{ X^{-1}YXY^{-1}X \t e^{-1}   X^{-1}Y \t f X^{-1}YX \t (ef )^{-1}ef YX^{-1}Y \t f^{-1} Y X^{-1}YX \t (ef)Y\t(ef)^{-1}=f}\\ \\ 
\s{X^{-1}YXY^{-1}X=YX^{-1}YXY^{-1}}\\ \\
\s{X^{-1}Y\t f^{-1}  X^{-1}YX \t (ef)   e^{-1}  X^{-1}YXY^{-1}X \t e =YX^{-1}Y \t f^{-1} Y X^{-1}YX \t (ef)Y\t(ef)^{-1} }
} \right.
\right \}
\end{equation}
Note that the first pair of equations which appear in the previous formula is  equivalent to the second one. In the case when $G$ is abelian, the previous formula simplifies to (passing to additive notation):
\begin{align*}
\H_\G(4_1)&=\# \left\{X \in G; e,f \in E|(X^2-3X+1)\t e= (-X^2+3X-1) \t f \right\}\\
&=\#E \#\left\{X \in G; e \in E|(X^2-3X+1)\t e=0 \right\},
\end{align*}
as it should, since the Alexander polynomial of the Figure of Eight Knot   is $\Delta(4_1)=X^2-3X+1$; {see \ref{Relation}.}

\begin{figure}
\centerline{\relabelbox 
\epsfysize 3cm
\epsfbox{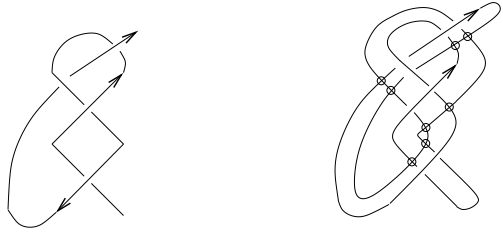}
\endrelabelbox}
\caption{The Figure of Eight Arc $4_1'$ and the welded virtual knot $c_1(4_1')$ obtained from it by adding a trivial 1-handle}
\label{AF8}
\end{figure}

The value of the crossed module invariant for the Figure of Eight Arc $4_1'$ (figure \ref{AF8}), for $G$ not necessarily abelian, can be obtained from equation (\ref{previous})  by making $f=1$, and inserting the relevant normalisation factors. This yields:
$$\H_\G(4_1')=\# E\#\left\{\s{X,Y \in G;{e \in E} }
\left | 
\substack {
\s{X^{-1}YXY^{-1}XYX^{-1}Y^{-1}X=Y }\\\\
\s{ X^{-1}YXY^{-1}X \t e^{-1}   X^{-1}YX \t e ^{-1}e  Y X^{-1}YX \t e Y\t e^{-1}=1}\\ 
} \right.
\right \}.
$$
Consider the welded virtual knot $c_1(4_1')$ obtained from the Figure of Eight Arc $4_1'$ by adding a trivial 1-handle to it{; see \ref{FMK}}. This welded virtual knot appears in figure \ref{AF8}. By using Theorem \ref{addtheo}, it thus follows that for any finite automorphic crossed module $\G$ we have $\H_\G(c_1(4_1'))=\H_\G(4_1')${; see also \ref{CA}.}
Recall  that by the discussion in \ref{fgc},  the knot groups of the welded virtual knots $4_1=c_2(4_1')$ and $c_1(4_1')$ are isomorphic.

Consider the crossed modules $\G_{{(n,p)}}$, where $p$ and $n$ are positive integers, {obtained from ${\rm GL}_n(\Z_p)$ acting on $(\Z_p)^n$}, defined in \ref{nonab}. Computations with {\it Mathematica}  prove that $\H_{\G_{{(n,p)}}}(c_1(4_1')) \neq \H_{\G_{{(n,p)}}}(4_1)$  for $p=3$ or $p=7$; see the following table. This proves that the welded {virtual knots}  {$4_1=c_2(4_1')$} and $c_1(4_1')$ are not equivalent, even though they have the same knot groups.

\begin{table}[ht]
\caption{}
\centering
\begin{tabular}{c c c c c c}
 \hline \hline
knot & $\H_{\G_{2,2}}$ & $\H_{\G_{2,3}}$   & $\H_{\G_{2,4}}$ &  $\H_{\G_{2,5}}$ &  $\H_{\G_{2,7}}$\ \\
\hline
$4_1$ & 48 & 3024  &  15360 & 228000& 1876896\\
$c_1(4_1')$ & 48 &3456   & 15360 & 228000 &2272032\\
 \end{tabular}
\end{table}

\subsection{The Solomon Seal Knot}\label{SSK}

\begin{figure}
\centerline{\relabelbox 
\epsfysize 8cm
\epsfbox{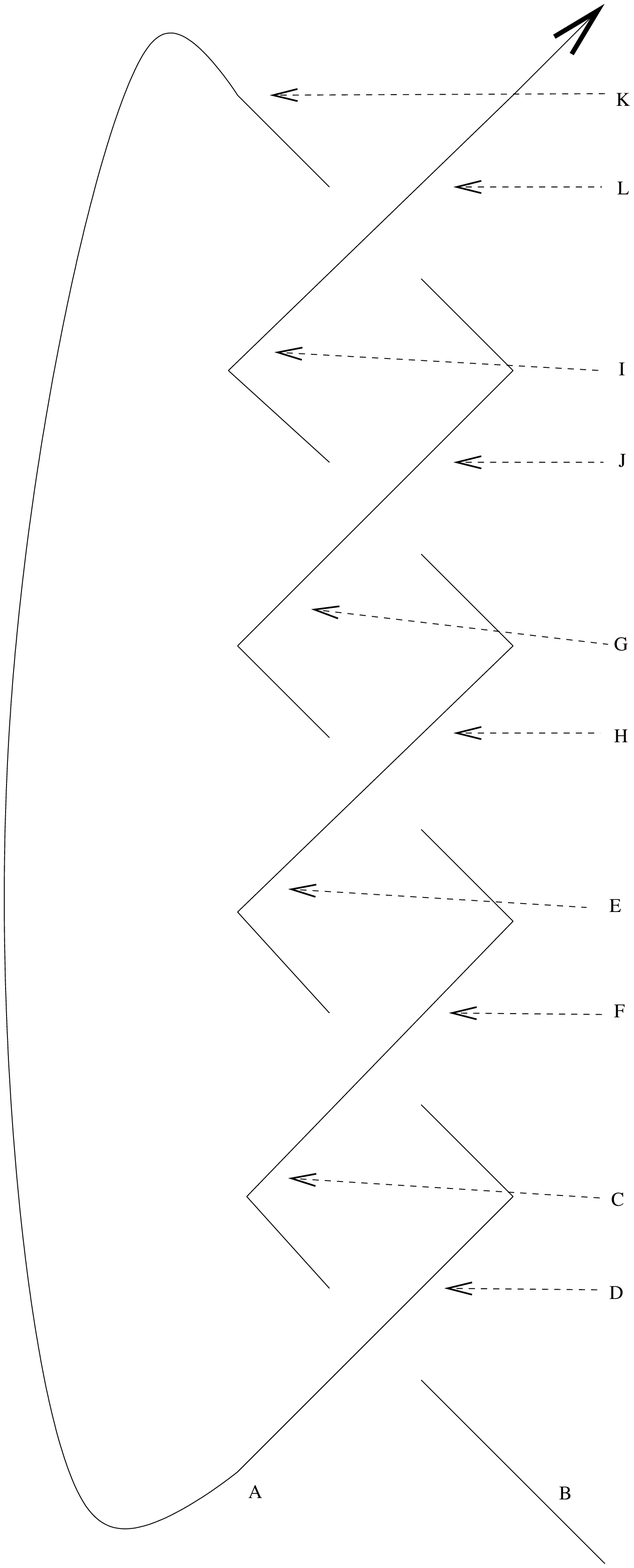}
\relabel{A}{$\s{\left (X,e  \right)}$}
\relabel{B}{$\s{\left (Y,f  \right)}$}
\relabel{C}{$\s{\left (XYX^{-1},X \t f   \right)}$}
\relabel{D}{$\s{\left (X,X\t f^{-1} ef  \right)}$}
\relabel{E}{$\s{\left (XYXY^{-1}X^{-1}, XY \t f^{-1}  XYX^{-1} \t (ef) \right)}$}
\relabel{F}{$\s{\left (XYX^{-1}, XY \t f  XYX^{-1} \t (ef)^{-1}ef \right)}$}
\relabel{G}{$\s{\left (XYXYX^{-1}Y^{-1}X^{-1}, XYX\t f XY \t (ef)^{-1} XYXY^{-1}X^{-1} \t (ef) \right)}$}
\relabel{H}{$\s{\left ( XYXY^{-1}X^{-1}  , XYX\t f^{-1} XY \t (ef) XYXY^{-1}X^{-1} \t (ef)^{-1} ef \right)}$}
\relabel{I}{$\s{\left (XYXYXY^{-1}X^{-1}Y^{-1}X^{-1}, 
XYXY\t f^{-1} 
XYXYX^{-1} \t (ef) 
XY \t (ef)^{-1} 
XYXYX^{-1}Y^{-1}X^{-1} \t (ef)
  \right)}$}
\relabel{J}{$\s{\left (XYXYX^{-1}Y^{-1}X^{-1},  
XYXY\t f 
XYXYX^{-1} \t (ef)^{-1} 
XY \t (ef) 
XYXYX^{-1}Y^{-1}X^{-1} \t (ef)^{-1} ef \right)}$}
\relabel{K}{$\s{\left (YXYXYXY^{-1}X^{-1}Y^{-1}X^{-1}Y^{-1}, YXYXY\t f 
 XYXY\t (ef)^{-1} 
YXY \t (ef)  XY \t (ef)^{-1} Y \t (ef)  \right)}$}
\relabel{L}{$\s{\left (XYXYXY^{-1}X^{-1}Y^{-1}X^{-1}, 
YXYXY\t f^{-1} 
 XYXY\t (ef) 
YXY \t (ef) ^{-1}
XY \t (ef) Y \t (ef)^{-1} ef 
 \right)}$}
\endrelabelbox}
\caption{Calculation of the crossed module invariant of the torus knot $5_1$.  In the top two colourings, we are using the fact that $XYXYX=YXYXY$.}
\label{K5}
\end{figure}

Let $\G=(E,G,\t)$ be an automorphic finite crossed module.
The crossed module invariant of the $(5,2)$-torus knot $5_1$ (the Solomon Seal Knot) is calculated in figure \ref{K5}. This permits us to conclude that:
$$\H_\G(5_1)=\#\left \{\s{X,Y \in G; e,f \in E} \left |\substack{ XYXYX=YXYXY \\\\ YXYXY\t f 
 XYXY\t (ef)^{-1} 
YXY \t (ef)  XY \t (ef)^{-1} Y \t (ef)=e }\right . \right \}.$$
Note that if the crossed module $\G=(E,G,\t)$ is  such that $G$ is abelian, then the previous expression simplifies to: 
\begin{align*}\H_\G(5_1)=\#E\#\left \{X \in G; e \in E  \left | X^4\t e-X^3 \t e +X^2 \t e-X \t e+e=0 \right . \right\},
\end{align*}
as it should, since the Alexander polynomial of the knot $5_1$ is $\Delta(5_1)=X^4-X^3+X^2-X+1$.

\begin{figure}
\centerline{\relabelbox 
\epsfysize 2.5cm
\epsfbox{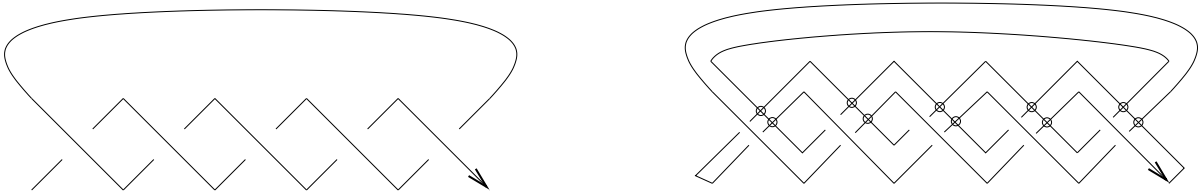}
\endrelabelbox}
\caption{The Solomon Seal arc $5_1'$ and the {welded virtual knot} $c_1(5_1')$ obtained by adding a trivial 1 handle to it.}
\label{ssa}
\end{figure}

The crossed module invariant of the Solomon Seal Arc $5_1'$, and the welded virtual knot $c_1(5_1')$ obtained from it by adding a trivial 1 handle, each presented in figure \ref{ssa}, can be obtained from {this calculation} by making $f=1$, and inserting the remaining {normalisation} factors.  Therefore it follows that:
$$\H_\G(5_1')=\#E\#\left \{\s{X,Y \in G; {e \in E}} \left |\substack{ XYXYX=YXYXY \\  \\
 XYXY\t {e}^{-1} 
YXY \t {e}  XY \t {e}^{-1} Y \t e=e }\right . \right \}.$$

Computations with  {\it Mathematica} show that $\H_{\G_{{(n,p)}}}(c_1(5_1'))  \neq \H_{\G_{{(n,p)}}}(5_1)$  for $n=2$ and $p=5$; see the following table. Therefore the pair $(5_1,c_1(5_1))$ is  a pair of welded virtual knots with the same knot group, but distinguished by their crossed module invariant.

\begin{table}[ht]
\caption{}
\centering
\begin{tabular}{c c c c c c}
 \hline \hline
knot & $\H_{\G_{2,2}}$ & $\H_{\G_{2,3}}$ & $\H_{\G_{2,4}}$   & $\H_{\G_{2,5}}$& $\H_{\G_{2,7}}$\\
\hline
$5_1$ & 24 & 432 &  1536 &168000 &98784 \\
$c_1(5_1')$ & 24 & 432 &1536 &204000 &98784
 \end{tabular}
\end{table}

\subsection{The 2-bridge knot $5_2$ {(Stevedore)}}

\begin{figure}
\centerline{\relabelbox 
\epsfysize 3cm
\epsfbox{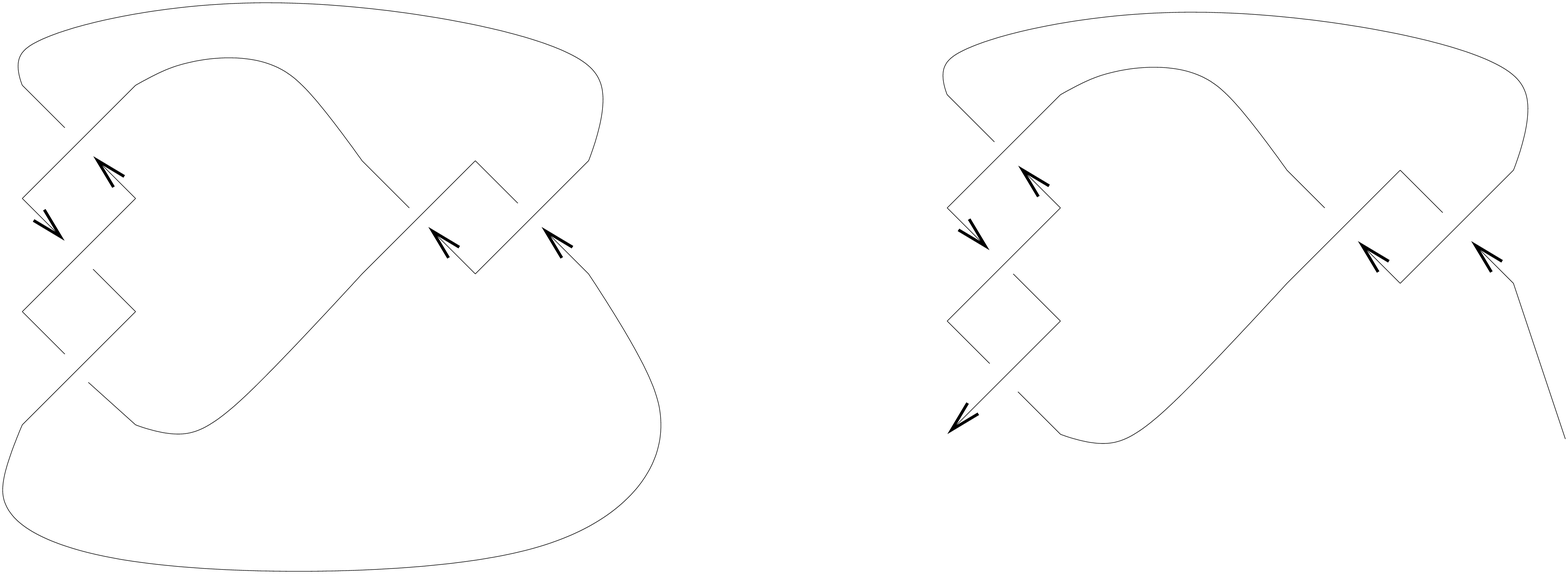}
\endrelabelbox}
\caption{The 2-{bridge} knot $5_2$ and the 2-{bridge} arc $5_2'$.}
\label{Knot}
\end{figure}

\begin{figure}
\centerline{\relabelbox 
\epsfysize 9cm
\epsfbox{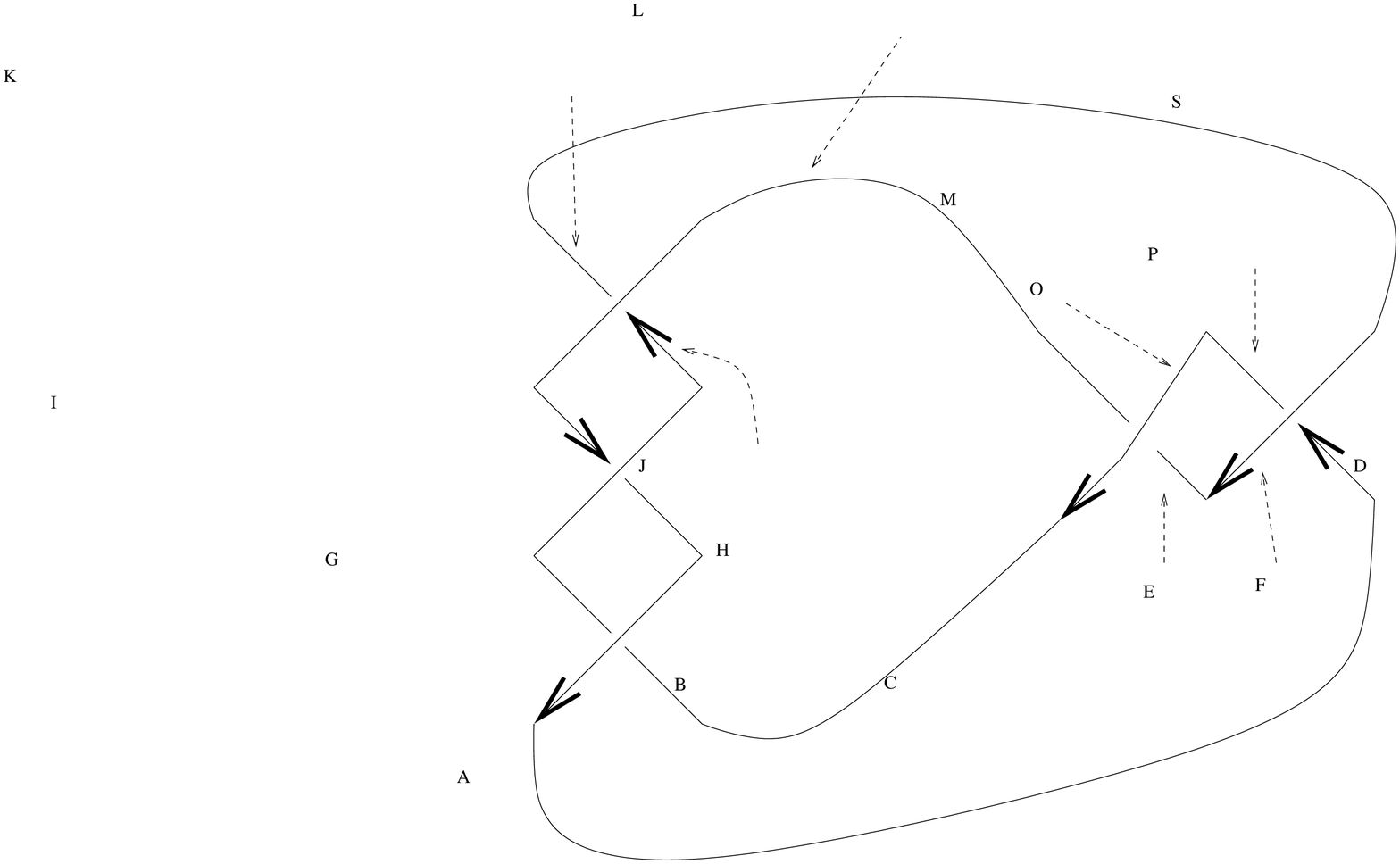}
\relabel{A}{$\ss{(X,f)}$}
\relabel{B}{$\ss{(Y,e)}$}
\relabel{C}{$\ss{(Y,e^{-1})}$}
\relabel{D}{$\ss{(X,f^{-1})}$}
\relabel{E}{$\ss{(Z,g)}$}
\relabel{F}{$\ss{(Z,g^{-1})}$}
\relabel{G}{$\ss{(X^{-1}YX,X^{-1} \t e)}$}
\relabel{H}{$\ss{(X,X^{-1} \t e^{-1} ef)}$}
\relabel{I}{$\ss{(X^{-1}YXY^{-1}X,X^{-1}Y \t e^{-1} X^{-1}YX\t(ef))}$}
\relabel{J}{$\ss{(X^{-1}YX,X^{-1}Y \t e X^{-1}YX\t(ef)^{-1}ef)}$}
\relabel{K}{$\ss{(X^{-1}YX^{-1}YXY^{-1}X,X^{-1}YX^{-1} \t e  X^{-1}Y\t(ef)^{-1} X^{-1}YX^{-1}Y^{-1}X\t (ef) )}$}
\relabel{L}{$\ss{(X^{-1}YXY^{-1}X, X^{-1}YX^{-1} \t e^{-1}  X^{-1}Y\t(ef) X^{-1}YX^{-1}Y^{-1}X\t (ef)^{-1}ef)}$}
\relabel{M}{$\ss{(Y^{-1}ZY,Y^{-1} \t g)}$}
\relabel{O}{$\ss{(Y,e^{-1}gY^{-1} \t g^{-1})}$}
\relabel{P}{$\ss{(Z^{-1}XZ,Z^{-1} \t f^{-1})}$}
\relabel{S}{$\ss{(Z,g^{-1}f^{-1}Z^{-1} \t f)}$}
\endrelabelbox}
\caption{Calculation of the crossed module invariant of the  2-{bridge} knot $5_2$.}
\label{KnotCalc}
\end{figure}

We now consider the 2-{bridge} knot $5_2$ and the 2-bridge arc $5_2'$, depicted in figure \ref{Knot}.
Let us calculate their crossed module invariant. Suppose that $\G=(E,G,\t)$ is a finite automorphic crossed module.  The calculation of $\H_\G(5_2)$ appears in figure \ref{KnotCalc}.
From this it follows that:
\begin{equation}
\H_\G(5_2)=\#\left \{\s{X,Y \in G; {e, f,g \in E}} \left |\substack { YX^{-1}YXY^{-1}XY^{-1}=X^{-1}YX^{-1}YXY^{-1}X\\  \\YX^{-1}YX^{-1} \t e^{-1}  YX^{-1}Y\t(ef) YX^{-1}YX^{-1}Y^{-1}X\t (ef)^{-1}Y\t (ef)=g^{-1}\\ \\ Z=YX^{-1}YXY^{-1}XY^{-1}\\\\ X^{-1}YX^{-1} \t e  X^{-1}Y\t(ef)^{-1} X^{-1}YX^{-1}Y^{-1}X\t (ef) =gf Z^{-1} \t f^{-1} \\\\    Y=Z^{-1}XZ \\\\ Z^{-1} \t f= e^{-1}gY^{-1} \t g^{-1}
  }\right.\right\}.
\end{equation}
Note that the last two equations in the previous formula follow from the remaining. 
When $G$ is abelian,  the previous expression reduces to:
$$\H_\G(5_2)=\# E\# \{X \in G; e \in E| 2X^2 \t e -3e+2 X^{-1}\t e=0\},
$$
as it should since the Alexander polynomial of the $5_2$ knot is $\Delta(5_2)=2X^2-3+2X^{-2}.$  The formula for the crossed module invariant of the arc $5_2'$ is:
\begin{equation}
\H_\G(5_2)=\# E\#\left \{\s{X,Y \in G; {e,g \in E}} \left |\substack { YX^{-1}YXY^{-1}XY^{-1}=X^{-1}YX^{-1}YXY^{-1}X\\  \\YX^{-1}YX^{-1} \t e^{-1}  YX^{-1}Y\t e YX^{-1}YX^{-1}Y^{-1}X\t e^{-1}Y\t e=g^{-1}\\ \\ Z=YX^{-1}YXY^{-1}XY^{-1}\\\\ X^{-1}YX^{-1} \t e  X^{-1}Y\t e^{-1} X^{-1}YX^{-1}Y^{-1}X\t e =g   }\right.\right\}.
\end{equation}

{Below there is a table comparing the value of $\H_{\G_{{(n,p)}}}(5_2)$ and $\H_{\G_{{(n,p)}}}(c_1(5_2'))$, for {$n=2$ and $p=2,3,4,5,7$.} Here as usual, $c_1(5_2')$ is obtained from the welded virtual {arc $5_1'$} by adding a trivial 1-handle to it. In particular it follows that the welded virtual knots $5_2=c_2(5_2')$ and $c_1(5_2')$ are not equivalent, even though they have the same knot groups.
\begin{table}[ht]
\caption{}
\centering
\begin{tabular}{c c c c c c}
 \hline \hline
knot & $\H_{\G_{2,2}}$ & $\H_{\G_{2,3}}$   & $\H_{\G_{2,4}}$ &  $\H_{\G_{2,5}}$ & $\H_{\G_{2,7}}$\ \\
\hline
$5_2$ & 24  & 864 &  1536 &  72000 &987840 \\
$c_1(5_2')$ & 24 & 864  &1536   & 84000& 1481760\\
 \end{tabular}
\end{table}
}

\subsection{The $(n,2)$-torus knot}

Let $n$ be an odd integer. An analogous calculation as in the case of the Trefoil Knot and the Solomon Seal Knot proves that  the crossed module invariant of the $(n,2)$-torus knot $K_n$ has the following expression (in additive notation):
$$\H_\G(K_n)=\# \left \{\s{X,Y \in G; e,f \in E}\left |\substack{ \displaystyle{\prod_{i=1}^n S_i=\prod_{i=1}^n S_{i+1}} \\ \\ \displaystyle{\left (\prod_{i=1}^n S_i\right)\t e - \sum_{k=2}^n {(-1)}^{k}\left( \prod_{i=k}^n  S_i \right) \t (e+f)  -f=0 }}\right .\right\} .$$
Here $S_i=X$ if $i$ is even and $S_i=Y$ if $i$ is odd. On the other hand, the crossed module invariant of the arc $A_n$, obtained from $K_n$ in the obvious way (see {subsection} \ref{SSK} for the case $n=5$) is:
$$\H_\G(A_n)=\# E \# \left \{\s{X,Y \in G; {f}\in E}\left |\substack{ \displaystyle{\prod_{i=1}^n S_i=\prod_{i=1}^n S_{i+1}} \\ \\ \displaystyle{- \sum_{k=2}^n {(-1)}^{k}\left( \prod_{i=k}^n  S_i \right) \t f  -f=0 }}\right .\right\} .$$

In the following table, we compare the value, for each positive {odd} integer {$n\leq 17$,} of  the crossed module invariants $\H_{\G_{(2,3)}}$ and $\H_{\G_{(2,5)}}$ for the pair of welded virtual knots $\left(K_n,c_1(	A_n)\right)$, where $c_1(A_n)$ is obtained from $A_n$ by adding a trivial 1-handle to it. Since the knot groups of $c_1(A_n)$ and of $c_2(A_n)=K_n$ are isomorphic, this gives some more examples of pairs of 1-component welded virtual knots with the same knot group, but distinguished by their crossed module invariant.
\begin{table}[ht]
\caption{}
\centering
\begin{tabular}{c c c c c c c c c}
 \hline \hline
 & $K_3$ & $K_5$   & $K_7$ &  $K_9$ & $K_{11}$ & $K_{13}$  & $K_{15}$  & $K_{17}$\\
\hline\
  $\H_{\G_{(2,3)}}$ & 4320 & 432  & 432 & 4320  & 432 & 432 & 4320 & 432\\
$\H_{\G_{(2,5)}}$ & 132000 & 168000  & 12000 & 132000  & 12000 & 12000 & 288000 & 12000\\

\\\hline \hline
& $c_1(A_{3})$ &$c_1(A_{5})$ &$c_1(A_{7})$ & $c_1(A_{9})$&  $c_1(A_{11})$ &$c_1(A_{13})$ & $c_1(A_{15})$ & $c_1(A_{17})$ \\
\hline
 $\H_{\G_{(2,3)}}$ &  4752  &   432 & 432 &  4752  & 432 & 432 & 4752 &432\\
$\H_{\G_{(2,5)}}$ & 168000 & 204000  & 12000 & 168000  & 12000 & 12000 & 360000 & 12000\\
 \end{tabular}
\end{table}

\subsection{{Final examples}}\label{Open}
{Let $m$ be a positive integer. We can define an  automorphic crossed module $\A_m=(\Z_m,\Z_2,\t)$, {where $\Z_2=\{-1,1,\times\}$,} and  the action of $\Z_2$ on $\Z_m$ is $1 \t a=a$ and $(-1) \t a=-a$, where $a \in \Z_m$. This generalises the crossed module $\A=\A_3$ {defined in subsection \ref{vchl}.}}

\begin{figure}
\centerline{\relabelbox 
\epsfysize 3cm
\epsfbox{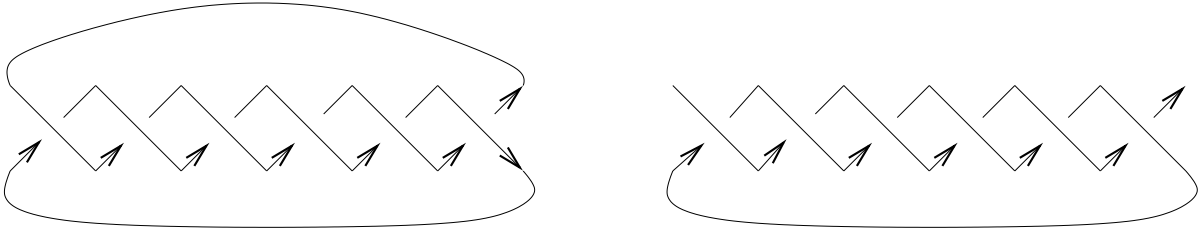}
\endrelabelbox}
\caption{The {link $P$ }and the associated arc $P'$.}
\label{pair}
\end{figure}

Consider the link $P$, as well as {the associated arc $P'$,}  shown in figure \ref{pair}. Let $\G=(E,G,\t)$ be a finite automorphic crossed module with $G$ abelian. An easy calculation shows that:

\begin{multline*}
\H_\G(P)=\# \{X,Y \in G;e,f \in E| \\-Y^{-3}X^{-3} \t f+ Y^{-3}X^{-2}\t (e+f) - Y^{-2}X^{-2} \t (e+f)+ Y^{-2}X^{-1}  \t (e+f) \\-Y^{-1}X^{-1} \t (e+f)+Y^{-1} \t (e+f) = e\}.
\end{multline*}
and
\begin{multline*}
\H_\G({{P'}})={\# E} \# \{X,Y \in G;e \in E| \\ Y^{-3}X^{-2}\t e - Y^{-2}X^{-2} \t e+ Y^{-2}X^{-1}  \t e \\-Y^{-1}X^{-1} \t e+Y^{-1} \t e = e\}.
\end{multline*}
{In the case of the automorphic crossed modules $\A_m$ defined above,  the previous  formulae simplify to (for each positive integer $m$):}

$$\H_{\A_m}(P)=m^2+2m\#\{a \in {\Z_m}| 2a=0\}+m\#\{a \in {\Z_m}|6 a=0\},$$
and: 
$$\H_{\A_m}({P'})=m \left (m+\#\{a \in {\Z_m}| 2a=0\}+m+\#\{a \in {\Z_m}|6 a=0\}\right),$$ 
thus 
 $\H_\A(P)=24$ and $\H_\A(P')=30$. {Here as usual $\A=(\Z_3,\Z_2,\t)$}.

Let $c_1(P')$ be {the}  welded virtual link obtained by adding a trivial 1-handle to the unclosed component of $P'$  (see \ref{fgc}), thus $P'$ and $c_1(P')$ have the same crossed module invariant{; see \ref{CA}}. Hence   $\left(P=c_2(P'), c_1(P')\right)$ is a pair of welded virtual links with the same knot group {(see \ref{fgc})}, but distinguished by their crossed module invariant $\H_\G$, where $\G=(E,G,\t)$ is a finite automorphic crossed module, {which can be chosen so that  $G$ is abelian.} In particular we have $\Alex(P)\cong\Alex(c_1(P'))$, but $\CM(P)\ncong\CM(c_1(P'))$; see {subsection \ref{universal}.}

\begin{figure}
\centerline{\relabelbox 
\epsfysize 2.8cm
\epsfbox{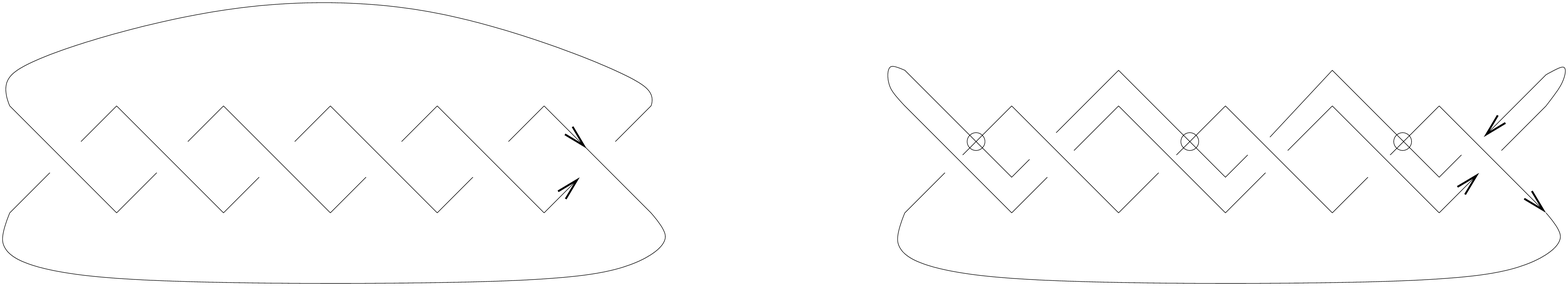}
\endrelabelbox}
\caption{Two {virtual links,} $P=c_2(P')$ and $c_1(P')$, with the same knot group but distinguished by their crossed module invariant.}
\label{EX}
\end{figure}
\begin{Exercise}
Prove directly that $P$ and $c_1(P')$ have the same knot group and are distinguished by their crossed module invariant. 
\end{Exercise} 

\begin{Exercise} The previous example can be generalised. For each positive odd integer $n$, let $P_n$ be the 3-dimensional torus link in $S^3$ with $2n$ crossings, similar to the link $P$ in figure \ref{pair}; {in other words, $P_n$ is the $(2,2n)$-torus link.} Let also $P_n'$ be the associated arc, and let $c_1(P_n')$ be the welded virtual {link obtained} by adding a trivial $1$-handle to the unclosed component of $P_n'$; see figures \ref{pair} and \ref{EX} for the case $n=6$. Prove that for any automorphic crossed module $\G=(E,G,\t)${, with $G$ abelian,} we have that: 
\begin{multline*}\H_\G(P_n)=\#  \{X,Y \in G; e,f \in E  | \\-X^{-n}Y^{-n} \t f +\sum_{k=1}^{n-1} (XY)^{-k} \t \left ( Y^{-1} \t (e+f) -(e+f)\right )+Y^{-1} \t (e+f)=e  \}, \end{multline*} 
and
\begin{equation*}\H_\G(P_n')=\# \left \{X,Y \in G; e \in E  \left | \sum_{k=1}^{n-1} (XY)^{-k} \t \left ( Y^{-1} \t e -e\right )+Y^{-1} \t e=e  \right .\right\}. \end{equation*} 
Thus if $n$ is odd then we have:
$${\H_{\A_m}(P_n)=m^2+2m\#\{a \in \Z_m| 2a=0\} +{m\#\{a \in \Z_m| 2na=0\}}},$$
and 
$${ \H_{\A_m}(P_n')=m \big (m+\#\{a \in \Z_m| 2a=0\} +m+\#\{a \in \Z_m| 2na=0\}\big),}$$
where as usual $\A_m=(\Z_m,\Z_2,\t)$ and $m$ is a positive integer.
{In particular it follows that} 
$${\H_{\A_n}(P_n)=2n^2+2n}$$ and $$\H_{\A_n}(c_1(P_n'))=\H_{\A_n}(P_n')=3n^2+n.$$ 

 {This provides an infinite sequence $\left(P_n,c_1(P_n')\right)$, where $n$ is an odd integer, of {pairs of 2-component welded virtual links with the same knot group, but distinguished by their  crossed module {invariant.}} {This sequence} includes not only {the} previous example, but also  the case of the Hopf Link and the Virtual Hopf Link in subsection \ref{vchl}.}
\end{Exercise}
{Note that, taking tubes, {the previous} example gives an infinite set of pairs of non-isotopic embeddings of a disjoint union of {two} tori $S^1 \times S^1$ into $S^4$ with the same fundamental group of the complement, but distinguished by their Crossed Module Invariant $I_\G$ of \cite{FM1}.}

\begin{figure} 
\centerline{\relabelbox 
\epsfysize 2cm
\epsfbox{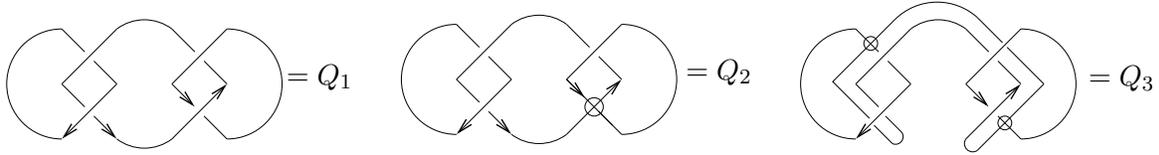}
\relabel{Q}{$=Q_1$}
\relabel{R}{$=Q_2$}
\relabel{S}{$=Q_3$}
\endrelabelbox}
\caption{Three Virtual Links  with the same knot group but distinguished by their crossed module invariant.}
\label{EX3}
\end{figure}

Another interesting example is provided by the virtual {links} $Q_1,Q_2$ and $Q_3$ shown in figure \ref{EX3}. The knot groups of $Q_1,Q_2$ and $Q_3$ are all  isomorphic to {$\{X,Y,Z\colon XY=YX, ZY=YZ\}$.}

Let $\G=(E,G,\t)$ be an automorphic crossed module. A simple calculation shows that:
$$\H_\G(Q_1)=\#\left \{X,Y,Z \in G; e,f,g \in E\left | \substack{-Y^{-1}X^{-1} \t f + Y^{-1} \t (e+f)=e\\- Z^{-1} Y \t g+ Z^{-1} \t (-f+g)=-f} \right .\right \},$$
$$\H_\G(Q_2)=\#\left \{X,Y,Z \in G; e,f,g \in E \left | \substack{-Y^{-1}X^{-1} \t f + Y^{-1} \t (e+f)=e\\Z^{-1} \t f=f} \right.\right \},$$
and:
$$\H_\G(Q_3)=\#\left \{X,Y,Z \in G; e,f,g \in E \left | \substack{ Y^{-1} \t e=e\\- Z^{-1} Y \t g+ Z^{-1} \t g=0} \right . \right \}.$$
Therefore the crossed module invariant $\H_\A$, where as usual $\A=(\Z_3,\Z_2,\t)${, separates} these $Q_1,Q_2$ and $Q_3$. 

\section*{Acknowledgements}

{JFM {was} financed by  Funda\c{c}\~{a}o para a Ci\^{e}ncia e Tecnologia
(Portugal), post-doctoral grant number SFRH/BPD/17552/2004, part of the
research project POCI/MAT/60352/2004 (``Quantum Topology''), also financed by FCT, cofinanced by the European Community fund FEDER.  LK thanks the National Science Foundation for support
under NSF Grant DMS-0245588.}

\end{document}